\documentclass[12pt,a4paper]{report}

\usepackage{amsmath}
\usepackage{amsthm}
\usepackage{amssymb}
\usepackage{amsfonts}
\usepackage{titleref}

\def\ncomm{\newcommand}
\def\rncomm{\renewcommand}

\def\nthm{\newtheorem}

\def\ni{\noindent}

\def\RL{\em{\mbf{\Rightarrow\!\Leftarrow}\tx{ contradiction}}}

\ncomm{\doit}[2]{\begin{#1} #2 \end{#1}}

\ncomm{\MAX}[2][{}]{\max_{#1}\!\set{#2}}
\ncomm{\MIN}[2][{}]{\min_{#1}\!\set{#2}}

\ncomm{\CUP}[1][{}]{\bigcup_{#1}}
\ncomm{\CAP}[1][{}]{\bigcap_{#1}}

\ncomm{\COP}[1][{}]{\coprod_{#1}}

\ncomm{\mbold}[1]{\boldsymbol{#1}}

\theoremstyle{plain}
\nthm{thm}{Theorem}[section]

\newtheorem{prop}[thm]{Proposition}
\newtheorem{corr}[thm]{Corollary}
\nthm{coro}[thm]{Corollary}

\theoremstyle{definition}
\newtheorem*{notation}{Notation}
\newtheorem{conj}{Conjecture}

\newtheorem*{acknowledgement*}{Acknowledgement} 

\theoremstyle{remark}

\ncomm{\addrem}[1]{\doit{rem}{#1}}
\ncomm{\addREM}[1]{\doit{rem}{#1}}

\def\numTHMsec{\numberwithin{thm}{section}}

\ncomm{\eqarab}[1][]{\renewcommand{\theequation}{#1\arabic{equation}}}
\ncomm{\eqroman}[1][]{\renewcommand{\theequation}{#1\roman{equation}}}
\ncomm{\eqAlph}[1][]{\renewcommand{\theequation}{#1\Alph{equation}}}
\ncomm{\eqalph}[1][]{\renewcommand{\theequation}{#1\alph{equation}}}

\ncomm{\boxedeqn}[1]{%
  \[\fbox{%
      \addtolength{\linewidth}{-2\fboxsep}%
      \addtolength{\linewidth}{-2\fboxrule}%
      \begin{minipage}{\linewidth}%
      \begin{equation}#1\end{equation}%
      \end{minipage}%
    }\]%
}

\newcommand{\eval}[2][\right]{\relax
  \ifx#1\right\relax \left.\fi#2#1\rvert}

\AtBeginDocument{\numTHMsec}
\ncomm{\dosp}[1]{\begin{split} #1 \end{split}}              

\ncomm{\doeq}[1]{\doit{equation}{#1}}       
\ncomm{\doEQ}[2]{\doit{equation}{\l{#1}#2}}
\ncomm{\doEQsplit}[2]{\doit{equation}{\doSplit{\l{#1}#2}}}
\ncomm{\doeqnn}[1]{\doit{equation*}{#1}}       
\ncomm{\doeqm}[1]{\doit{math}{#1}}

\ncomm{\doeqnnQED}[1]{\[#1\QED\]}       

\ncomm{\mnote}[1]{\marginpar{#1}}

\ncomm{\doeqsA}[2]{\subeq{#1 \doAlign{#2}}}

\ncomm{\subeq}[1]{#1}

\ncomm{\geBY}[1]{\byeq{#1}{\ge}}
\ncomm{\leBY}[1]{\byeq{#1}{\le}}
\ncomm{\gnBY}[1]{\byeq{#1}{>}}
\ncomm{\lnBY}[1]{\byeq{#1}{<}}
\ncomm{\neBY}[1]{\byeq{#1}{\ne}}

\ncomm{\eqBY}[1]{\byeq{#1}{=}}
\ncomm{\ggBY}[1]{\byeq{#1}{\gg}}
\ncomm{\llBY}[1]{\byeq{#1}{\ll}}
\ncomm{\subseteqBY}[1]{\byeq{#1}{\subseteq}}
\ncomm{\nsubseteqBY}[1]{\byeq{#1}{\nsubseteq}}
\ncomm{\notsubseteqBY}[1]{\byeq{#1}{\notsubseteq}}
\ncomm{\RightarrowBY}[1]{\byeq{#1}{\Longrightarrow}}
\ncomm{\LeftarrowBY}[1]{\byeq{#1}{\Longleftarrow}}
\def\RightarrowL{\Longrightarrow}

\ncomm{\rightarrowUP}[1]{\UP{#1}{\rightarrow}}
\ncomm{\mapstoUP}[1]{\UP{#1}{\mapsto}}

\ncomm{\dostam}[1]{}

\ncomm{\set}[1]{\left\{#1\right\}}
\ncomm{\sset}[1]{\left|\left\{#1\right\}\right|}
\ncomm{\SER}[3][1]{{#3}_{#1},\ldots,{#3}_{#2}}
\ncomm{\SET}[3][1]{\set{\SER[#1]{#2}{#3}}}
\ncomm{\VEC}[3][1]{\lr{\SER[#1]{#2}{#3}}} 

\ncomm{\dotop}[2][]{#1\atop #2}
\ncomm{\stack}[1]{\substack{#1}}
\ncomm{\stackl}[1]{\begin{subarray}{l}#1\end{subarray}}
      
\ncomm{\UP}[2][]{\stackrel{#1}{#2}}
\ncomm{\os}[2][]{\overset{#1}{#2}}
\ncomm{\us}[2][]{\underset{#1}{#2}}

\ncomm{\byeq}[2]{\stackrel{\re{#1}}{#2}}
\ncomm{\byr}[2]{\stackrel{\r{#1}}{#2}}
\ncomm{\byre}[2]{\stackrel{\re{#1}}{#2}}

\def\em{\ensuremath}

\ncomm{\xydiagram}[8]{
\xymatrix{ 
    {#1} \ar[r]^{#2} \ar[d]_{#4} & {#3} \ar[d]^{#5} \\
    {#6} \ar[r]_{#7}       & {#8}
}}

\ncomm{\tn}[1]{\textnormal{#1}}
\ncomm{\tb}[1]{\textbf{#1}}
\ncomm{\tu}[1]{\underline{#1}} 
\ncomm{\te}[1]{\emph{#1}}
\ncomm{\ti}[1]{\textit{#1}}
\ncomm{\tx}[1]{\text{#1}}
\ncomm{\txx}[1]{\texttt{#1}}
\ncomm{\txby}[1]{\tx{(by \re{#1})}}

\ncomm{\tin}[1]{{\tiny #1}}

\ncomm{\tl}[1]{{\large #1}}
\ncomm{\tll}[1]{{\Large #1}}
\ncomm{\tL}[1]{{\LARGE #1}}
\rncomm{\th}[1]{{\huge #1}}
\ncomm{\tH}[1]{{\HUGE #1}}

\def\mbf{\mathbf}	

\ncomm{\tup}[1]{\textup{#1}}

\ncomm{\mT}[1]{{\textstyle #1}}
\ncomm{\mD}[1]{{\displaystyle #1}}
\ncomm{\mS}[1]{{\scriptstyle #1}}
\ncomm{\mSS}[1]{{\scriptscriptstyle #1}}
\ncomm{\mBF}[1]{\mbox{\boldmath$#1$}}

\DeclareMathOperator{\PSL}{PSL}

\DeclareMathOperator{\Prob}{Prob}

\DeclareMathOperator{\diam}{diam}
\DeclareMathOperator{\End}{End}
\DeclareMathOperator{\Ker}{Ker}
\DeclareMathOperator{\Fix}{Fix}
\DeclareMathOperator{\SL}{SL}
\DeclareMathOperator{\GL}{GL}
\DeclareMathOperator{\M}{M}
\DeclareMathOperator{\Tr}{Tr}
\DeclareMathOperator{\tr}{tr}
\DeclareMathOperator{\Prod}{Prod}
\DeclareMathOperator{\Diag}{Diag}

\DeclareMathOperator{\Id}{Id}
\DeclareMathOperator{\mult}{mult}

\def\Stab{{\operatorname{Stab}}}
\def\rank{{\operatorname{rank}}}

\def\Char{{\operatorname{char}}}
\ncomm{\dochar}[1]{\Char(#1)}
\ncomm{\doChar}[1]{\op{Char}(#1)}

\def\supp{{\operatorname{supp}}}
\def\dim{{\operatorname{dim}}}
\def\Cay{{\operatorname{Cay}}}

\ncomm{\doIm}[1]{\op{Im}(#1)}
\def\opIm{\op{Im}}

\def\Span{{\hbox{\rm span}}}
\ncomm{\dospan}[1]{\Span(#1)}
\ncomm{\doSpan}[1]{\op{Span}(#1)}

\ncomm{\doOP}[1]{\defOP{#1}}
\ncomm{\doOPx}[1]{\defOPx{#1}}

\ncomm{\myOP}[1]{\operatorname{#1}}
\ncomm{\myOPx}[1]{\operatorname*{#1}}
\ncomm{\op}[1]{\operatorname{#1}}

\def\opC{\op{C}}
\def\opN{\op{N}}
\def\opE{\op{E}}

\ncomm{\mb}[1]{\mbox{#1}}
\def\mq{\quad}
\ncomm{\mqbq}[1]{\mq\mb{#1}\mq}
\ncomm{\mbq}[1]{\mb{#1}\mq}
\ncomm{\mqb}[1]{\mq\mb{#1}}

\ncomm{\vs}[1][5]{\vskip #1mm}
\ncomm{\vsk}[1][5]{\vskip #1mm}

\ncomm{\doproof}[1]{\doit{proof}{#1}}
\def\QED{\qedhere}

\ncomm{\dothm}[1]{\begin{thm} #1 \end{thm}}
\ncomm{\docorr}[1]{\begin{corr} #1 \end{corr}}
\ncomm{\docoro}[1]{\begin{corr} #1 \end{corr}}
\ncomm{\dofact}[1]{\doit{fact}{#1}}
\ncomm{\dosfact}[1]{\doit{sfact}{#1}}
\ncomm{\dodef}[1]{\doit{defn}{#1}}

\ncomm{\doquest}[1]{\doit{quest}{#1}}
\ncomm{\doprob}[1]{\begin{pro} #1 \end{pro}}
\ncomm{\doconj}[1]{\begin{conj} #1 \end{conj}}

\ncomm{\doprop}[1]{\begin{prop} #1 \end{prop}}
\ncomm{\dolemm}[1]{\doit{lemma_}{#1}}
\ncomm{\dolem}[1]{\doit{lemma_}{#1}}

\ncomm{\donota}[1]{\begin{notation} #1 \end{notation}}

\def\thethmarab{
\setcounter{equation}{0}
\renewcommand{\theequation}{\thethm.\arabic{equation}}
}

\def\thethmaleph{
\setcounter{equation}{0}
\renewcommand{\theequation}{\thethm\alph{equation}}
}

\ncomm{\doTHM}[1]{\thethmaleph \dothm{#1}}
\ncomm{\doCORO}[1]{\thethmaleph \docoro{#1}}
\ncomm{\doPROP}[1]{\thethmaleph \doprop{#1}}
\ncomm{\doLEM}[1]{\thethmaleph \dolem{#1}}
\ncomm{\doFACT}[1]{\thethmaleph \dofact{#1}}
\ncomm{\doSFACT}[1]{\thethmaleph \dosfact{#1}}
\ncomm{\doDEF}[1]{\thethmaleph \dodef{#1}}

\ncomm{\doPROOF}[1]{\thethmarab\doproof{#1}}
\ncomm{\doPROOFof}[2]{\thethmarab\doproofof{#1}{#2}}

\ncomm{\doenum}[1]{\begin{enumerate} #1 \end{enumerate}}
\ncomm{\doitem}[1]{\doit{itemize}{#1}}

\ncomm{\nnum}{\nonumber}
\ncomm{\ntag}{\notag}

\ncomm{\nt}{\notag}
\ncomm{\nn}{\nonumber}

\ncomm{\pageB}{\pagebreak}
\ncomm{\lineB}{\linebreak}

\ncomm{\TBD}[1]{}
\ncomm{\old}[1]{}

\ncomm{\donone}[1]{#1}

\ncomm{\lr}[1]{\left(#1\right)}
\ncomm{\lrv}[1]{\left|#1\right|}
\ncomm{\lrn}[1]{\left\|#1\right\|}
\ncomm{\lrb}[1]{\left[#1\right]}
\ncomm{\lrs}[1]{\left\{\!{#1}\!\right\}}
\ncomm{\lrp}[1]{\left(\!{#1}\!\right)}
\ncomm{\lrg}[1]{\left<\!{#1}\!\right>}
\ncomm{\lre}[1]{\left.{#1}\right.}
\ncomm{\lrf}[1]{\lfloor{#1}\rfloor}

\ncomm{\domatrix}[1]{\doit{matrix}{#1}}
\ncomm{\dosmatrix}[1]{\lr{\doit{smallmatrix}{#1}}}

\ncomm{\dobmatrix}[1]{\dopmatrix{#1}}
\ncomm{\dopmatrix}[1]{\doit{pmatrix}{#1}}
\ncomm{\dobdiag}[1]{\dopmatrix{{#1}&0\\0& #1^{-1}}}
\ncomm{\dosbdiag}[1]{\dosmatrix{{#1}&0\\0& #1^{-1}}}

\ncomm{\docases}[1]{\doit{cases}{#1}}
\ncomm{\dodarrow}[2]{\left. #1 \right \downarrow_{#2}}
\ncomm{\domid}[2]{\left. #1 \right|_{#2}}

\ncomm{\dosideset}[3][]{\sideset{#1}{#2}#3}
\ncomm{\dolim}[2]{\limits_{#1}^{#2}}

\ncomm{\xright}[2][]{\xrightarrow[#1]{#2}}
\ncomm{\xleft}[2][]{\xleftarrow[#1]{#2}}

\ncomm{\tfr}[2][1]{\tfrac{#1}{#2}}
\ncomm{\fr}[2][1]{\frac{#1}{#2}}
\ncomm{\dfr}[2][1]{#1\!/#2}

\ncomm{\sqr}[1]{\sqrt{\!\!#1}}

\def\ST{such that}
\def\WRT{with respect to}
\def\RESP{respectively}

\def\wLOG{without loss of generality}
\def\WLOG{Without loss of generality}

\def\tFHS{the following holds}

\def\TWH{Then we have}

\def\TE{There exist}

\def\TES{There exists}

\ncomm{\stam}[1]{}
\ncomm{\none}[1]{#1}
\ncomm{\gen}[1]{\langle #1\rangle}

\def\good{\em{\surd}}

\def\l{\label}

\ncomm{\re}[1]{\eqref{#1}}

\def\r{\ref}
\def\Sr{\S\ref}

\ncomm{\cf}[1]{(cf. \cite{#1})}
\ncomm{\cfr}[2][XXX]{(cf. \cite[#1]{#2})}
\ncomm{\ci}[1]{~\cite{#1}}
\ncomm{\CITE}[3][]{\tn{(#1\cite[#2]{#3})}} 

\ncomm{\doenv}[2]{\noindent \tb{{#1}:}{#2}}

\ncomm{\doearray}[1]{\doit{eqnarray*}{#1}}
\ncomm{\doeArray}[1]{\doit{eqnarray}{#1}}

\ncomm{\eqL}[1]{\lefteqn{#1}}

\ncomm{\dm}[1]{\doit{displaymath}{#1}}

\ncomm{\marray}[1]{\doit{math}{\doit{array}{#1}}}

\ncomm{\dmarray}[1]{\doit{displaymath}{\doit{array}{#1}}}

\ncomm{\doalign}[1]{\doit{align*}{#1}}
\ncomm{\doAlign}[1]{\doit{align}{#1}}

\ncomm{\doAlignSplit}[2]{\doAlign{\l{#1}\doSplit{#2}}}

\ncomm{\doalignL}[1]{\doit{flalign*}{#1}}
\ncomm{\doAlignL}[1]{\doit{flalign}{#1}}

\ncomm{\doAlignSplitL}[2]{\doAlignL{\l{#1}\doSplit{#2}}}

\def\addtx{\intertext}
\def\addTX{\intertext}

\def\qu{\quad}

\ncomm{\tqq}[1][]{\tx{#1}\qquad}
\ncomm{\tq}[1][]{\tx{#1}\quad}
\ncomm{\qqt}[1][]{\qquad\tx{#1}}
\ncomm{\qtq}[1][]{\quad\tx{#1}\quad}

\ncomm{\qmq}[1]{\;\qu{#1}\qu\;}
\def\eqB{\qmq{=}}
\def\leB{\qmq{\le}}
\def\geB{\qmq{\ge}}
\def\lnB{\qmq{<}}
\def\gnB{\qmq{>}}
\def\llB{\qmq{\ll}}
\def\ggB{\qmq{\gg}}
\def\RightarrowB{\qmq{\Rightarrow}}

\def\RightarrowLB{\qmq{\RightarrowL}}

\def\iffB{\qmq{\iff}}
\def\subseteqB{\qmq{\subseteq}}

\ncomm{\doalignAT}[2][2]{\doit{alignat*}{{#1}#2}}
\ncomm{\doAlignAT}[2][2]{\doit{alignat}{{#1}#2}}

\ncomm{\doSplit}[1]{\doit{split}{#1}}
\ncomm{\eqSplit}[2]{\doeq{\l{#1}\doSplit{#2}}}
\ncomm{\eqsplit}[2]{\doeqnn{\doSplit{#2}}}
\ncomm{\doGather}[1]{\doit{gather}{#1}}

\ncomm{\vsp}[1][1]{\vspace{#1cm}}

\rncomm{\sec}[1]{\chapter{#1}}
\ncomm{\ssec}[1]{\section{#1}}	
\ncomm{\sssec}[1]{\subsection*{#1}} 

\def\BTA{By the assumption}

\def\IP{In particular}

\def\FTOH{On the one hand}

\def\OTOH{On the other hand}

\def\IOW{In other words}

\def\sWAD{so we are done}
\def\SWAD{So we are done}

\def\sWADW{\sWAD\ with }
\def\sWADB{\sWAD\ by }
\def\SWADW{\SWAD\ with }

\def\tWAD{then we are done}
\def\TF{Therefore}

\def\OW{Otherwise}
\def\IAC{In any case}

\def\eps{\em{\varepsilon}}
\ncomm{\epss}[1][]{\varepsilon\!^{#1}}
\ncomm{\Oeps}[1][]{O(\epss[#1]\,\!)}

\def\Oepss{\Oeps[2]}
\def\Om{\Omega}
\def\om{\omega}
\def\Omeps{\Om(\eps)}
\def\Omepss{\Om(\eps^2)}

\ncomm{\NMIDu}[2][]{#1\!\!\nmid\!\!_{\UP{#2}}}
\ncomm{\MIDu}[2][]{#1\!\!\mid\!\!_{\UP{#2}}}

\ncomm{\NMID}[2][]{#1\!\!\nmid\!_{\!\UP{#2}}}
\ncomm{\MID}[2][]{#1\!\!\mid\!_{\!\UP{#2}}}

\def\s{{ }}

\def\WLOG{Without loss of generality}
\def\wLOG{without loss of generality}

\def\char{charecteristic}

\def\mbb{\mathbb}

\def\R{\mathbb{R}}
\def\N{\mathbb{N}}
\def\C{\mathbb{C}}
\def\P{\mathbb{P}}
\def\Z{\mathbb{Z}}
\def\F{\mathbb{F}}
\def\E{\mathbb{E}}
\def\K{\mathbb{K}}

\def\G{\em{{\mathcal{G}}}}

\def\cP{\em{{\mathcal{P}}}}

\def\SS{S\cup S^{-1}}

\ncomm{\seq}[3][1]{#3_{#1},\ldots,\,#3_{#2}}

\def\ddot{\!\cdot\!}

\ncomm{\n}[1]{\em{\|#1\|}}

\ncomm{\Lo}[1]{\em{\n{#1}_{1}}}
\ncomm{\Lt}[1]{\em{{\n{#1}}_2}}
\ncomm{\Lp}[1]{\em{\n{#1}_p}}
\ncomm{\Lq}[1]{\em{\n{#1}_q}}
\ncomm{\Lr}[1]{\em{\n{#1}_r}}  
\ncomm{\Li}[1]{\em{{\n{#1}_{\infty}}}}
\ncomm{\Lm}[1]{\em{{\n{#1}_{\infty}}}}

\ncomm{\No}[1]{\em{\n{#1}_1}}
\ncomm{\Nt}[1]{\em{{\n{#1}}_2}}
\ncomm{\Np}[1]{\em{\n{#1}_p}}
\ncomm{\Nq}[1]{\em{\n{#1}_q}}
\ncomm{\Nr}[1]{\em{\n{#1}_r}}  
\ncomm{\Ni}[1]{\em{{\n{#1}_{\infty}}}}
\ncomm{\Nm}[1]{\em{{\n{#1}_{\infty}}}}

\ncomm{\Ltt}[1]{\em{{\n{#1}}_2^2}}

\def\ol{\overline}

\def\over{/\!\!}

\def\C{\mathbb{C}}

\def\lam{\lambda}
\def\Lam{\Lambda}
\def\del{\delta}

\def\OB{orthonormal basis}

\def\pmI{\em{\set{\pm I}}}
\def\pmO{\em{\set{\pm 1}}}

\def\SM{\!\setminus\!}
\ncomm{\BS}[1]{_{{#1}\backslash}}
\ncomm{\FS}[1]{_{\slash #1}}

\ncomm{\doframe}[2]{
	\doit{frame}{
		\frametitle{#1}
		#2
	}
}

\ncomm{\doblock}[2]{
	\doit{block}{{#1}{#2}}
}

\ncomm{\hide}[1]{}
\ncomm{\showit}[1]{#1}

\def\Plunnecke{Pl\"unnecke}

\begin{document}

\vsp[6]

\doit{center}{

\th{\tb{Expansion properties\\ of finite simple groups}}

\vsp[5]

Thesis submitted for the degree of \\
``Doctor of Philosophy''\\ 
by \\
Oren Dinai

\vsp[3]

Submitted to the Senate of the Hebrew University\\ 
September 2009
}

\newpage
\noindent 
This work was carried out under the supervision of \\  
Prof. Alex Lubotzky

\newpage

\ni
To my mother, who taught me that talent is useless without the ability to carry it out.
\vsk[2]

\ni
To my father, who showed me that imagination has no boundaries, and whom I deeply miss.


\newpage
{\ni \tl{\tb{Acknowledgements}}}

\vsk[2]

Most importantly, I am grateful for working under the supervision of the great mathematical riddle maker, Alex Lubotzky.
Other than his broad mathematical perspective, I had the privilege to get to know a real Mensch\footnote{a good man in Yiddish}.
 
I am thankful to Avi Wigderson and Udi Hrushovski for encouraging me in this work and for many helpful discussions. 

I thank Elon Lindenstrauss for the fruitful conversations and the ongoing encouragement.
These conversations enriched my mathematical thinking.

I wish to thank also Michael Larsen who introduced me to the striking application of the Invariant Theory of section \r{invariant theory} to my work. 

I am grateful to Harald Helfgott for his guidance through his remarkable work and for many illuminating discussions.

I wish to thank also Peter Sarnak, Laci Babai and Alex Gamburd for many useful conversations.

Last but not least, I would like to express my profound thanks to Inna Korchagina 
for many valuable talks and for good advise.

\newpage
\sec{Abstract}

\ssec{Diameter and Growth of Cayley graphs}

A family of finite groups $\{G_n\}_{n\in\N}$ is said to have \ti{poly-logarithmic diameter} 
if for some absolute constants $C,d>0$, 
for every $G_n$ and every subset $S_n\subseteq G_n$ generating $G_n$, 
we have 
$$\diam(\Cay(G_n,S_n))\le C \log^d(|G_n|),$$
where $\diam(\Cay(G,S))$ is the diameter of the Cayley graph of $G$ with respect to $S$.

\vsk[1] 
A well know conjecture of Babai \cite{bs2} asserts that all the non-abelian finite simple groups have poly-logarithmic diameter.
In this work we investigate the family of groups $\SL_2$ (and $\PSL_2$) over finite fields, 
and we prove the conjecture for this family of groups.

In fact, we investigate a stronger Growth property that would imply in particular the poly-logarithmic diameter bounds.
By this, we extend the techniques that were developed by Helfgott \cite{helf} who dealt with the family of groups $\SL_2$ (and $\PSL_2$) over finite fields of prime order.

\ssec{The main results}			

Our main result asserts that the family 
$$\set{\SL_2(\F_{p^n}): p\tx{ prime}; n\in \N}$$
has poly-log diameter. 
Note that this result holds uniformly for all finite fields regardless of their \char.
This result holds also for the family $\PSL_2$ over finite fields.

\vsk[6]
By using results from Additive Combinatorics, we proved the following stronger Growth property:

\TES\ $\eps>0$ \ST\ \tFHS\ for any finite field $\F_q$. 
Let $G$ be the group $\SL_2(\F_q)$ (or $\PSL_2(\F_q)$) and let $A$ be a generating set of $G$.
Then we have, 
$$|A \ddot A \ddot A|\ge \MIN{|A|^{1+\eps},|G|}.$$

Our work extends the work of Helfgott \cite{helf} who proved similar results for the 
family $\{\SL_2(\F_p): p\tx{ prime}\}$.



\newpage
\tableofcontents
\newpage

\sec{Introduction}

\ssec{Background}

Let us define the \ti{directed diameter} of a finite group $G$ with respect to a set of
generators $S$ to be the minimal number $l$ for which any
element in $G$ can be written as a product of at most $l$ elements
in $S$. We denote this number by $\diam^+(G,S)$. 
Define the (undirected) \ti{diameter} of a finite group $G$ with respect to a set of
generators $S$ to be $\diam(G,S):= \diam^+(G,S\cup S^{-1})$.

The diameter of groups has many applications. Aside from group theory (see
\textbf{\cite{BKL_diam1,Larsen_diam3,LS_diam4}}) and
combinatorics(see \cite{Dixon,ER,ET1,ET2}) the diameter of groups
shows up in computer science areas such as communication networks
(see \cite{Sto,PV}), generalizations of Rubik's puzzles (see
\cite{DF,McK}), algorithms and complexity (see \cite{EG,Je}). For
a detailed review see \cite{BHKLS_diam2}.

Since we are interested in the ``worst case generators'', we define 
$$\diam(G):=max\{diam(G,S):G=\gen{S}\}.$$
A family of finite groups $\set{G_n:n\in \N}$ is said to have \tb{poly-log diameter} (resp. \tb{log diameter}) 
if for any $n\in\N$ we have $$\diam(G_n)\leq C \log^d(|G_n|)$$ for some constants $C,d>0$ (resp. for $d=1$).

In \cite{dinai}, the author shows (with an effective algorithm) that for any fixed $p,m\in \N$ with $p$ a prime and $p>m\geq 2$, 
the family $$\G_{m,p}:=\set{\SL_m(\Z/\!p^n\Z):n\in \N}$$ has poly-log diameter. 
Abert and Babai  \cite{abert-babai} showed that for any fixed prime $p_0$, 
the family $\set{C_{p_0} \wr C_p: p\tx{ prime; }p\neq p_0}$ has logarithmic diameter.

A long standing conjecture of Babai \cite{bs2} asserts that the family of non-abelian finite simple groups has a poly-logarithmic diameter. 
Very little is known about this conjecture.
See \cite{bs1} and \cite{bs2} for some partial results concerning the alternating groups.

A breakthrough result of Helfgott \cite{helf} proves the conjecture for the family $\set{\SL_2(\F_p):p\tx{ prime}}$. 
The main goal of this paper is to extend Helfgott work to the family $\set{\SL_2(\F_{p^n}):p\tx{ prime}; n\in \N}$. We follow the basic strategy of Helfgott (with some short cuts following \cite{bg su(2)}) and in particular we also appeal to additive combinatorics and sum-product theorems. 
The new difficulty is that unlike fields of prime order, general finite fields have subfields, and subsets which are ``almost'' subfields - which are ``almost'' stable \WRT\ sum and product.

\ssec{Main results}
Our main results are the following.
\doTHM{[See Theorem \r{main theorem} in \Sr{ssec:main}]\l{thm:A^{(3)}:growth}
\TES\ $\eps\in \R_+$ \ST\ \tFHS\ for any finite field $\F_q$. 
Let $G$ be the group $\SL_2(\F_q)$ and let $A$ be a generating set of $G$.
\TWH\footnote{The same assertion holds for $\PSL_2(\F_q)$.}, 
$$|A\ddot A\ddot A|\ge \MIN{|A|^{1+\eps},|G|}.$$
}
From this we easily get the following.

\doCORO{[See corollary \r{corr:diameter} in \Sr{ssec:main}]\l{corr:main:diameter}
\TE\ $C,d\in \R_+$ \ST\ \tFHS\ for any finite field $\F_q$. 
Let $A$ be a subset of generators of $G=\SL_2(\F_q)$.
\TWH, 
\doeqnn{\diam^+(G,A) < C \log^d(|G|)}	
and for any $\del\in \R_+$ we have,
\doeqnn{|A|>|G|^{\del} \Rightarrow \diam^+(G,A) < C\lr{\tfr{\del}}^d. } 
}
 
\ssec{Organization of the manuscript}
The manuscript is organized as follows:
In \Sr{sec:Preliminaries} we bring notations and definitions, which are required for this work, as well as mathematical background. 
In \Sr{sec:Tools from Additive combinatorics} we collect useful facts from Additive Combinatorics to be used later.
In \Sr{sec:Useful properties of SL_2(F)} we prove some useful facts about $\SL_2(\F_q)$.
In \Sr{sec:Growth properties of SL_2} we extend few of the main ingredients from the proof of Helfgott, from $\SL_2(\F_p)$ to $\SL_2(\F_q)$. 
In \Sr{sec:Main results} we show how to use all the previous sections in order to prove the main results of this manuscript. 
In \Sr{sec:Further conjectures and questions} we present some questions/conjectures.

\sec{Preliminaries}\l{sec:Preliminaries}

\ssec{Notations}

We will use the following notations. 
$\log x$ will stand for $\log_2 x$, $\log$ in the base $2$. 
We will always use $p$ for a prime number and $q$ for a prime power. 
For a subset $A\subseteq B$ and $x\in B$ denote for short $A\SM \{x\}$ by $A\SM x$
and similarly $A\cup x :=A\cup \{x\}$. For a field $\F$, denote by
$\overline{\F}$ some fixed algebraic closure of $\F$. 
We denote $$(G,\cdot)$$ a multiplicative group which is not necessarily 
commutative and $$(G,+)$$ will stand for a commutative additive group.


%


\doDEF{
Let $G$ be a group and let $A,B,A_1,\ldots,A_n\subseteq G$ be non-empty subsets of $G$. 
For $k\in\Z$ denote 
\doalign{
	\mbf{A^k}			&:=	\set{a^k: a \in A}\\
	\mbf{A^{\pm 1}}	&:=	A\cup A^{-1}
}
Define the product-set, 
$$\mbf{A\ddot B}:=\set{a\ddot b:a\in A,b\in B}$$
and for $x\in G$ define $x\ddot A := \set{x}\ddot A$ and $A\ddot x := A \ddot \set{x}$. 
Denote the product set of $A_1,\ldots,A_n$ by
$$\prod_{i=1}^n A_i :=\set{a_1\cdots a_n :\forall 1\le i\le n, a_i \in A_i}$$
and the product set of one set with itself $n$-times by 
$$\mbf{A^{(n)}}	:= \prod_{i=1}^n A.$$
The most important notations in this manuscript will be 
\doalign{
	A^{[0]} &:=\set{1}\\
	A^{[1]} &:= A^{\pm} \cup 1	\\
	\mbf{A^{[n]}} &:= (A^{[1]})^{(n)}
}
the set of words of length at most $n$ in the letters $A^{\pm}:=A\cup
A^{-1}$. Note that in general we have only the containments $$A^n \subseteq A^{(n)} \subseteq A^{[n]}.$$
}

\doSFACT{\l{sfact:A^{xyz}}
Since we have three possible operations on the subsets\footnote{Note the these operations \tb{on subsets} of $G$ are not induced from operations on elements of $G$.}, $A^{[m]},A^{(n)}$ and $A^k$, we use the following ``group action'' notation $A^{gh}=(A^g)^h$. 
For example, 
\doalign{
	A^{(n)[m]}	&:=	(A^{(n)})^{[m]}\\
	A^{[m](n)}	&:=	(A^{[m]})^{(n)}.	
}
Similarly $$A^{xyz}:=((A^x)^y)^z$$ when $x,y,z$ is any of these operation e.g., $A^{k(n)[m]}:=((A^{k})^{(n)})^{[m]}$.
Note the these operations on subsets is associative 
$$A^{(xy)z}=A^{x(yz)}=((A^{x})^y)^z.$$

Note that in general we have only the containments, $$A^{(n)[m]}\subseteq A^{[m](n)}=A^{[nm]}.$$
We can write these properties as a table of relations between the operations as 
$[n][m]=[nm]$ and $(n)(m)=(nm)$ and $[mn]=[m](n)\neq (n)[m]$.

Note that if $\gen{A}$ is abelian then $A^{k(n)}=A^{(n)k}$ and similarly 
\doeqnn{A^{k[m]}=A^{[m]k}.}	
}



\doDEF{\l{def:[g,h]}
Let $G$ be a group and let $g,h\in G$. We will denote by
\doalign{
	\mbf{g^{h}}	&:= h^{-1} g h\\
	\mbf{\lrb{g,h}}	&:=	g^{-1}g^h = g^{-1}h^{-1}gh
	}
For subsets $A,B \subseteq G$ we denote by
\doeqnn{A^B:=\set{a^b:a\in A,b\in B}}	
and $x^B:=\set{x}^B$ for short. 
For commutator of two subsets we will write 
\doeq{\mbf{[A,B]_{set}}:=\{[a,b]:a\in A;b\in B\}. \l{eq:[A,B]}} 
Note that we have only containments 
\doeqnn{[A,B]_{set}\subseteq A^{-1}A^B\subseteq A^{-1}B^{-1}AB.}
}

\doDEF{
Let $G$ be a group and let $A,B\subseteq G$. Define
\doeqnn{\opC_B(A)	:= \set{b\in B:a^b=a \tx{ for all }a\in A}.} 
}

\doSFACT{\l{sfact:A^{g(n)}}
Note that using these notations we always have $g(n)=(n)g$ for $g\in G$ and $n\in\N$. I.e.,
\doeqnn{A^{g(n)}=A^{(n)g}=(A^g)^{(n)}=(A^{(n)})^g} 
and similarly $g[m]=[m]g$ and $kg=gk$.
So conjugation (or any other automorphism) commutes with the operations $A^{[m]},A^{(n)},A^k$.
}

\doDEF{
We will use the generation notation $\gen{A}$
depending on the category we are working. The categories that will
be involved in the manuscript will be groups and rings.
}

\doDEF{\l{def:Prod and Diag}
Let $g=\dopmatrix{a & b \\c & d}\in \SL_2(\F)$. Denote, 
\doalign{
	\mbf{\Prod(g)}	&:=	a\ddot d\\
	\mbf{\Diag(g)}		&:=	\lr{a,d}.
}
Extend these functions to $\Prod(V)$ and $\Diag(V)$ for subset $V\subseteq \SL_2(\F)$.
}


\doDEF{\l{def:D_X}
Let $g=\dopmatrix{a & b \\c & d}\in \SL_2(\F)$ and $x\in \F^{\times}$. Denote, 
\doalign{
	\mbf{D_g}				&:=	\dobmatrix{a & 0 \\  0 & d}\\
	\mbf{D_{(a,d)}}	&:= \dobmatrix{a &0 \\0 &d}\\
	\mbf{D_x}				&:= \dobmatrix{x &0 \\0 &x^{-1}}
}
Extends these notations to subsets in the obvious way $D_X:=\set{D_x:x\in X}$ where $X$ is either $X\subseteq \F^{\times}$ or $X\subseteq \F\times \F$ or $X\subseteq \SL_2(\F)$.
}


\doDEF{ 
For positive real-valued functions, we write $f\ll g$ if $f=O(g)$. Similarly we write $f\gg
g$ if $g\ll f$, and $f\approx g$ if $f\ll g\ll f$. Similarly we
will use the dual notation $f=\Om(g)$ for $g=O(f)$.
Denote also $$f\sim g \iff \fr{2}f \leq g \leq 2f.$$
}

\doSFACT{\l{sfact:Oeps}
Let $\eps\in \R_+$ be real number with $\eps<\fr{2}$.
Then we have 
\doearray{
	1-\eps&<\tfr{1+\eps}<&1-\tfr{2}\eps\\
	1+\eps&<\tfr{1-\eps}<&1+{2} \eps.
}
Therefore for any $X,Y\in \R_+$ we have in the $\Om$-language:
$$X \ll Y^{1+\Oeps} \iff X^{1-\Oeps} \ll Y$$ 
and similarly 
$$X^{1+\Omeps}\ll Y \iff X\ll Y^{1-\Omeps}.$$


}


\doDEF{\l{def:non zero-divisors}
Let $R$ be a ring (not necessarily commutative) and let $a\in R$. 
Define the endomorphisms $L_a$ and $R_a$ by $L_a(b)=ab$ and $R_a(b)=ba$.
Then $L_a$ is endomorphism of the right\footnote{the action of the scalars is from the right.} $R$-module $R$ and $R_a$ is endomorphism of the left\footnote{the action of the scalars is from the left.}  $R$-module $R$.
Denote the right ideal $\Ker(L_a)$ and the left ideal $\Ker(R_a)$ by
\doearray{
	\Ker(L_a)&:=&\set{b:ab=0}\\
	\Ker(R_a)&:=&\set{b:ba=0}.
}	
Now suppose $R$ is commutative ring. 
Denote the set of non zero-divisors in $R$ by $R^{\times}$: 
\doeqnn{a\in R^{\times} \iff L_a\tx{ is injective }\iff \Ker(L_a)=0.}
}

If $A$ is a subset of a commutative ring $R$ we will need different notations to distinguish the product-set 
$A\ddot A=\set{ab:a,b\in A}$ 
and the sum-set 
$A+A=\set{a+b:a,b\in A}$. 
Therefore we will need in some situations the following definitions.

\doDEF{\l{def:sum-set}
Let $A$ be a subset of an additive (semi) group $G$ and let 
$n\in \N$.
Denote by
\doearray{
	\sum_n A &:=\{a_1+\ldots +a_n :\forall i, a_i \in A\}.
}
}


\doDEF{\l{def:mult(G)}
Let $\Gamma\subseteq X\times Y$ be a directed graph.
Denote the inverse (opposite) graph $\Gamma^{-1}\subseteq Y\times X$ (or $\Gamma^{\tx{op}}$) by
$$\Gamma^{-1}:=\set{(y,x):(x,y)\in \Gamma}.$$
Let $A\subseteq X$ and $a\in X$.
Denote
\doearray{
	\Gamma_a &:=& \set{y\in Y: (a,y)\in \Gamma}\\
	\Gamma(A)&:=& \bigcup_{a\in A}\Gamma_a
}
Denote $$\deg(\Gamma):=\max_{x\in X}\set{|\Gamma_x|}.$$
Clearly $\deg(\Gamma)\leq d \Rightarrow |\Gamma(A)|\le d|A|$ for any $A\subseteq X$.
We will say that $\Gamma$ is $d$-regular\footnote{or we write for short, $\Gamma$ is $(1:d)$.}, if $|\Gamma_x|=d$ for all $x\in X$. 
We define the multiplicity of $\Gamma$ to be $$\mult(\Gamma):=\deg(\Gamma^{-1}).$$
}

We will use the previous definition with the following simple observations. 
\doSFACT{\l{sfact:mult(f)}
A function $f\in Y^X\subseteq X\times Y$ is a directed graph which is $1$-regular graph. 
Therefore we get, 
$$\mult(f)\leq n \RightarrowB |f(A)|\geq |A|/n\tx{ for any }A\subseteq X.$$

For example, any one variable polynomial $0\neq f(x)\in \F[x]$ of degree $d$ defines a substitution map $f_s:\F\rightarrow\F$ \ST
$$\mult(f_s)\le \deg(f).$$

Similarly if $0\neq f(x,x^{-1})\in \F[x,x^{-1}]$ with $deg_x(f)+deg_{x^{-1}}(f)=d$ then $\mult(f_s)\leq d$ where 
$$f_s:\F^{\times}\rightarrow \F.$$
E.g., $f(x)=x^2+x^{-3}$ has multiplicity $\leq 5$.
By abuse of notation we will write for $f\in \F[x,x^{-1}]$,
\doeqnn{\mbf{\mult(f)}:=\mult(f_s).} 
}

\ssec{Uniform poly-logarithmic diameter bounds}

\doDEF{\l{def:diam(graph)}
For a finite (undirected) graph $\Gamma=(V,E)$ define $\diam(\Gamma)$, the \tb{diameter} of the graph $\Gamma$, to be the
minimal $l$ such that \ti{any} two vertices are connected by a path
with at most $l$ edges. Set the diameter to be $\diam(\Gamma)=\infty$ if the graph is not
connected.
}

\doDEF{\l{def:diam(group)}
For a finite group $G$ and a subset $S$ of $G$, define
\doeqnn{\diam(G,S):=\diam(\Cay(G,S)).}
For a finite group $G$ and a set of generators $S$ of $G$, we have
\doeqnn{\diam(G,S)=\min\{k:S^{[k]}=G\}}
Define the maximal diameter of $G$ to be
\doeqnn{\diam_{\max}(G):=\max\{\diam(G,S): S\subseteq G, \gen{S}=G\}}
or just for short $\mbf{\diam(G)=\diam_{\max}(G)}$.
For a finite group $G$ and a set of generators $S$ of $G$, define
\doeqnn{\mbf{\diam^+(G,S)}:=\min\{k: (S\cup 1)^{(k)}=G\}.}
} 

\addREM{
It is easy to see that for a set of generators $S$ of $G$ with
$s:=|\SS|$ we have
\[\log_s |G|-1 \leq \diam(G,S)\leq|G|-1\]
by a simple count of words in $G$ with the letters $S\cup S^{-1}$.
Still there is an exponential gap between these two bounds. 
So usually the goal is to find either an \ti{upper} logarithmic or
a poly-logarithmic diameter bound for $\diam_{\max}(G)$. 
This bound is of interest when each group, in the family of groups, 
can be generated by a subset of bounded size.
}


A well known conjecture of Babai asserts the following \cf{bs1,bs2}.

\begin{conj}[Babai]\label{conj babai}
There exist $C,d\in \R_{+}$ \ST\ for any non-abelian finite simple group $G$ we have 
$$\diam_{\max}(G)\le C \log^d|G|.$$
This bound may even be true for $d=2$, but not
$d<2$, as the groups $Alt(n)$ demonstrate.
\end{conj}

The first step towards proving Babai's conjecture was made by Helfgott \cfr[\S1.2 Main Theorem]{helf}.

\doTHM{[Helfgott]\label{thm:helf}
Denote the family of groups 
$$\G=\set{\SL_2(\F_p): p \tx{ prime}}.$$ 
There exist $C,d\in \R_{+}$ \ST\ for any $G\in \G$ we have, 
$$\diam_{max}(G)\leq C \log^d|G|.$$
}

We extend this theorem to all finite fields to get the following.
\doTHM{[See corollary \r{corr:diameter} in \Sr{ssec:main}]\label{thm:oren_diameter_bound}
Denote the family of groups 
$$\G=\set{\SL_2(\F_{p^n}): p \tx{ prime}; n\in \N}.$$ 
There exist $C,d\in \R_{+}$ \ST\ for any $G\in \G$ we have, 
$$\diam^+_{\max}(G)\leq C \log^d|G|.$$
}

The main idea in Helfgott's work is to show
an expansion property of subsets w.r.t the product operation in
the group. For this he reduced the problem to an expansion
property of the addition and multiplication operations in the
underline fields. One advantage of these results is their,
relatively, elementary proofs. One disadvantage of these results
is that they do not supply a algorithm(/method) for actually
calculating such a short paths(/products) in the graphs(/groups).

\sec{Tools from Additive combinatorics}\l{sec:Tools from Additive combinatorics}

\ssec{The fundamental tools} 

\sssec{Ruzsa triangle inequality}

The following useful lemma of Ruzsa
allows one to pass from control of sum-set to control of
minus-sets (cf. \cite[Lemma 2.6]{tv} and \cite[\S 2.3 Lemma 2.1]{helf}).

\doLEM{[Ruzsa]\l{Ruzsa triangal inequality}	
Let $G$ be a group and let $A,B,C\subseteq G$ be finite subsets. 
Then we have, 
\doeq{\l{lemma_ruzsa_bound_1} 
	|AB||C|\leq |AC^{-1}||CB|.
}
}

\doPROOF{ 
Define the product map $p:AC^{-1} \times CB \rightarrow G$ by $p(x,y)=xy$. 
Then for any $a\in A$ and $b\in B$ we have, 
$$p^{-1}(ab) \supseteq \set{(ac^{-1},cb):c\in C}$$ so $|p^{-1}(ab)|\geq |C|$. 
Therefore $|AC^{-1}||CB|\geq |p^{-1}(AB)|\geq |C||AB|$ \sWAD.
}

In particular by taking $B=C=A^{-1}$ we get the following Corollary.
\doCORO{\l{coro:ruzsa:sum to minus}
Let $G$ be a group and let $A \subseteq G$ be a finite subset.
For any $1\leq K \in \R$ we have, 
$$|A\ddot A|\leq K|A|\Rightarrow |AA^{-1}|\leq K^2|A|.$$
}

\doPROOF{
By lemma \r{Ruzsa triangal inequality} we get
\doalign{
	\fr[|A\ddot A^{-1}|]{|A|}	&\leBY{lemma_ruzsa_bound_1} \fr[|A\ddot A|]{|A|} \fr[|A^{-1}\ddot A^{-1}|]{|A|}\\
											&\eqB \lr{\fr[|A\ddot A|]{|A|}}^2\\
											&\leB K^2.	\QED	
	}
}

%

\doDEF{\l{def:D(A,B)}
Let $G$ be a group and let $A,B\subseteq
G$ be finite non empty subsets. Define, 
$$D(A,B):=\frac{|AB^{-1}|}{|A|^{1/2}|B|^{1/2}}.$$
Define the \tb{Ruzsa distance} between $A$ and $B$ to be 
$$d(A,B):=\log(D(A,B)).$$
}

It is easy to see that the following properties hold.
\doSFACT{\l{sfact:dilation invariant}
Let $G$ be a group and let $\emptyset\neq A,B\subseteq G$ be finite subsets. 
Then for any $x,y\in G$ we have, $$d(A,B)=d(B,A)=d(xA,yB)=d(Ax,Bx)\geq 0.$$
}

As an immediate consequence of lemma \r{Ruzsa triangal inequality} we get,
$$D(A,B)\leq D(A,C)D(B,C)$$ therefore we get the following Triangle inequality.

\doSFACT{\l{sfact:is quasi-metric}
Let $G$ be a group and let $\emptyset\neq A,B,C\subseteq G$ be finite subsets. 
Then we have, $$d(A,B)\leq d(A,C)+d(B,C).$$
Therefore $d(A,B)$ is quasi-metric\footnote{actually $d(A, B)=0 \iff A,B$ are both left cosets of some finite subgroup $H\le G$ (see \cite[Proposition 2.38]{tv}).} on the set of finite subsets of $G$.
}

\sssec{\Plunnecke-Ruzsa inequality}\label{sssec:Plunnecke-Ruzsa}

The following theorem of \Plunnecke-Ruzsa allows one to pass from
control of sum-set to control of iterated sum-sets (cf. \cite[\S 6.5, Corollary 6.29]{tv}).

\doTHM{[\Plunnecke-Ruzsa] \label{thm:Plunnecke-Ruzsa}
Let $(G,+)$ be an additive group and let $A,B\subseteq G$ be finite subsets.
Suppose
\doeq{|A+B|\leq K |B|}
for some $1\leq K\in \R$. Then for any $n,m\in\N$ we have, 
\doeq{\label{ieq_PR_1}|\sum_n A|\le K^n|B|\mqbq{and}|\sum_n A-\sum_m A|\le K^{n+m}|B|.}
} 

\IP\ we get the following.
\doCORO{\l{coro:plunnecke:from minus to plus} 
Let $(G,+)$ be an additive group and let $A\subseteq G$ be finite subset.
Then for any $1\leq K\in \R$ we have,
\doeq{\l{eq:plunnecke:from minus to plus} |A - A|\leq K|A| \RightarrowB |A + A|\leq K^2|A|.}
}

\doPROOF{
By taking $B=-A$, we are done by theorem \r{thm:Plunnecke-Ruzsa}.
}

Another special case of theorem \r{thm:Plunnecke-Ruzsa} is the following result.
\doCORO{ 
Let $R$ be a commutative ring and $A\subseteq R$ a finite subset and let $b\in R^{\times}$.
Suppose 
\doeqnn{|A+bA|\leq K|A|}
for some $1\leq K\in \R$. 
Then we have, 
\doeq{\label{plunnecke from twisted sum to sum and minus}|A + A|\leq K^2|A|\mqbq{and}|A - A|\leq K^2|A|.}
}
\doPROOF{
By taking $B=bA$, we are done by theorem \r{thm:Plunnecke-Ruzsa}.
}

\addREM{
Note that actually we only used the fact that the addition in $R$ is commutative and that $|A|=|bA|$. 
Therefore this statement is true also for non-commutative rings provided that $L_b$ is injective. 
}

\sssec{From large growth to large tripling}

In corollary \r{coro:plunnecke:from minus to plus} one cannot drop the additive assumption to get polynomial bound  
like  \re{eq:plunnecke:from minus to plus} (cf. \cite[\S 2]{helf}). 
However, one can deduce easily from Lemma \r{Ruzsa triangal inequality} the following result.

\doLEM{\CITE{\S 2.3 Lemma 2.2}{helf}\l{lem:from growth to tripling}
Let $G$ be a group and let $A \subseteq G$ be a finite subset. 
\doAlign{
\addtx{For any $1\leq K\in \R$ and $x_1,x_2,x_3\in \set{\pm 1}$ we have,}
	|A^{(3)}|	\leq K|A| &\RightarrowB	|A^{x_1} A^{x_2}A^{x_3}|	\leq K^3|A|. \l{eq:from growth to tripling:1}\\
\addtx{For any $3\leq n \in \N$ and $1\leq K\in \R$ we have,}
	|A^{[3]}|	\leq K|A|	&\RightarrowB	|A^{[n]}|	\leq K^{n-2}|A|. \l{eq:from growth to tripling:2}\\
\addtx{\IP\ for any $3\leq n\in \N$ and $1\leq K\in \R$ we get,}
	|A^{[n]}| > K |A| &\RightarrowB |A^{(3)}| > \tfr{2} {\sqrt[3n]{K}}|A|. \l{eq:from growth to tripling:3}
	}
}

\doPROOF{
\BTA, 
\doeqnn{
	|A^{-1}A^{-1}A^{-1}|=|AAA|\leq K|A|.
}
Therefore by Lemma \r{Ruzsa triangal inequality} we get,
\doearray{
	|AAA^{-1}|	& \byeq{lemma_ruzsa_bound_1}{\leq}& \fr{|A|}|AAA||A^{-1}A^{-1}|\\
							&\leq& \fr{|A|}|AAA||A^{-1}A^{-1}A^{-1}|\\
							&=&  \lr{\frac{|AAA|}{|A|}}^2|A|\\
							&\leq& K^2|A|
}
Therefore we get also, 
\subeq{
\doeq{\l{proof_ruzsa for k products_bound_a}
	|AA^{-1}A^{-1}|=|AAA^{-1}|\leq K^2|A|.
}
By repeating the previous argument but now with $\mbold{A}=A^{-1}$ (i.e., $A^{-1}$ in the roll of $A$) we get
\doeqnn{\l{proof_ruzsa for k products_bound_b}
	|A^{-1}A^{-1}A|,|A^{-1}AA|\leq K^2|A|.
}
\OTOH,
\doearray{
|A^{-1}AA^{-1}|	&	=	&	|AA^{-1}A|\\
\nn		& \byeq{lemma_ruzsa_bound_1}{\leq}	& \fr{|A|}|AA^{-1}A^{-1}||AA|\\
\nn		&\leq & \frac{|AA^{-1}A^{-1}|}{|A|}\frac{|AAA|}{|A|}|A|\\
\nn		&\byeq{proof_ruzsa for k products_bound_a}{\leq} & K^3|A|. 
}
}
Therefore we are done with the bound  \re{eq:from growth to tripling:1}. \good

By induction for $n\geq 3$ we get from Lemma \r{Ruzsa triangal inequality} that,
\doearray{
	|A^{[n+1]}| &\byeq{lemma_ruzsa_bound_1}{\leq}&	\fr{|A|}|A^{[n-1]}A||A^{-1}A^{[2]}|\\
							&\leq& \frac{|A^{[n]}|}{|A|}\frac{|A^{[3]}|}{|A|}|A|\\							
							&\leq& K^{n-1}|A|. 
}
\sWAD\ with  \re{eq:from growth to tripling:2}.\good

If we combine  \re{eq:from growth to tripling:1} and  \re{eq:from growth to tripling:2} we get for any $n\geq 3$,
\doalign{
								|A^{(3)}|		&\le	K|A|\\
\RightarrowL\qu	|A^{[3]}|		&<		(2K)^3|A|\\
\RightarrowL\qu	|A^{[n]}|		&\le	(2K)^{3(n-2)}|A|\\
\RightarrowL\qu	|A^{[n]}|		&<		(2K)^{3n}|A|.	
}
Therefore by negating the inequalities we get,
\doeqnn{|A^{[n]}| \geq K|A|	\qu\RightarrowL\qu |A^{(3)}|	> \fr{2}K^{1/(3n)}|A|}
\sWAD\ with  \re{eq:from growth to tripling:3}.\good
}
%
%

\ssec{Expansion properties in fields}

When dealing with fields one can use the following Sum-Product theorem \CITE[cf. ]{\S 2.8}{tv} 
which is a slight improvement of \cite{bkt,bkt1}.

\doTHM{\CITE{Theorem 2.52}{tv}\label{freiman thm}	
There exists an absolute $C>0$ \ST\ \tFHS\ for any $1\leq K\in\R$ and any field $\F$.
Let $A\subseteq \F$ be a finite subset and suppose
\doeqnn{|A+A|+|A\ddot A|\leq K|A|.}
Then either
$|A|< CK^C$ or for some subfield $\E\leq\F$ and $x\in\F^{\times}$ we have, 
\doeqnn{|\E|\leq CK^C|A|\qtq[and]|A\SM x\E|\le CK^C.}
}

The power of this quantitative theorem is that if a set is almost stable under the two field's operations
then as a set it is almost a field, \tb{up to a polynomial
lost}. We will be interested in subsets with large growth:
$$\max\set{|A+A|,|A\ddot A|} \sim |A+A|+|A\ddot A| \gg |A|^{1+\eps}.$$
Therefore we will use the following definition.

\doDEF{[Almost fields]\label{def:almost field}  
Let $\F$ be a field and let $A\subseteq \F$ be a finite subset and let $\eps\in \R_+$.
We will say that $A$ is $\mbold{\eps}$\tb{-almost field}, or $\mbold{\eps}$\tb{-field} for short, 
if for some subfield $\E\le \F$ and $x\in\F^{\times}$ we have,
\doeq{\l{eq:almost field}|\E|\le |A|^{1+\eps}\mqbq{and}|A\SM x\E|\le|A|^{\eps}.}
If the above holds then we will say that that $A$ is $\mbold{\eps}$\tb{-field} $\mbold{\E}$.
Define $A$ to be \tb{pure} $\eps$-field if 
\doeq{\l{eq:pure}|\E|\le |A|^{1+\eps}\mqbq{and}A\subseteq \E.}
If  \re{eq:almost field} holds but  \re{eq:pure} does not hold then we will say that $A$ is an impure $\eps$-field. 
\IOW, $A$ is \tb{impure} $\eps$-field if  \re{eq:almost field} holds and also
\doeq{\l{eq:impure}|A\SM \E|>0.}
}

\doDEF{[Almost stable subsets]\label{def:almost stable} 
Let $\F$ be a field, $A\subseteq\F$ be a finite set and let $\eps \in \R_+$. 
We will say that $A$ is $\mbf{\eps}$\tb{-close}, or $\mbf{\eps}$\tb{-stable}, if 
\doeq{\l{eq:closed}|A\ddot A|+|A+A| \leq |A|^{1+\eps}.}
Otherwise, we will say that $A$ has $\eps$-expansion, or $\eps$\tb{-growth}.
}

Let's restate Theorem \r{freiman thm} using this terminology.

\doTHM{\l{thm:freiman2}
There exists $C>0$ such that the following holds for any $\eps\in \R_+$ with $\eps<\tfr{C}$. 
Let $\F$ be a field and let $A\subseteq \F$ be a finite subset of size $|A|>C^{1\!/\eps}$. 
Then we have,
\doAlignAT[3]{
	&A\tx{ is \eps-field}\qu 	&&\RightarrowL&& \qu A\tx{ is $C$\eps-stable.}\l{eq:thm: freiman 2.1}\\
	&A\tx{ is \eps-stable}\qu &&\RightarrowL&& \qu A\tx{ is $C$\eps-field.}\l{thm: freiman 2.2}
	}
}

\addREM{
The statement  \re{eq:thm: freiman 2.1} is trivial, as we shall see in the proof below. 
The important part of the theorem is  \re{thm: freiman 2.2}.
The theorem can be stated as follows: 
For any $\eps>0$ which is small enough, if $A$ is big enough (depending on $\eps$), both  \re{eq:thm: freiman 2.1} and  \re{thm: freiman 2.2} hold.
}

\doPROOF{
Suppose $A$ is $\eps$-field. Therefore by \r{def:almost field} we get,
\doeqnn{|\E|\le |A|^{1+\eps}\mqbq{and}|A\SM x\E|\le|A|^{\eps}.}
Denote $X:=A\SM x\E$ and so we get,
\doalign{
	|A + A|	&\le |(x\E\cup X)+(x\E\cup X)|\\
					&\le |\E|+|\E||X|+|X|^2\\
					&\le 3|A|^{1+\eps}
}
and similarly the same bound for $|A\ddot A|$. Therefore if $|A|^{\eps}\ge 6$
we get 
\doalign{
	|A + A|+|A\ddot A| &\le 6|A|^{1+\eps}\\
										 &\le |A|^{1+2\eps}
}
\sWADW  \re{eq:thm: freiman 2.1}.\good

Now suppose $A$ is $\eps$-stable. Denote $K:=|A|^{\eps}$ so by \r{def:almost stable} we get,
$$|A\ddot A|+|A+A| \leq K|A|.$$
Therefore by Theorem \r{freiman thm} \tFHS\ for some absolute\footnote{The constant $C>0$ from Theorem \r{freiman thm} is absolute and do not depend in $\eps$.} $C_1>0$. 
Either 
\doeq{\l{eq:C_1K^{C_1}}|A|< C_1K^{C_1}}
or for some subfield $\E\leq\F$ and $x\in\F^{\times}$ we have, 
\doeq{\l{eq:subfield}|\E|\leq C_1K^{C_1}|A|\mqbq{and}|A\SM x\E|\le C_1K^{C_1}.}
Therefore if $\eps$ is small enough, say $\eps<\tfr{2C_1}$, and $|A|^{\eps}$ is big enough, say $|A|^{\eps}>C_1$, then 
$$C_1K^{C_1}=C_1|A|^{C_1\eps}<|A|^{2C_1\eps}<|A|.$$
Therefore  \re{eq:C_1K^{C_1}} does not hold and from  \re{eq:subfield} we get that $A$ is $2C_1\eps$-field, \sWADW  \re{thm: freiman 2.2}.
}

We can state the non trivial part of theorem  \re{thm: freiman 2.2} in the $\Om$-language as follows:
\doCORO{\l{coro:freiman}
There exists $C>0$ \ST\ for any $\eps>0$ small enough the following hold for any finite subset $A\subseteq \F$ which is big enough.
$$A\mb{ is not }\eps\mb{-field}\qu\RightarrowL\qu A\mb{ has }\Omeps\mb{-growth}.$$
}
 
\ssec{Expansion functions in fields} \l{ssec:Expansion functions in fields}

We begin by introducing some new notations.
\doDEF{\l{def:Tr_g}
Let $\F$ be a field and let $g\in \GL_n(\F)$. 
Define, 
\doalign{
	\Tr_g(A,B)	&:=	\Tr(AB^g)\\						
	\addtx{for any $A,B\in \M_n(\F)$ and denote for $V\subseteq M_n(\F)$,}
	\Tr_g(V)		&:=\set{\Tr(AA^g):A\in V}.
	}
}

\doDEF{\l{def:tr_g}
Let $\F$ be a field, $x,y\in \F^{\times}$ and $g\in \SL_2(\F)$.
Define $\tr:\F^{\times}\rightarrow \F$ and $\tr_g:\F^{\times}\times \F^{\times}\rightarrow \F$ by
\doalign{
	\tr(x)			&:= \Tr(D_x)\\
	\tr_g(x,y)  &:= \Tr(D_x (D_y)^g).
	}
Extend these definitions to $\tr(X)$ and $\tr_g(X,Y)$ for subsets $X,Y\subseteq \F^{\times}$.	
}

We immediately get the following equivalent definition.
\doSFACT{\l{sfact:tr_g}
Let $g=\dobmatrix{a & b\\ c & d}\in \SL_2(\F)$ and $x,y\in \F^{\times}$.
\TWH,
\doalign{
	\tr(x)			&= x+x^{-1}\\
	\tr_g(x,y) 	&= ad\ddot \tr(xy)-bc \ddot \tr(x/y).
	}
}

\doDEF{\l{def:tr_t}
Let $x,y\in \F^{\times}$ and let $t\in \F$. Define,
\doeqnn{\tr_t(x,y) := t\ddot \tr(xy)+(1-t)\ddot \tr(x/y)}
}


As a consequence of \r{sfact:tr_g} and \r{def:tr_t}, we immediately deduce the following.

\doSFACT{\l{sfact:tr_g(x,y)} 
Let $x,y\in \F^{\times}$ and $g\in \SL_2(\F)$ with $t=\Prod(g)$. 
Then we have,
\doEQsplit{eq:tr_g(x,y)}{	
	\Tr(D_x (D_y)^g) 	&= \tr_g(x,y)\\
			&= \tr_t(x,y)\\
			&=	t\ddot \tr(xy)+(1-t)\ddot \tr(x/y).
	}
}
%
\addREM{
\IP\ from  \re{eq:tr_g(x,y)}, we get that
\doalign{ 
	\Prod(g)=1 &\Rightarrow \tr_g(x,y) = \tr(xy)\\
	\Prod(g)=0 &\Rightarrow \tr_g(x,y)=\tr(x/y).
}
Note that 
\doeqnn{\Prod(g)=1 \iff g \tx{ is triangular}} 
i.e., $g$ is either upper triangular or lower triangular.
}

We make the following easy observations in any field $\F$.

\doSFACT{\l{fact:trace identity for scalars}
Let $\F$ be a field and let $G=\F^{\times}$ be its multiplicative group.
Let $x,y\in \F^{\times}$ and $X,Y \subseteq \F^{\times}$. Then we have, 	
\doAlign{
	\tr(x)\tr(y)	&= \tr(xy)+\tr(xy^{-1})	\l{fact:eq:trace identity for scalars}\\
	\addTX{and therefore,}
	\tr(X)\tr(Y) 	&\subseteq 	\tr(XY)+\tr(XY^{-1}).\nt\\
	\addTX{and in particular,}	
	\tr(X)\tr(X) &\subseteq \tr(X^{[2]})+\tr(X^{[2]}).\l{fact:eq:prod_sum_bound} 
	}
}

\doPROOF{
 \re{fact:eq:trace identity for scalars} is trivial from the definition of $\tr(x)=x+x^{-1}$.
The two other equations immediately follow from  \re{fact:eq:trace identity for scalars}.
}

%
%
%
%

The following striking reduction of Helfgott\footnote{The following proof is due to Helfgott and is different from his original proof.}
allows one to gain large expansion from
the non commutativity in the group by twisting properly some
commutative sets (cf. \cite[\S 3]{helf} and \cite[\S 4]{bg su(2)}).

\doTHM{[Helfgott]	
\l{helfgott: reduction from trace function to sum product} 
There exists $C>0$ \ST\ the following holds for any field $\F$.
Let $X\subseteq \F^{\times}$ be a finite subset and suppose 
\doeq{\l{prop.3.18.1.if}|\{a_1\ddot \tr(xy)+a_2\ddot \tr(xy^{-1}):x,y\in X^{[4]}\}|<K|\tr(X)|.}
for some $1\leq K\in \R$ and $a_1,a_2 \in \F^{\times}$. 

Then we have, 
\doeq{\l{prop.3.18.1.then}|\tr(X^2)\tr(X^2)|+	|\tr(X^2)+\tr(X^2)|<CK^C|\tr(X)|.}

Let $V\subseteq \SL_2(\F)$ be a finite subset of diagonal matrices and suppose
\doeq{\l{trace growth}|\Tr(V^{[4]} \ddot V^{g[4]})|< |\Tr(V)|^{1+\eps}}
for some $g\in \SL_2(\F)$ with\footnote{i.e., $g$ has no zero entries.} $Prod(g)\notin \set{0,1}$ and some $\eps\in \R_+$. 

Then we have,
\doeq{\l{prop.3.18.2.then}|\Tr(V^2)\ddot \Tr(V^2)|+|\Tr(V^2)+\Tr(V^2)| <C|\Tr(V^2)|^{1+C\eps}.}
}

\doPROOF{
Denote $N:=|\tr(X)|$ and for $x,y\in\F^{\times}$ denote
$$\tr_{(a_1,a_2)}(x,y):=a_1\ddot \tr(xy)+a_2\ddot \tr(xy^{-1}).$$
\BTA\  \re{prop.3.18.1.if} we get, 
\doAlignSplit{bg proof 1}{
			|\tr_{(a_1,a_2)}(X^{[4]},X^{[4]})| 	&= |\tr_{(1,a_1/a_2)}(X^{[4]},X^{[4]})|\\
																					&< K|\tr(X)|\\ 
																					&= KN.
}
%
%

Now for any subset $Y\subseteq \F^{\times}$ we get the following. For any $z,w\in Y$ we have 
$x:=zw,y:=zw^{-1} \in Y^{[2]}$ and $t:=xy=z^2,s:=xy^{-1}=w^2 \in Y^2$. 
Therefore we get,
\doeq{\l{bg trick}\{(t,s):t,s\in Y^2\}\subseteq \{(xy,xy^{-1}):x,y\in Y^{[2]}\}.}
Now set $Y:=X^{[2]}$ which satisfy 
\doeqnn{X^{[2]2}=Y^2\subseteq Y^{(2)}=Y^{[2]}=X^{[4]}.}
Therefore by  \re{bg proof 1} and  \re{bg trick} we get,
$$|\{\tr(t)+a\ddot \tr(s):t,s\in Y^2\}|<KN.$$
Denote $Z:=Y^2=X^{[2]2}$ so we got 
$$|\tr(Z)+ a \ddot \tr(Z)|<KN.$$
Since $\mult(\tr(x^2))\leq 4$, we have 
$$N=|\tr(X)|\leq |X|\leq |Y^2|\leq 4|\tr(Y^2)|= 4|\tr(Z)|.$$ 
Therefore 
\doeq{\l{bound trx22}
	|\tr(Z)|\leq |\tr(Z)+a\ddot \tr(Z)|<KN\leq 4K|\tr(Z)|.
}

Therefore, by \Plunnecke-Ruzsa  \re{plunnecke from twisted sum to sum and minus} with 
$$A=B=\tr(Z)=\tr(X^{[2]2})$$ we get
\doAlignSplit{plus bound}{
		|\tr(X^2)+ \tr(X^2)|	&\leq |\tr(X^{[2]2})+ \tr(X^{[2]2})|\\
													&< 4^2 K^2 |\tr(X^{[2]2})|\\ 
													&= 2^4 K^2 |\tr(Z)|
}

Now by fact \r{fact:trace identity for scalars} applied to $W=X^2$ we get that
$$|\tr(X^2)\ddot \tr(X^2)| \leBY{fact:eq:prod_sum_bound} |\tr(X^{2[2]})+\tr(X^{2[2]})|.$$
But since $X\subseteq \F$ we have $Z=X^{2[2]}=X^{[2]2}$ we get by  \re{plus bound} that
\doAlignSplit{bound prod}{
		|\tr(X^2)\ddot \tr(X^2)| &\leq 	|\tr(X^{2[2]})+\tr(X^{2[2]})|\\
														 &= 		|\tr(X^{[2]2})+\tr(X^{[2]2})|\\
														 &\leq 	2^4 K^2 |\tr(Z)|.
}
Therefore by combing  \re{bound prod} and  \re{plus bound} we get
\doearray{
	|\tr(X^2)+\tr(X^2)|+|\tr(X^2)\ddot \tr(X^2)| &\le& 2^5 K^2 |\tr(Z)|\\ 
																												 &\lnBY{bound trx22}& 2^5 K^3 N\\
																												 &=& 2^5 K^3 \tr(X)
}
\sWAD\s with  \re{prop.3.18.1.then}.\good 

Set $X:=\set{x\in \F: D_x \in V}$ (i.e., $V=D_X$). 
\BTA\  \re{trace growth} and by fact \r{sfact:tr_g(x,y)} 
 we get
\doearray{
|\{ad\ddot \tr(xy)-bc\ddot \tr(xy^{-1}):x,y\in X^{[4]}\}| &\eqBY{eq:tr_g(x,y)}& |\Tr(V^{[4]}\ddot V^{g[4]})|\\
																													&\lnBY{trace growth}& |\Tr(V)|^{1+\eps}\\
																													&=&	|\tr(X)|^{1+\eps}
} 
Therefore by  \re{prop.3.18.1.then}, we have,
\doeqnn{|\tr(X^2)\tr(X^2)|+|\tr(X^2)+\tr(X^2)|\ll |\tr(X)|^{1+O(\eps)}.}
\IOW\ we have,
\doeqnn{|\Tr(V^2)\Tr(V^2)|+|\Tr(V^2)+\Tr(V^2)|\ll |\Tr(V)|^{1+O(\eps)}} 
\sWAD\s with  \re{prop.3.18.2.then}.\good 
}

Now let us see some very simple observations that we will use later.

\doLEM{\l{lem:expansion of products}
There exists $c>0$ such that the following holds. 
Let $\F$ be a field and let $g \in \SL_2(\F)$.
Let $V\subseteq \SL_2(\F)$ be a finite subset of diagonal matrices. 
Suppose $\Tr(V^{[4]})\subseteq \E$ for some subfield $\E\leq \F$.	

If $\Prod(g) \notin \E$ then we have,
\doeq{|\Tr(V^{[4]} \ddot V^{[4]g})|>c|\Tr(V)|^2. \l{expansion of fields}}	
If $\Prod(g) \neq 1$ then we have,
\doeq{|\Tr([V,g])|> c|\Tr(V)|. \l{expansion of commutator}}
}

\doPROOF{
Denote $g=\dobmatrix{a & b\\ c & d}$. Set $X:=\set{x\in \F: D_x\in V}$ and set 
\doalign{
	T	&\eqB	\Tr(V^{[4]}V^{[4]g})\\
		&\eqBY{eq:tr_g(x,y)}\{ad \ddot \tr(xy)-bc\ddot \tr(x/y):x,y\in X^{[4]}\}.
	}
Therefore we get, 
\doeqnn{T':=\{ad\ddot \tr(t)-bc \ddot \tr(s):t,s\in X^{[2]2}\} \byeq{bg trick}{\subseteq} T.} 
Set $f(z,w):= ad \ddot z+(1-ad)\ddot w$ and since
$ad-bc=1$ we get 
\doeqnn{T'=f(\tr(X^{[2]2}),\tr(X^{[2]2})).} 
Now if $\Prod(g)=ad\notin \E$ then $f|_{\E\times\E}$ is injective. 
Indeed: if we set $t=ad$ then by solving $tz+(1-t)w=tz'+(1-t)w'$,
we get that $t(z-z')=(1-t)(w'-w)$. 
Since $t\neq 0,1$ we get that either $z-z'=w'-w=0$ or $\frac{1-t}{t}=t^{-1}-1\in \E$ which contradicts our assumption that that $t=ad\notin\E$. 
Note that by the same way $f\MID{x\E\times x\E}$ is injective for any coset of $\E$.
\BTA\ 
$$\tr(X^{[2]2})\subseteq \tr(X^{[4]})=\Tr(V^{[4]})\subseteq \E$$ therefore
\doeqnn{|T|\geq |T'| = |\tr(X^{[2]2})|^2 \geq |\tr(X^2)|^2 \geq (\fr{4}|\tr(X)|)^2} 
\sWAD\s with  \re{expansion of fields}.\good 

Now if $\Prod(g)=ad\neq 1$ then we get by fact \r{sfact:tr_g(x,y)} that,
\doalign{
	|\Tr([V^{[4]},g])| &\eqB	|\set{\Tr(v^{-1} v^g):v\in V^{[4]}}|\\
							&\eqBY{eq:tr_g(x,y)}	|\set{2ad+(1-ad)\tr(x^2):x\in X^{[4]}}|\\
							&\eqB |\tr(X^{[4]2})|\\
							&\geB \fr{4}|X^{[4]}|
} 
\sWAD\ with  \re{expansion of commutator}.\good  
}

\doSFACT{ 
Let $V$ and $g$ be as in Lemma \r{lem:expansion of products} and let $x,y\in \F^{\times}$. 
\TWH, 
\doEQsplit{tr(xy)=tr(x/y)}{
	\tr(xy)=\tr(x/y) &\iff \tx{either }x^2=1\tx{ or }y^2=1\\ 
									 &\iff \tx{either }x=\pm 1\tx{ or }y=\pm 1.
}

If $\Tr(V^{[2]})\subseteq \E$ and $V\nsubseteq \set{\pm I}$ then, 
\doeq{\Prod(g)\in \E \iff \Tr(VV^g)\subseteq \E.\l{rem:Tr(VV^g) subseteq E}}
}

\doPROOF{
Note that $$\tr(x)=2 \iff x=1$$ and $\tr(x)=-2 \iff x=-1$. 
Moreover for any $c\neq \pm 2$, 
$$|\tr^{-1}(c)|=2$$ since $\tr(x)=\tr(x^{-1})$ and $x\neq x^{-1}$. 
Therefore 
\doalign{
	\tr(x)=\tr(y) &\iff x\in\set{y^{\pm 1}}\\
								&\iff xy=1\qu\tx{or}\qu x/y=1.
}
\IP\  \re{tr(xy)=tr(x/y)} follows.\good

By fact \r{sfact:tr_g(x,y)} we get that 
$$\Tr(D_x D_y^g) \eqBY{eq:tr_g(x,y)} \Prod(g)(\tr(xy)-\tr(x/y))+\tr(x/y).$$  
Therefore if $D_x\neq \pm I$ and $D_y\neq \pm I$ and $\tr(xy),\tr(x/y) \in \E$ then 
$$\Tr(D_x D_y^g)\in \E \iff \Prod(g)\in \E.$$
\TF\ we immediately get \re{rem:Tr(VV^g) subseteq E}.\good
}
 

\sec{Useful properties of $\SL_2(\F)$}\l{sec:Useful properties of SL_2(F)}

\ssec{Bounded generation of large subsets}

In the following section, we will prove few Growth properties of large subsets of finite (quasi) simple groups.
First we give some background concerning the regular representation (and the convolution of functions). 
We will follow the techniques which were developed by Gowers (cf. \cite{gowers}) and later were expanded by Babai, Nikolov and Pyber (cf. \cite{BNP,Nikolov-Pyber}). 
%

\sssec{The spectral decomposition} 


\doDEF{
Let $G$ be a finite group. 
We identify the group ring $\C[G]$ with $\C^G$ so instead of writing $\sum a_g g \in \C[G]$ 
we write $X\in \C^G$ with $X(g)=a_g$ for all $g\in G$.
 
On the other hand we identify subsets $A\subseteq G$ as the indicators functions 
$1_A\in \C^G$ and similarly elements $g\in G$ as the indicators functions $1_g\in \C^G$. 

In the algebra $\C[G]$ we have the usual inner product and convolution product.
For $X,Y\in \C[G]$ we have 
$$\gen{X,Y}=\sum_g X(g)\overline{Y(g)}$$ 
and the (convolution) product $\mbf{X* Y}$, or for short just $\mbf{X.Y}$, is defined by, 
\doeqnn{(X.Y)(g)=\gen{X.Y,g}=\sum_{xy=g} X(x)Y(y).}
}

\doSFACT{\l{fact:conjugate}
Let $G$ be a finite group and let $X,Y,Z\in \C[G]$.
Define $X^T\in \C[G]$ by $$\mbf{X^T}(x):=X(x^{-1})$$ and 
$X^*\in \C[G]$ by 
$$\mbf{X^*}(x)=\overline{X(x^{-1})}.$$
We will be interested mainly in functions in $\R[G]$ so there will be no difference in these notations.
Then we have,
$$\gen{X.Y,Z}=\gen{Y,X^*Z}=\gen{X,ZY^*}.$$
}

\doPROOF{
For any $x,y,z\in G$ we have,
$$\gen{xy,z}=\gen{y,x^{-1}z}=\gen{x,zy^{-1}}$$
therefore by linearity we get,
$$\gen{X.Y,Z}=\gen{Y,X^*.Z}=\gen{X,Z.Y^*}.$$
}

\doDEF{\l{def:L(),R()}
Let $G$ be a finite group and let $X,Y\in \C[G]$.
Let $L(\cdot)$ and $R(\cdot)$ be the left and the right regular representations of $G$, 
\doearray{
	L(X)(Y)	&:=&	X.Y\\
	R(X)(Y)	&:=& Y.X^*.
} 
}

\doSFACT{
Let $G$ be a finite group, $X,Y\in \C[G]$ and let $L(\cdot)$ and $R(\cdot)$ be the left and the right regular representations of $G$.
Then we have,
\doearray{
	L(X.Y)		&=& L(X)L(Y)\\
	R(X.Y)		&=& R(X)R(Y).
}
Clearly $L(\cdot)$ and $R(\cdot)$ commutes,
\doeqnn{L(X)R(Y)	=R(Y)L(X).}
Moreover, we have
\doearray{
	L(X)^*	&=& L(X^*)\\	
	R(X)^*	&=& R(X^*).
}
i.e., $$\gen{X.u,v}=\gen{L(X)u,v}=\gen{u,L(X)^*v}=\gen{u,X^*.v}$$ for any $u,v\in\C[G]$ (and similarly for $R(X)$).
}
\doPROOF{
All follows immediately from the definitions of $L()$ and $R()$ in \r{def:L(),R()} and fact \r{fact:conjugate}.
}

\doSFACT{
Let $V=\C[G]$ and denote by $U\lr{V}$ the group of unitary transformations of $V$.
Then $L(G)$ and $R(G)$, the left and the right regular representations, and also $X\mapsto X^T$, 
are all in $U(\C[G])$.
}
\doPROOF{
Clearly for any $g\in G$, $L(g)$ and $R(g)$ and $X\mapsto X^T$, are linear maps which permute the \OB\ $\set{h:h\in G}$.
}

\doSFACT{\l{fact:X(1)}
Let $G$ be a finite group of size $N$ and $X\in \C[G]$.
Then we have,
\doeqnn{X(1)=\tfr{N}\Tr(L(X))=\tfr{N}\Tr(R(X^*)).}
}


\doPROOF{
For any $g\in G$ we have, 
\doearray{
	X(1) &=& \gen{X.1,1}\\
			 &=& \gen{X.g,g}\\
			 &=& \gen{g.X,g}
}			 
Therefore $L(X)$ and $R(X^*)$ have the same diagonal \WRT\ the \OB\ $\set{g:g\in G}$.  
}


\doDEF{\l{def:Tr(X)}
Let $G$ be a finite group of size $N$ and let $X\in \C[G]$. Denote,
$$\mbf{\Tr(X)}:=\Tr(L(X)).$$ 
Therefore by \r{fact:X(1)} we get,
\doeqnn{X(1)=\tfr{N}\Tr(X).}
}

\doSFACT{
Let $G$ be a finite group of size $N$ and $X\in \C[G]$. Then,
$$\|X\|^2=X^*.X(1)=X.X^*(1)$$
and 
$$\|X\|^2=\tfr{N}\Tr(X^*.X)=\tfr{N}\Tr(X.X^*).$$
}

\doPROOF{
We have,
\doearray{
	\|X\|^2	&=& \gen{X.1,X.1}\\
					&=& X^*.X(1)\\
					&=& X.X^*(1)
}
Therefore by \r{def:Tr(X)} with $\mbold{X}=X.X^*$ (and $\mbold{X}=X^*.X$) we are done.
}

\doTHM{[SD\footnote{The spectral decomposition} of real symmetric endomorphism]\l{thm:SD:A}
Let $G$ be a finite group of size $N$ and let $A \in \End(\R[G])$.
Suppose $A$ is be a symmetric endomorphism i.e.\footnote{see fact \r{fact:conjugate}.}, $A=A^T$. 
Then there exist an \OB\ $\ol{\alpha}=(\alpha_i)$ 
of $\R[G]$,
and $\lam_1\ge \lam_2 \ldots \ge \lam_N$ in $\R$ \ST
\doeq{\l{eq:SD:A}\gen{A\alpha_i,\alpha_j}=\delta_{ij}\lam_i}
for any $1\leq i,j\le N$.
}

\doPROOF{
This is a standard theorem in linear algebra for symmetric matrix $T\in \M_n(\R)$.
}

\doCORO{
[Rayley inequality]\l{coro:Rayley}
Let $G$ be a finite group of size $N$ and let $A\in \End(\R[G])$ (not necessarily symmetric).
Then there exist 
\OB\ $\ol{\beta}$ of $\R[G]$, and $\lam_1\ge \lam_2 \ldots \ge \lam_N \ge 0$ in $\R$ \ST
\doeArray{\l{eq:SD:A^TA}\gen{A\beta_i,A\beta_j}&=&\delta_{ij}\lam_i^2}
for any $1\leq i,j\le N$.
Let $1\le k\le N$ and suppose $v\in \C[G]$ with $v\bot \beta_i$ for all $i<k$. Then we have,
\doeq{\l{eq:Rayley}\n{Av}\le \lam_k\n{v}.}
}

\doPROOF{
Since $AA^T,A^TA\in \End(\R[G])$ are symmetric we can decompose $A^TA$ and $AA^T$ by theorem \r{thm:SD:A}.
Moreover $AA^T,A^TA\geq 0$ (i.e., they are positive-semidefinite) therefore they have the same, non negative, eigen values.
Therefore there exist 
\OB\ $\ol{\beta}$ of $\R[G]$, and $\lam_1\ge \lam_2 \ldots \ge \lam_N \ge 0$ in $\R$ \ST
\doearray{
	\gen{A^TA\beta_i,\beta_j}	&=&	\delta_{ij}\lam_i^2
}
for any $1\leq i,j\le N$. 
\SWADW  \re{eq:SD:A^TA}.\good

Let $1\le k\le N$ and $v\in \C[G]$ and suppose $v\bot \beta_i$ for all $i<k$. 
Then we have,
\doearray{
	\n{Av}^2 &=& \gen{Av,Av}\\
					 &=& \gen{A^TAv,v}\\
					 &=& \lrg{\sum_i\gen{v,\beta_i}A^TA\beta_i,\sum_j\gen{v,\beta_j}\beta_j}\\
					 &=& \sum_{1\le i,j\le N} \gen{A^TA\beta_i,\beta_j}\gen{v,\beta_i}\ol{\gen{v,\beta_j}}\\
					 &\eqBY{eq:SD:A^TA}& \sum_{i=k}^N \lam_i^2|\gen{v,\beta_i}|^2\\
					 &\le& \lam_k^2\n{v}^2
}
\sWAD.
}
 
\doDEF{\l{def:lam_i & m_i}
Let $G$ be a finite group of size $N$ and let $A\in \End(\R[G])$.
By corollary \r{coro:Rayley} 
there exist \OB\ $\ol{\beta}$ of $\R[G]$,
and $0\leq \lam_i \in \R$, in decreasing order, s.t.
\doearray{
	\gen{A^TA\beta_i,\beta_j}&=&\delta_{ij}\lam_i^2.
}
Denote $\lam_i(A):=\lam_i$ and by $m_i(A)$ the multiplicity of $\lam_i(A)$.
I.e., 
$$\mbf{m_i(A)}:=\dim(\Ker(A^TA-\lam_i^2\Id)).$$
Denote $\mbf{\lam(X)}:=\lam_2(X)$ and $\mbf{m(X)}:=m_2(X)$.
}



\sssec{Rapid mixing and Mixing Growth}

\doDEF{\l{def:Fix(G)}
Let $G$ be a group and let $\F$ be a field and let $(\rho,V)$ be finite dimensional representation of $G$.
Denote the fix points of $(\rho,V)$ by 
\doalign{
	\mbf{\Fix(\rho(G))}	&:=\set{v\in V:\rho(g)v=v\tx{ for any }g\in G}\\
											&:= \CAP[g\in G] \Ker(\rho(g)-\Id)
}	
and if the action is clear from the context we will abbreviate and write $\mbf{\Fix(G)}$.
We will say that $(\rho,V)$ is a \tb{trivial} representation if $$\Fix(G)=V.$$
}

\doDEF{\l{def:M(G)}
Let $G$ be a finite group and let $\F$ be a field. 
Define
$$M(G,\F):=\min\set{deg(\rho):\rho \tx{ is a non-trivial irreducible }\F\tx{-representation of }G}.$$
Since $M(G,\C)$ and $M(G,\R)$ will be more relevant for our purposes when investigating finite groups, we abbreviate 
$$\mbf{M(G):}=M(G,\R)$$ the minimal degree of non-trivial \ti{real} representation of it. 
}

\doDEF{\l{def:Prob[G]}
Denote by $\mbf{\Prob[G]}$ the elements $X\in \R[G]$ with $X(g)\geq 0$ for any $g\in G$ and with $\Lo{X}=1$.
Denote by $\mbf{U_X}$ the uniform probability on the support of $X$.
I.e., if $A=\supp(X)$ then $U_X=\tfr{|A|}1_A$.
Denote by $\mbf{U}=U_G\equiv \tfr{N}$ the uniform probability on $G$. 
}

\doSFACT{\l{sfact:Lt(Y-U):ge}
Let $G$ be a group of size $N$ and let $Y\in \Prob[G]$.
Then we have,
\doalign{
	\|Y-U\|^2 &\eqB \|Y\|^2 - \tfr{N}. 	
	}
\IP\ $$\|Y\|^2 \ge \tfr{N}$$ with equality if and only if $Y=U$.
Moreover $$\|Y\|^2 \ge \fr{|\supp(Y)|}$$ with equality if and only if $Y=U_Y$.
}

\doPROOF{
Since $Y-U \bot U$ and $Y-U_Y\bot U_Y$ we get
\doalign{
	\|Y\|^2 &\eqB \|Y-U\|^2 + \|U\|^2\\
					&\eqB \|Y-U_Y\|^2 + \|U_Y\|^2
}
therefore the claim follows.
}

\doPROP{[Young inequality]\l{prop:Young}
Let $1\le r,p,q \le \infty$ and suppose $\tfr{p} + \tfr{q} = 1+ \tfr{r}$.
Let $G$ be a finite group and let $X,Y\in \C[G]$.
Then we have,
\doeq{\l{eq:Young}\|X.Y\|_r\le \|X\|_p \|Y\|_q.}
We will call such a triple $(r,p,q)$ a \tb{Young triple}.
}

\doDEF{\l{def:operator norm}
Let $G$ be a finite group and $A\in \End(\C[G])$. For any $p,q\ge 1$ denote the \tb{operator norm} $\mbf{\|A\|_{p,q}}$ by 
\doalign{
	\|A\|_{p,q}	&=	\max_{v\neq 0}\frac{\|A(v)\|_p}{\|v\|_q}\\ 
							&=\max_{\|v\|_q=1} \|A(v)\|_p
}
Denote by $\mbf{\Lam(A)}$ the \tb{spectrum} of $A$ and by $\mbf{\rho(A)}$ the \tb{spectral radius} of $A$.
}

\doSFACT{\l{def:rho(X)}
For any $X\in \R[G]$, the operators 
$$L(X),L(X^T),L(X^*),R(X),R(X^T),R(X^*)$$ have the same spectral radius, the same spectrum and the same operator norms.
Therefore we write for short 
$$\mbf{\|X\|_{p,q}}=\|L(X)\|_{p,q}$$ and 
$$\mbf{\rho(X)}=\rho(L(X))$$ and the spectrum of $L(X)$ by
$$\mbf{\Lam(X)}=\Lam(L(X)).$$
\IP, for any Young triple $(r,p,q)$ and $X\in \C[G]$ we get
$$\rho(X)\leq \|X\|_{r,p}\leq \|X\|_q.$$
}

\doSFACT{\l{sfact:lam:ge}
Let $G$ be a group of size $N$ and let $X,Y\in \Prob[G]$. Then we have,
\doAlignSplit{eq:lam:ge}{	
	\|X.Y- U\| &\le \lam(X)\lrn{Y- U}\\ 
	\|X.Y- U\| &\le \lam(Y)\lrn{X- U} 
	}
}

%


\doPROOF{

\FTOH\ $$U.g=g.U=U$$ for any $g\in G$ so we get that $$X.U=U.X=U$$
and so $\lam_1(X) \ge 1$.

\OTOH\ by Young inequality \r{prop:Young} with $(r,p,q)=(2,1,2)$ we get,
$\|X.Y\|\le \|Y\|$ for any $Y\in \C[G]$ i.e.,
$$\rho(X) \le \|X\|_{2,2}\le \n{X}_1=1,$$
therefore $$\lam_1(X)=1,$$
and by corollary \r{coro:Rayley}, \wLOG\ we can assume $\beta_1\equiv \tfr{\sqr{N}}$.
Now since $X \bot Y - U$, we get by Rayley inequality \r{coro:Rayley} that
\doearray{
	\|X.Y-U\| &=& \lrn{X.(Y-U)}\\
						&\leBY{eq:Rayley} & \lam_2(X)\lrn{Y-U}
}
\sWADW with the first inequality of  \re{eq:lam:ge}.

Now since $$\n{X.Y-U}=\|Y^T\!\!.X^T-U\|,$$ we can apply the first bound with $\mbf{X}=Y^T$ and $\mbf{Y}=X^T$, 
\sWAD\ with the second inequality of  \re{eq:lam:ge}.
}

%

\doSFACT{\l{sfact:Fix(G)}
Let $G$ be a finite group. Then
$$\Fix(L(G))=\Fix(R(G))=\Span(U).$$
\IOW\ for any linear subspace $0\neq W\le \C[G]$ we have,
\doalign{
				&	G.W\subseteq W\\
	\iffB &W.G\subseteq W\\
	\iffB &W =\Span(U).}
}

\doPROOF{
Since for any $g,h\in G$ we have
\doalign{
	X(hg^{-1})	&\eqB	\gen{X,hg^{-1}}\\
							&\eqB	g.X(h)\\
							&\eqB	X.g(h)
}							
we get that 
\doalign{
	X\in \Fix(G) &\iffB X\equiv X(1).\QED
	}
}
%
%

\doSFACT{\l{sfact:m_i:geq:M}
Let $G$ be a finite group of size $N$ with $M=M(G)$ and let $X\in \Prob[G]$. 
Then for any $1<i\le N$ we have, 
$$m_i(X)\ge M.$$
}


\doPROOF{
Set $A:=L(X^T\!\!.X)$ and for any $1\le i \le N$ set
\doalign{
	A_i	&:=A-\lam_i^2\Id\\
	V_i	&:=\Ker(A_i).
}
Since $A_i$ commutes with all the elements of $R(G)$ and $A_i \in \End(\R[G])$ we get that $V_i$ is a real representation of $G$ (with the right action of $G$ on $V_i$).
If $i\neq 1$ then $\beta_i \bot U$ and so by \r{sfact:Fix(G)},
$V_i$ is non trivial real representation of $G$ so 
\doeqnnQED{m_i=\dim(V_i)\ge M.}
}


\addREM{
Note that there is no a priori assumption that $\lam_2(X)\neq 1$.
Actually if $\lam_2(X)=1$ then by the same argument we get that 
$$m(X)=m_2(X)\ge M+1.$$
}

\doSFACT{\l{sfact:lam:le}
Let $G$ be a group of size $N$ and let $M=M(G)$.
Then for any $Y\in \Prob[G]$ we have,
\doeq{\l{eq:lam:le}\lam(Y)\leq \sqr{\tfr[N]{M}}\lrn{Y-U}.}
}


\doPROOF{
Since $m_2(Y)\geq M$ and $Y - U \bot U$ we get
\doalign{
	\|Y-U\|^2 &= 		\|Y\|^2-\|U\|^2\\
						&=		 \tfr{N}(\Tr(Y^T.Y)-1)\\
						&= 		\tfr{N}\sum_{i=2}^N\lam_i^2(Y)\\
						&\ge 	\tfr[M]{N}\lam^2(Y). \QED
	}
}

\sssec{Corollaries}

The following Corollary is a slight modification of an argument of \cite{Nikolov-Pyber,BNP,BNP_new} (which followed and extended results of \ci{gowers}).


\doCORO{\CITE{Theorems 2.1 and Corollary 2.2}{BNP_new}\l{coro:Lt(X.Y-U):le} 
Let $G$ be a group of size $N$ with $M=M(G)$ and let $X,Y\in \Prob[G]$.
Then we have,
\doeq{\l{eq:Lt(X.Y-U):le}
	\|X.Y-U\|\le\sqr{\tfr[N]{M}}\|Y-U\|\|X-U\|.
	}
Inductively we get for any $n\in \N$ and $X_1,\ldots,X_{n+1} \in \Prob[G]$,
\doeq{\l{eq:Lt(X.Y-U):le:many} 
	\|X_1\ldots X_{n+1}-U\|\le (\tfr[N]{M})^{n/2} \prod_{i=1}^{n+1} \|X_i-U\|.
	}
}
\doPROOF{
By facts \r{sfact:lam:ge} and \r{sfact:lam:le} we get,
\doalign{
	\|X.Y-U\|^2 &\leBY{eq:lam:ge} \lam(X)\lam(Y)\|Y-U\|\|X-U\|\\
							&\leBY{eq:lam:le} \tfr[N]{M}\|Y-U\|^2\|X-U\|^2 \QED
}
}
%
%
%



\doCORO{\CITE{Corollary 2.3}{BNP_new}\l{coro:L_i(X.Y.Z-U):ln} 
Let $G$ be a group of size $N$ with $M=M(G)$ and let $X,Y,Z\in \Prob[G]$.
Then we have,
\doeq{\l{eq:L_i(X.Y.Z-U):ln} 
	\Li{X.Y.Z-U}<\sqr{\tfr[N]{M}}\|X\|\|Y\|\|Z\|.
	}
Inductively we get for any $n\in \N$ and $X_1,\ldots,X_{n+2}\in \Prob[G]$ that,
\doeq{\l{eq:L_i(X.Y.Z-U):ln:many} 
	\Li{X_1\ldots X_{n+2}-U} 	< (\tfr[N]{M})^{n/2} \prod_{i=1}^{n+2} \|X_i\|
	}	
}

\doPROOF{
By proposition \r{prop:Young} with $(r,p,q)=(\infty,2,2)$ and corollary \r{coro:Lt(X.Y-U):le} we get,
\doalign{
	\Li{X.Y.Z-U}	&\eqB  \Li{(X.Y-U).Z}\\ 
								&\leBY{eq:Young}	\n{X.Y-U}\n{Z}\\
								&\leBY{eq:Lt(X.Y-U):le}	\sqr{\tfr[N]{M}}\|X-U\|\|Y-U\|\|Z\|\\
								&\lnB	\sqr{\tfr[N]{M}}\|X\|\|Y\|\|Z\|.	\QED
	}
}




Now let us the implications of the properties above \CITE[cf. ]{Corollaries 2.5 and 2.6 and Theorem 2.14}{BNP_new}.

\doTHM{[Babai-Nikolov-Pyber] \l{thm:bnp:subsets}
Let $G$ be a finite group of size $N$ with $M=M(G)$. 
Let $A_1,\ldots,A_t\subseteq G$ be subsets of size $|A_i|=K_i\frac{N}{M}$ where $K_i\in \R_+$.
\TWH,
\doeq{|A_1 A_2|	>	\fr{2}\MIN{K_1 K_2\tfr[N]{M},N} \l{eq:bnp:subsets:two}}
and if $t\geq 3$ then we have\footnote{The case $t=3$ was proved in \cite[Theorem 3.3]{gowers}.},
\doeq{\prod_{i=1}^t K_i \geq M^2 \RightarrowLB \prod_{i=1}^t A_i=G.\l{eq:bnp:subsets:many}}
}



\doPROOF{
For any $1\le i\le t$ set $X_i \in \Prob[G]$ by
$$X_i:=U_{A_i}=\tfr{|A_i|} 1_{A_i}.$$
Since\footnote{One can denote the convolution either as $X_1*\ldots *X_t$ or $X_1\ldots X_t$ or just by $X_1\cdots X_t$ since this is the product in the algebra $\C[G]$. We use in this manuscript the middle way.} 
$$\supp(X_1\ldots X_t)=\supp(X_1)\cdots \supp(X_t),$$ 
we get by corollary \r{coro:Lt(X.Y-U):le} that,
\doalign{
	\fr{|A_1 A_2|} 	&\eqB \fr{|\supp(X_1.X_2)|}\\
									&\leB \n{X_1.X_2}^2\\
									&\eqB \n{X_1.X_2-U}^2 + \n{U}^2\\									
									&\leBY{eq:Lt(X.Y-U):le} \tfr[N]{M}\|X_1-U\|^2\|X_2-U\|^2 + \tfr{N}\\
									&\lnB \tfr[N]{M}\|X_1\|^2\|X_2\|^2 + \tfr{N}\\
									&\eqB \fr[N]{M}\fr{|A_1||A_2|} + \fr{N}\\
									&\eqB \fr[M]{N}\fr{K_1 K_2} + \fr{N}\\
									&\eqB \fr[M]{N}(\fr{K_1 K_2} + \fr{M}).\\
									&\leB 2\fr[M]{N}\MAX{\fr{K_1 K_2},\fr{M}}
}
Therefore by rearranging the inequalities we are done with  \re{eq:bnp:subsets:two}.

Now by corollary \r{coro:L_i(X.Y.Z-U):ln} we get,
\doalign{
	\Li{X_1.X_2.X_3-U}	&\lnBY{eq:L_i(X.Y.Z-U):ln}	\sqr{\tfr[N]{M}}\|X_1\|\|X_2\|\|X_3\|\\
											&\eqB	\sqr{\tfr[N]{M}}(|A_1||A_2||A_3|)^{-1/2}\\
											&\eqB	\tfr[M]{N} (K_1 K_2 K_3)^{-1/2}.
	}
Therefore if $K_1 K_2 K_3 \ge M^2$ then	$\Li{X_1.X_2.X_3-U}<\tfr{N}$ so $$A_1\ddot A_2\ddot A_3=\supp(X_1.X_2.X_3)=G$$
\sWADW  \re{eq:bnp:subsets:many} for $t=3$.

Similarly by corollary \r{coro:L_i(X.Y.Z-U):ln} we get for any $n\in \N$,
\doalign{
	\Li{X_1\ldots X_{n+2}-U} 	&\lnBY{eq:L_i(X.Y.Z-U):ln:many} (\fr[N]{M})^{n/2} \prod_{i=1}^{n+2} \|X_i\|\\
														&\eqB	\fr[M]{N} (\prod_{i=1}^{n+2}K_i)^{-1/2}
	}
\sWADW  \re{eq:bnp:subsets:many} for $t\ge 3$.	
}


As a special case of the previous we get immediately the following Corollary.

\doCORO{\CITE{Corollary 2.11}{BNP_new}\l{coro:bnp}	
Let $G$ be a finite group of size $N$ with $M=M(G)$. 
Let $A \subseteq G$ be a subset of size $|A|=K\frac{N}{M}$.
Then we have, 
\doeqnn{|A^{(2)}|	\gnB \tfr{2}\MIN{N,K|A|}.}
And for any $t\ge 3$ we have,
$$|K|\ge M^{\dfr[2]{t}}	\RightarrowLB A^{(t)}=G.$$
}

%
%

%
%
%

\doTHM{\l{BNP} 

There exist $C\in \R_+$ such that the following holds.
Let $\F_q$ be a finite field and let $A$ be a subset of $G=\SL_2(\F_q)$. 
Then we have,
\doeq{\l{BNP_bound_1a}|A| > Cq^{2\frac{2}{3}} \Rightarrow A^{(3)}=\SL_2(\F_q).} 
For any $3\leq m\in \N$ we have,
\doeq{\l{BNP_bound_1b}|A|> Cq^{2\frac{2}{m}} \Rightarrow A^{(m)}=\SL_2(\F_q).}
For any $0<\delta \leq \tfr{2}$ we have,
\doeq{\l{BNP_bound_1c}|A|> q^{2+\delta} \Rightarrow |A^{(2)}|> \fr{C}q^{2+2\delta}.}
}

\doPROOF{
By a well known fact (which was first proved by Frobenius) for any finite field $\F_q$ and $G=\SL_2(\F_q)$ we have,
\doalign{
	M(G,\R)	&\eqB \tfr{2}(q-1)\\
					&\geB \tfr{2}q(1-o(1))\\
					&\ggB q.  
	}
Therefore if $N=|G|=q(q^2-1)$ and $M=M(G)$ then
\doalign{
	\fr[N]{M}					&\eqB 2q(q+1)\\
										&\leB 2q^2(1+o(1))\\
										&\llB q^2\\ 
\addtx{and for any $m\ge 3$,}
	\fr[N]{M^{1-2/m}}	&\llB q^{2\frac{2}{m}}. 
	}
Therefore the claim follows immediately by corollary \r{coro:bnp}. 
}


%
%

\addREM{
\IP\ the theorem guarantee bounded generation for \ti{any} large subset $A$ of $G=\SL_2(\F_q)$. 
\IP, any subgroup $H<G$ has large index 
$$[G:H]\gg q.$$
}

\ssec{Symbolic generation of traces} \label{invariant theory}

The Invariant theory of tuples of matrices under various actions
was developed over fields of zero characteristic. We will actually be interested in 
the positive \char\ (cf. \cite{Procesi: Invariant theory of n by n matrices},
\cite{Concini-Procesi: A characteristic free approach to invariant theory},
\cite{Donkin: Invariants of several matrices}). 
 
\doDEF{ 
For $m\geq 2$ denote by $\op{R}_{2,m}$ the ring of
invariants of $m$-tuples of $2\times 2$ generic matrices $(X_1,\ldots,X_m)$ over a infinite
field $\F$ under the simultaneous conjugation action of the general linear group.
To be precise, we have $4m$ variables $x_1,y_1,z_1,w_1,\ldots,x_m,y_m,z_m,w_m$ which we denote 
by $\ol{X_i}=(x_i,y_i,z_i,w_i)$ and $\ol{X}=(\ol{X_1},\ldots,\ol{X_m})$.
Each matrix $X_i=\dobmatrix{x_i & y_i\\z_i & w_i}$ is a formal matrix with four variables $\ol{X_i}$ for $1\leq i\leq m$. 
We define an action of $g\in \GL_2(\F)$ on $f(X_1,\ldots,X_m)\in \F[\ol{X}]$ by 
$$f^g(X_1,\ldots,X_m):=f(X^g_1,\ldots,X^g_m).$$
We define the algebra of invariants of this polynomial ring under the action of $\GL_2(\F)$ by 
$$\op{R}_{2,m}(\F):=\set{f\in \F[\ol{X}]: f^g=f\tx{ for any }g\in \GL_2(\F)}.$$
}



We will use the following results of Procesi and Domokos-Kuzmin-Zubkov 
(cf. \cite{Procesi} and \cite[\S 4]{Domokos-Kuzmin-Zubkov}).

\doTHM{\CITE{Corollary 4.1}{Domokos-Kuzmin-Zubkov}

If $\Char(\F)\neq 2$ then,
\doeqnn{\set{det(X_i),\tr(X_{i_1}\cdots X_{i_s}):1\le i\le m;1\le s\le 3; 1\le i_1<\ldots<i_s\le m}} 
is a minimal system of generators of $\op{R}_{2,m}(\F)$. 

If $\Char(\F) = 2$ then, 
\doeqnn{\set{det(X_i),\tr(X_{i_1}\ddot\ldots\ddot X_{i_s}):1\le i,s \le m; 1\le i_1<\ldots<i_s\le m}} 
is a minimal system of generators of $\op{R}_{2,m}(\F)$. 
}

From this we get immediately the following result.
\doLEM{[Trace generation]\label{trace generation}
Let $\F$ be a field and let $A \subseteq \SL_2(\F)$ be a subset of size $2\leq |A| \leq m$.
Then we have the ring generation,
\doeqnn{\gen{\Tr(A^{[m]})}=\gen{\Tr(\gen{A})}.}
Moreover if $\Char(\F)\neq 2$ then we have the ring generation,
\doeqnn{\gen{\Tr(A^{[3]})}=\gen{\Tr(\gen{A})}.}
}

\IP\ we get the following.
\doCORO{
Let $\F$ be a finite field and let $A \subseteq \SL_2(\F)$ be a subset of size $|A| \leq m$.
Suppose $\gen{A}=\SL_2(\F)$. \TWH, 
\doeqnn{\l{trace_A_k}\gen{\Tr(A^{[m]})}=\F} 
and if $\Char(\F)\neq 2$ 
\doeqnn{\l{trace_A_3}\gen{\Tr(A^{[3]})}=\F.}
The same assertion holds under the weaker assumption 
\doeqnn{\gen{\Tr(\gen{A})}=\F.} 
Similarly if $\E$ is a subfield of $\F$ then 
$$\gen{\Tr(\gen{A})}=\E \RightarrowLB \gen{\Tr(A^{[m]})}=\E.$$
}

\addREM{
There are various possible generations types,
depending on the category of objects which are involved: groups,
rings, algebras, vector spaces, modules and fields. In the invariant
context, rings and groups operations are involved. E.g., generation
as $\F_p$ vectors spaces is stronger then rings and for finite
fields there is no difference between ring generation and field
generation. Here the meaning is ring generation in the outer
bracket and group generation in the internal bracket. Explicitly: $\gen{\Tr(A^{[m]})}_{ring}=\gen{\Tr(\gen{A}_{group})}_{ring}$. 
}

\ssec{Size of Minimal generating sets of $\PSL_2(\F_q)$}

By Lemma \r{trace generation} we got that for any finite field $\F=\F_q$ with $\Char(\F)\neq 2$ 
and any subset of generators $\gen{A}=\SL_2(\F_q)$ we have a ``Bounded Generation of Trace Generators'' i.e., 
$$\gen{\Tr(A^{[3]})}=\F.$$ 
In this section we want to extend it to $\Char(\F)=2$ as well. 
The main theorem of this section, and the only part that we will use later, is Theorem \r{Tr(A^{[6]})} which asserts, $$\gen{\Tr(A^{[6]})}=\F.$$

\doDEF{
Let $G$ be a finitely generated group. 
Let us call a subset $A$ of a group $G$ a \tb{minimal generating set} if $\gen{A}=G$ but for any proper subset $A'\subsetneq A$ we have $\gen{A'}\neq G$. 
Let us call a subgroup $H$ of $\PSL_2(\F_q)$ a \tb{subfield subgroup} if $H\cong \PSL_2(q')$ for some subfield $\F_{q'}$ of $\F_q$. 
}

Saxl and Whiston proved the following result about the size of minimal generating sets of $\PSL_2(\F_q)$
\CITE[cf. ]{Theorem 3 and Theorem 7 with its proof}{SW}.

\doTHM{\CITE{Theorems 3,7}{SW} \label{theorem_saxl_whiston} 
Let $G=\PSL_2(\F_q)$ with $q=p^r$ a prime power and let $A=\set{g_1,\ldots,g_m}$ be a minimal set of generators of $G$. 

If $r=1$ then $|A|\leq 4$. If $r>1$ let $r = p_1^{e_1}\ldots p_n^{e_n}$ be the prime decomposition of $r$ and let 
$$A_i:=A\SM g_i \qtq[and]H_i:=\gen{A_i}.$$  
If $|A|>6$ then up to some reordering of the $g_i$'s and the $p_j$'s one of the following hold.
\doenum{
\item For any $i\geq 3$, $H_i$ is a subfield subgroup and there exists a unique $j$ for which 
$$H_i\leq G_j\cong \PSL_2(p^{r/p_j}).$$
\item For any $i\geq 2$, $H_i$ is a subfield subgroup. For any $1\leq j\leq n$, let $S_j$ be the set of subfield subgroups $H_i$ for which $j$ is minimal subject to $H_i\leq G_j\cong \PSL_2(p^{r/p_j})$. 
Then $|S_1|\leq 2$ and $|S_j|\leq 1$ for any $j\geq 2$.

\item For any $i\geq 1$, $H_i$ is a subfield subgroup. 
For any $1\leq j\leq n$, let $S_j$ be the set of subfield subgroups $H_i$ for which $j$ is minimal subject to $H_i\leq G_j\cong \PSL_2(p^{r/p_j})$. 
Then $|S_1|\leq 3$ and $|S_j|\leq 1$ for any $j\geq 2$.
}
}

As an immediate corollary we get the following claim.
\doCORO{
Let $q$ be a prime power and $G=\PSL_2(\F_q)$ and $A=\set{g_1,\ldots,g_m}$ be a minimal set of generators of $G$.
Let $H_i:=\gen{A\SM \set{g_i}}$. 

If $|A|\geq 7$ then the subgroups $H_i$ which are subfield subgroups 
$$H_i\cong \PSL_2(\F_{q_i})$$ 
satisfy that their underlying fields $\F_{q_i}$ are generating the whole field $\F_q$.
}

\doPROOF{
Let us use the same notations of the previous theorem. 
Let $q=p^r$ and $r = p_1^{e_1}\ldots p_n^{e_n}$ be the prime decomposition of $r$.
By the previous theorem we have three cases to consider. 
In all the cases we get that for any $S_j$ there exist $i=i_j$ and $H_i$ and $r_i$ \ST
$$H_i\cong \PSL_2(p^{r_i})\notin S_j.$$ 
\IOW\ for any $1\le j\le n$, $r_{i_j} \nmid (r/p_j)$. 
Therefore the l.c.m. of these $r_i$'s is 
$$\op{lcm}(r_{i_1},\ldots,r_{i_n})=r$$
so we are done. 
}

Now let us use this corollary to prove the following Theorem.

\doTHM{\l{Tr(A^{[6]})}
Let $\F_q$ be a finite field of order $q$, $G=\SL_2(\F_q)$ and $A$ be a set of generators of $G$. 
Then we have,
\doeqnn{\gen{\Tr(A^{[6]})}=\F_q.}
}

\doPROOF{
By Lemma \r{trace generation} we got that if $\Char(\F)\neq 2$ then 
$$\gen{\Tr(A^{[3]})}=\F,$$ 
so we are only left with the case that $\Char(\F)=2$ and 
$$G=\SL_2(\F_q)=\PSL_2(\F_q)$$ 
with $q=2^r$. 
By taking a subset $A'$ of $A$ if needed, \wLOG\ $A$ is minimal generating set. 
If $|A|\leq 6$ then by Lemma \r{trace generation} we get $\gen{\Tr(A^{[6]})}=\F_q$. 

Now by induction on $r$, and the previous theorem, if $r=1$ then $|A|\leq 4$ and so 
$$\gen{\Tr(A^{[4]})}=\F_q.$$
Otherwise, let $r = p_1^{e_1}\ldots p_n^{e_n}$ be the prime decomposition of $r$. 
Now if $|A|\geq 7$ then by the previous corollary we get proper subfield subgroups 
$$H_i\cong \SL_2(2^{r_i})$$ 
\ST\ the subfields $\F_{2^{r_i}}$ generate $\F_{2^r}$. 
By the induction hypothesis on these $H_i$ which are generated by $A_i=A\SM g_i$, we get 
$$\gen{\Tr(A_i^{[6]})}=\F_{2^{r_i}}.$$
Therefore $\gen{\Tr(A^{[6]})}=\F_q$ as we wanted.
}

\ssec{Avoiding certain traces}

We first start with a useful identity that we will use many times.

\doLEM{\l{lem:trace identity for matrices}
Let $\F$ be a field and $g,h\in \SL_2(\F)$. Then we have,
\doeq{\l{eq:trace identity for matrices}\Tr(g)\Tr(h)=\Tr(gh)+\Tr(gh^{-1}).}
}

\doPROOF{
From the Cayley-Hamilton identity $h^2-\Tr(h)h+I=0$, we get by multiplying by $gh^{-1}$, the matrix identity
$$gh-\Tr(h)g+gh^{-1}=0.$$
Therefore by taking the trace and reordering the identity we are done.
}

\doDEF{\l{def_A_mid_}
Let $G$ be a linear group and let $A\subseteq G(\F)$ and let $X\subseteq \F$. Denote, 
\doalign{
\mbf{A\MID{X}}	&:=	\{g\in A: \Tr(g)\in X\}\\
\mbf{A\NMID{X}}	&:=	\{g\in A: \Tr(g)\notin X\}
}
As usual, when $X=\set{x}$ is singleton we will write just $x$ instead of $X$ and we write $\pm x$ instead of $\set{\pm x}$.
I.e., 
\doalign{
	\mbf{A\MID{x}}			&:=	A\MID{\set{x}}\\
	\mbf{A\MID{\pm x}}	&:=	A\MID{\set{\pm x}}
}
and similarly for $\mbf{A\NMID{x}}$ and $\mbf{A\NMID{\pm x}}$.
}

\doDEF{\l{def_projective_and_fix}
Let $\F$ be a field and let $V(\F)=\F^2\SM \set{\dosmatrix{0\\0}}$. Let 
$$\mBF{\P(\F)} := \ol{V(\F)}= V(\F)/\!\!\sim$$ 
be the \tb{projective line} over $\F$
where for any $u,v\in V(\F)$, 
\doeqnn{\ol{u}=\ol{v} \iff u\sim v \iff \Span(u)=\Span(v).} 
Now let $V=V(\ol{\F})=\ol{\F}^2\SM \set{\dosmatrix{0\\0}}$ and let $G=\SL_2(\ol{\F})$ act on $V$ by left multiplication.
We will be interested in the action of $G$ on $\P(\ol{\F})$ which is induced from the action of $G$ on $V$.
Note that $$g\ol{v}=\ol{v} \iff gv=\lam v\tx{ for some }\lam \in \ol{\F}^{\times}.$$ 
For $g\in G$ denote,
\doeqnn{\mbf{\Fix(g)}:=\set{\ol{v}\in \P(\ol{\F}): g\ol{v}=\ol{v}},} 
the \tb{fix points} of $g$ \WRT\ the action on $\P(\ol{\F})$.
}

The following simple fact is stated also as a definition.
 
\doSFACT{\l{def_unipotent}
Let $G=\SL_2(\F)$. Denote by $\mbf{G_u}$ the non trivial $\pm$ unipotent elements in $G$:

$u\in G_u \iff$ there exist $w\in \SL_2(\F)$ and $a\in \pmO$ and $x\in \F^{\times}$ \ST\ 
$$u^w=a\dobmatrix{1 & x\\0 & 1}=a(I+xE_{12}).$$
If we denote the two columns of $w$ by $w=(w_1,w_2)$ then 
$$\Fix(u)=\set{\ol{w_1}}.$$
We have,
\doalign{
	G_u	&=	G\MID{\pm 2}\SM \pmI\\
			&=	\set{u\in G:|\Fix(u)|=1}.
	}
\IOW\footnote{We write for short $x=\pm y \iff x\in \set{\pm y}$.} 
$$G_u=\set{u\in G: \Tr(u)=\pm 2}\SM \pmI,$$ 
are the elements with exactly one fix point in $\P(\ol{\F})$.
For $A\subseteq G$ denote $$\mbf{A_u}:=A\cap G_u.$$
}

The following simple fact is stated also as a definition.

\doSFACT{\l{def_semi_simple}
Let $G=\SL_2(\F)$. Denote by $\mbf{G_s}$ the semi simple elements in $G$:

$s\in G_s \iff$ there exist $w\in \SL_2(\ol{\F})$ and $y\in \ol{\F} \SM \pmO$ \ST\ 
$$u^w= D_y=\dobmatrix{y & 0\\0 & y^{-1}}.$$
If we denote the two columns of $w$ by $w=(w_1,w_2)$ then 
$$\Fix(s)=\set{\ol{w_1},\ol{w_2}}.$$
We have, 
\doalign{ 
	G_s	&=	G\NMID{\pm 2}\\
			&=\set{s\in G:|\Fix(s)|=2}.
	}
\IOW\footnote{We write for short $x\neq \pm y \iff x\notin \set{\pm y}$.} 
$$G_s=\set{s\in G: \Tr(u)\neq \pm 2},$$ 
are the elements with exactly two fix points in $\P(\ol{\F})$.
For $A\subseteq G$ denote $$\mbf{A_s}:=A\cap G_s.$$
}


\doDEF{
For $A\subseteq G$ we denote for short,
\doalign{
	\opC(A) &=	\opC_G(A)= \set{g\in G: a^g=a \tx{ for any }a\in A}\\
	\opN(A) &=	\opN_G(A)= \set{g\in G: A^g=A}.
	}
}

\doSFACT{\l{sfact:C(g)}
Let $G=\SL_2(\F)$ and let $s\in G_s$ and $u\in G_u$. Then we have, 
\doalign{
	\opC(s) &\subseteq G_s \cup \pmI \\
	\opC(u) &\subseteq G_u \cup \pmI
	}
In fact, 
\doalign{
	\opC(s) &=	\set{s'\in G:\Fix(s')=\Fix(s)}\cup \pmI \\	
	\opC(u)	&= 	\set{u'\in G:\Fix(u')=\Fix(u)}	\cup \pmI.	
	}
}

\doSFACT{\l{sfact:N(C(g))}
Let $G=\SL_2(\F)$ and let $s\in G_s$ and $u\in G_u$. 
Then we have,  
\doalign{	
	\opN(\opC(s))	&= \set{g\in G:g(\Fix(s))=\Fix(s)}\\	
	\opN(\opC(u))	&= \set{b\in G:\Fix(u)\subseteq \Fix(b)}.	
	}
\IOW\ if $\Fix(s)=\set{\ol{w}_1,\ol{w}_2}$ then 
$$g \in \opN(\opC(s)) \iff \tx{ either $g$ fix both $\ol{w}_i$ or $g$ flips between them}.$$
Similarly if $\Fix(u)=\set{\ol{w}_1}$ then 
$$g \in \opN(\opC(u)) \iff g\tx{ fix }\ol{w}_1.$$
}


%

\doDEF{\l{def_FIX}
For a subset $V\subseteq \SL_2(\F)$ denote 
\doeqnn{\Fix(V):=\bigcap_{g\in V}\Fix(g).}
}

The following Lemma is a slight modification of an argument of Helfgott for producing many semi-simple elements (cf. \cite[\S 4.1 Lemma 4.2]{helf}).

\doLEM{[Helfgott]\l{lem:many s.s_}
Let $\F$ be a field, let $G=\SL_2(\F)$ and let $A \subseteq G$ be a finite subset.
Suppose $\gen{A}$ is a non-abelian subgroup\footnote{or we could write for short $[A,A]\neq 1$.}of $G$.  
Then we have,
$$|A^{[3]}\cap G_s|\geq \tfr{4}|A|.$$
}

\doPROOF{
Let $A':=A\SM \pmI$. Then $A'=A_u \cup A_s$ and by the assumption $|A'|\ge 2$.
If $A_u =\emptyset$ then $|A_s|\ge 2$ so $|A_s|\ge \tfr{2}|A|$ so we are done. 
\OW\ let $g\in A_u$ and set $C=\opC_G(g)$ and $B=A \SM C$.
\BTA, $B\neq \emptyset$.

If $h\in B_u$ then for some $x,y\in\F^\times$ and $w\in \SL_2(\F)$ and $a,b\in \set{\pm 1}$ 
we have	
\doalign{
	g^w	&=	\dobmatrix{a&x\\0&a} = aI+x\opE_{12}\\
	h^w	&=	\dobmatrix{b&0\\y&b} = bI+y\opE_{21}
	} 
and $(h^{-1})^w=bI-y\opE_{21}$.
Therefore  			
$$\Tr(gh^{\pm 1})=2ab\pm xy.$$
Now if $\Char(\F)=2$ then both $gh, gh^{-1}$ are semi simple elements
and if $\Char(\F) \neq 2$ then at least one of $gh, gh^{-1}$ is semi simple. 

\TF\ we get that for any $h\in B$, either $h\in G_s$ or $gh\in G_s$ or $gh^{-1}\in G_s$. 
\TF\ $A^{[2]}$ contains at least $\fr{2}|B|$ semi-simple elements so
\doeq{\l{bound_many s.s_1}|A^{[2]}\cap G_s|\geq \fr{2}|B|=\fr{2}(|A|-|\opC_A(g)|).}

\OTOH, if $h\in A \SM C$ then $h\ddot \opC_A(g)\subseteq A^{[2]}\SM C$. 
Set $B'= A^{[2]}\SM C$ and so $|B'|\geq |C_A(g)|$.
\TF\ by applying the previous argument  \re{bound_many s.s_1} with $\mbf{B}=B'$ we get that,
\doeq{\l{bound_many s.s_2}|A^{[3]}\cap G_s|\geq \fr{2}|B'|\ge \fr{2}|\opC_A(g)|.}

Putting together  \re{bound_many s.s_1} and  \re{bound_many s.s_2} we get,
\eqsplit{bound_many s.s_3}{ 
	|A^{[3]}\cap G_s|	&\ge \tfr{2}\MAX{|A|-|C_A(g)|,|C_A(g)|}\\
										&\ge \tfr{4}|A|.	\QED
	}
}

The following Lemma is a slight variant of Lemma \r{lem:many s.s_}.

\doLEM{\l{lemma_many non zero trace_}
Let $\F$ be a finite field. 
Let $G=\SL_2(\F)$ and let $A \subseteq G$ and suppose $\gen{A}=G$. 	
Then we have,
\doeqnn{|A^{[3]} \NMID{0}|\geq \fr{4}|A|.}
}


\doPROOF{
If $\Char(\F)=2$ then $G_s=G\NMID{0}$ \sWADB Lemma \r{lem:many s.s_}. 
\OW\ $\Char(\F)\neq 2$ and therefore $G\MID{0}\subseteq G_s$. 
If $0\notin \Tr(A^{[3]})$ \tWAD. 
\OW\ fix $g\in A^{[3]}|_0$ and let $\om \in \ol{\F}$ with $\om^2=-1$. 
\TF\footnote{We denote $Spec_{\ol{\F}}(g)$ to emphasize that we take all the eigen values in $\ol{\F}$.}\ $\Lam(g)=Spec_{\ol{\F}}(g)=\set{\pm \om}$.

Note that 
\eqsplit{eq:nil trace}{
	\Tr(g)=0	&\iff g^2=-I\\ 
						&\iff g^{-1}=-g.
	}
Denote $C=\opC_G(g)$ and $N=\opN_G(C)$. 
\BTA\ and by fact \r{sfact:N(C(g))} $$A \nsubseteq N.$$
Set $B=A\SM N\neq \emptyset$ and let $h\in B$. 
If $\Tr(h)=0$ then 
\eqsplit{bound_many non zero trace_1}{
	\Tr(gh)=0 &\iff ghgh=-I\\
						&\iff g g^h=I\\ 
						&\iff g^h=g^{-1}.
	}
\TF\ $\Tr(gh)=0 \Rightarrow h\in N$ \RL!\ (since we we took $h\notin N$).
Therefore we got that either $\Tr(h)\neq 0$ or $\Tr(gh)\neq 0$. 
So
\doeq{\l{bound_many non zero trace_2} |A^{[2]}\NMID{0}|\geq \fr{2}|B|=\fr{2}(|A|-|A\cap N|).}
\OTOH\ if $h\in A \SM N$ then $h(A \cap N)\subseteq A^{[2]}\SM N$ therefore, 
$$|A^{[2]}\SM N|\geq |A\cap N|.$$
Therefore by applying the previous argument  \re{bound_many non zero trace_2} with $\mbf{B}=B'=A^{[2]}\SM N$ we get that 
\doeq{\l{bound_many non zero trace_3} |A^{[3]}\NMID{0}|\geq \fr{2}|B'|\ge \fr{2}|A\cap N|.}
Combining  \re{bound_many non zero trace_2} and  \re{bound_many non zero trace_3} we get
\eqsplit{bound_many non zero trace_4}{
	|A^{[3]}\NMID{0}|	&\geq \fr{4}|A|.	\QED
	}
}

\doLEM{\l{lemma_many traces outside a field_}
Let $\F$ be a finite field and let $G=\SL_2(\F)$. 
Suppose $A\subseteq G$ with $\gen{A}=G$ and let $\E<\F$ be a proper subfield.  
Then we have,
  \doeqnn{|A\NMID{\E}|>0 \RightarrowLB |A^{[4]}\NMID{\E}| \geq \fr{12}|A|.}
}

\doPROOF{
Denote $B=A^{[3]}$. 
If $|B\NMID{\E}|\geq \tfr{12}|A|$ then we are done so assume 
$$|B\NMID{\E}|<\fr{12}|A|.$$
From Lemma \r{lemma_many non zero trace_} we get that 
\doeqnn{|B\NMID{0}|\ge \fr{4}|A|.}
Therefore 
\doeqnn{|B\MID{\E^{\times}}|>(\fr{4}-\fr{12})|A|=\fr{6}|A|.}
From Lemma \r{lem:trace identity for matrices} if $g \in G\NMID{\E}$ and $h\in G\MID{\E^{\times}}$ then, 
$$\tq[either]\Tr(gh^{-1})\notin \E \qtq[or] \Tr(gh)\notin \E.$$
By the assumption there is $g\in A\NMID{\E}$ therefore we get
$B':=gB\subseteq A^{[4]}$ and so 
\doalign{
	|A^{[4]}\NMID{\E}|	&\geB |B'\NMID{\E}|\\
											&\geBY{eq:trace identity for matrices} \fr{2}|B\MID{\E^{\times}}|\\
											&\gnB\fr{12}|A|.	\QED
	}
}


Therefore we get immediately the following result.
\doCORO{\l{corr_escaping all subfields_}
Let $\F$ is a finite field and let $G=\SL_2(\F)$. 
Let $A\subseteq G$ and suppose $\gen{A}=G$ and $\gen{\Tr(A)}=\F$. 
Then for \ti{any} proper subfield $\E<\F$ we have, 
\doeqnn{|A^{[4]}\NMID{\E}|\geq \fr{12}|A|.}
}

\doCORO{\l{corr_escaping from all subfields_}
Let $\F$ is a finite field and let $G=\SL_2(\F)$.
Let $A\subseteq G$ and suppose $\gen{A}=G$.
Then for any proper subfield $\E<\F$ we have, 
\doeqnn{|A^{[9]}\NMID{\E}|\geq \fr{12}|A|.}
}
\doPROOF{
By Lemma \r{Tr(A^{[6]})}, $\gen{\Tr(A^{[6]})}=\F$ therefore 
\doeqnn{|A^{[6]}\NMID{\E}|>0.}
Now as in the proof of Lemma \r{lemma_many traces outside a field_} we get that either $|A^{[3]}\NMID{\E}|\geq \fr{12}|A|$ (and then we are done) or
\doeqnn{|A^{[3]}\MID{\E^{\times}}|>\fr{6}|A|.}
Therefore if take $b\in A^{[6]}\NMID{\E}$ and $B':=A^{[3]}\MID{\E^{\times}}$ and $B'':=bB' \subseteq A^{[9]}$ 
then we get
\doeqnn{|A^{[9]}\NMID{\E}|\ge |B''\NMID{\E}| \geBY{eq:trace identity for matrices} \fr{2}|B'|>\fr{12}|A|. \QED}
}

%

\sec{Growth properties of $\SL_2(\F_q)$}\l{sec:Growth properties of SL_2}		

\ssec{Some useful Growth properties}

\doDEF{
Let $G$ be a group and $g,h\in G$. Define the conjugacy class equivalence by 
$$g\sim h \iff g^G=h^G.$$ 
I.e., $g\sim h \iff g^x=x^{-1}gx=h$ for some $x\in G$.
Given a subset $A\subseteq G$ denote $$\widetilde{A}=A\over\sim.$$
By abuse of notation we will view $\widetilde{A}\subseteq A$ as a set of representatives so: 
$$\forall a\in A,\ \exists !b\in \widetilde{A} \tx{ \ST\ } a\sim b.$$
}

The following useful Lemma connects growth and commutativity. 

\doLEM{\CITE{\S 4.1 Proposition 4.1}{helf}\l{lemma_group_elements_}	
Let $G$ a finite group and let $\emptyset\neq A\subseteq G$.
Then there exists $a\in A$ such that,
\doeq{\l{lemma_group_elements_bound_1} |\opC_{A^{-1}A}(a)|\ge \frac{|\widetilde{A}||A|}{|A^{-1}AA|}.}

If $\gen{A}=G$ then for any proper subgroups $H,K<G$ we have, 
\doeq{\l{lemma_group_elements_bound_2} |A^{[4]}\SM(H\cup K)|>\fr{4}|A|.}
}

\doPROOF{
Let $a,b\in A$ and $g\in G$ and suppose $g^a=g^b$. Then we have,
\doeqnn{\l{bound_group_elements_1} ba^{-1}\in \opC_{Aa^{-1}}(g)\subseteq \opC_{AA^{-1}}(g).}
Therefore we get,
\doeqnn{b\in \opC_{Aa^{-1}}(g)a \subseteq \opC_{AA^{-1}}(g)a.} 
Therefore for any $g\in A$ we get that $g^A\subseteq A^{-1}AA$ and 
\doeq{\l{bound_group_elements_2} 
	\frac{|A|}{|\opC_{AA^{-1}}(g)|}\leq |g^A|.
	} 
\OTOH\ if we denote $\Lam=\widetilde{A}$ then
\doeq{\l{bound_group_elements_3}
	\fr{|\Lam|}\sum_{g\in \Lam} |g^A| = \fr[|\Lam^A|]{|\Lam|}\le \frac{|A^{-1}AA|}{|\Lam|}.
	} 
Therefore there exists $g\in \widetilde{A}\subseteq A$ s.t.
\doeqnn{
	\frac{|A|}{|\opC_{AA^{-1}}(g)|} \leBY{bound_group_elements_2} |g^A| \leBY{bound_group_elements_3} \frac{|A^{-1}AA|}{|\widetilde{A}|}
}
so by arranging the inequality we are done with  \re{lemma_group_elements_bound_1}.\good

Now suppose $\gen{A}=G$.
Since $A \SM H \neq \emptyset$ we get that for $a\in A \SM H$, $a (A\cap H) \subseteq A^{[2]}\SM H$ therefore 
\doeqnn{|A^{[2]}\SM H | \ge \MAX{|A\SM H|,|A\cap H|}\ge \fr{2}|A|.}

If $H=K$ then we are done. If $A \subseteq H\cup K$ then there exists $a,a'\in A$ \ST\ $a\in H\SM K$ and $a'\in K\SM H$ therefore 
$a a'\in A^{[2]}\SM (H\cup K)$. 
\IAC\ there exists $b\in A^{[2]}\SM (H\cup K)$.
Denote $B=A^{[2]}\SM H$ so $b \in B\SM K$ therefore $b(B\cap K)\subseteq A^{[4]}\SM K$ therefore 
\doeqnn{|A^{[4]}\SM(H\cup K)|\ge \MAX{|B\SM K|,|B\cap K|} \ge \fr{2}|B|\ge \fr{4}|A|}
\sWAD\ with  \re{lemma_group_elements_bound_2}. \good
}


\doCORO{\l{corr_large centrelizer_}
Let $\F$ be a field. Let $G$ be a subgroup of $\GL_n(\F)$ and let $A\subseteq G$ be a finite subset. 
Let $B\subseteq A$ with $|B|\ge c|A|$ for some $c\in \R_+$.
Then there exists $b\in B$ such that,
\doeq{\l{remark large centrelizer}	|\opC_{AA^{-1}}(b)| \ge c\frac{|\Tr(B)||A|}{|A^{-1}AA|}.}
}

\doPROOF{
Since conjugate elements have the same trace we get, 
\doeqnn{|\widetilde{A}|\ge |\Tr(A)|.}
Therefore by Lemma \r{lemma_group_elements_} there exists $a\in A$ such that,
\doeqnn{|\opC_{AA^{-1}}(a)|\geBY{lemma_group_elements_bound_1} \frac{|\Tr(A)||A|}{|A^{-1}AA|}.}
\TF\ if $B\subseteq A$ and $|B|\ge c|A|$ then there exists $b\in B$ such that\footnote{I want to thank H.Helfgott for helpful discussion concerning this variant.},
\doalign{
	|\opC_{AA^{-1}}(b)| 	&\ge	|\opC_{BB^{-1}}(b)|\\
										&\ge	\frac{|\Tr(B)||B|}{|B^{-1}BB|}\\
										&\ge	c\frac{|\Tr(B)||A|}{|A^{-1}AA|}. \QED	
	}
}

A variant of the following Lemma was proved in \cite[Proposition 4.10]{helf}. 
Here, we will show another way of proving it.

\doLEM{\label{large trace set}
Let $\F$ be a field and let $G=\SL_2(\F)$.
Let $g\in G_s$ be a semi simple element.
Let $h\in G$ and suppose $Fix(h)\SM Fix(g)\neq \emptyset$.
Define the function
$F:\SL_2(\F)\rightarrow \F^3$ by 
\doeqnn{F(b)=(\Tr(b),\Tr(gb),\Tr(hb)).} 
Then $\mult(F)\leq 2$. 
\IP, for any subset $B\subseteq G$,
\doeq{\l{remark large trace}\fr{2}|B| \le |F(B)| \leq |\Tr(B)||\Tr(gB)||\Tr(hB)|.}
}


\doPROOF{
There exists $w\in \SL_2(\ol{\F})$ such that 
\doalign{
	g	&=	\dobmatrix{\alpha & a\\0 &\alpha^{-1}}^w\\
	h	&=	\dobmatrix{\beta & 0 \\ b & \beta^{-1}}^w
	} 
with $b\in \ol{\F}^{\times}$ and $\alpha\notin \pmO$. 
Let $g'=\dobmatrix{x & y \\z & w}^w\in \SL_2(\F)$. 

We need to show that for any $c_1,c_2,c_3$ there are at most two $g'$ with
\doeqnn{\docases{\det(g')	&=1  \\ F(g')	&=(\Tr(g'),\Tr(gg'),\Tr(hg'))=(c_1,c_2,c_3) }}
By opening trace equalities we get the linear system 
\doeqnn{
	\dopmatrix{1 & 1 & 0 & 0\\ \alpha & \alpha^{-1} & 0 & a\\ \beta & \beta^{-1} & b & 0}
	\dopmatrix{x \\ w \\ y \\ z}=
	\dopmatrix{c_1 \\ c_2 \\ c_3}.
}
Denote $A=\dopmatrix{1 & 1 & 0 & 0\\ \alpha & \alpha^{-1} & 0 & a\\ \beta & \beta^{-1} & b & 0}$
and $\ol{x}=\dopmatrix{x \\ w \\ y \\ z}$ and $\ol{c}=\dopmatrix{c_1 \\ c_2 \\ c_3}$.
Therefore, from our assumption on $b$ and $\alpha$, 
$$\rank(A)=3$$ so the set of solutions $A^{-1}(\ol{c})$ is either empty, or a one dimensional affine linear subspace\footnote{i.e., $A^{-1}(\ol{c})$ is a dilation of a one dimensional linear subspace of $\ol{\F}^4$.} of $\ol{\F}^4$.
Note that for any $z$ there is exactly one triple $(x,w,y)$ \ST\ $g'$ is a solution.
\OTOH, $g'\in \SL_2(\F)$ so $xw-yz=1$ and therefore there at most two solutions $g'$ on the affine line 
$A^{-1}(\ol{c})$ with $\det(g')=1$. 
\IOW
\doeqnn{|A^{-1}(\ol{c})\cap \SL_2(\ol{\F})|\leq 2.	\QED}
}



%

\ssec{Avoiding subvarieties}


\doDEF{
Let $\F$ be a field. Let $G$ be a group and let $(V,\rho)$ be a finite dimensional representation of $G$ over $\F$. 
When the action will be clear from the context we will write the \tb{linear action} on $V$ simply by $gv$ instead of $\rho(g)v$.
Let $W_1,\ldots W_m <V$ be proper subspaces of $V$ and let \doeqnn{W=\bigcup_{i=1}^m W_i.}
We will assume that the above union is \tb{non trivial} in the sense that 
\doeqnn{W_i \leq W_j \Rightarrow i=j.}
We will call $W$ a \tb{linear variety} with \tb{decomposition}\footnote{if the union is non trivial then the decomposition is unique.} $W=\bigcup_{i=1}^m W_i$.
Denote \doeqnn{\mbf{\Stab_G(W)}=\set{g\in G: gW=W}.}
We will sometimes abbreviate and write 
$$\mbf{G_W}=\mbf{\Stab(W)}=\Stab_G(W)$$ 
when the group $G$ is clear from the context.
Denote,
\doalign{
	\mbf{\dim(W)}		&:=	\max_{i}\set{\dim(W_i)}\\
	\mbf{\deg_d(W)}	&:=	|\set{i: \dim(W_i)= d}|\\
	\mbf{\deg(W)}		&:=	\deg_{\dim(W)}(W).
}
}

The following ``escaping Lemma'' will be useful. 
The following proof is a slight modification of \cite[\S 4.2 Lemma 4.4]{helf}.

\doLEM{[Helfgott]\label{lemma_escape from orbits_}
For any $n,m\in \N_+$ there exists $k\in \N_+$ such that the following holds. 
Let $G$ be a group and let $(V,\rho)$ be a finite dimensional representation of $G$ over a field $\F$. 
Let $W_1,\ldots W_m\leq V$ be subspaces of $V$ and suppose $W=\CUP[i] W_i$ is a linear variety with $\dim(W)\leq n$.
Let $A$ be a subset of generators of $G$. 
Let $0\neq w\in V$ and denote the orbit of $w$ by $O:=Gw$ and 
$$V_w:=\F[G]w=\Span(O).$$
Suppose $O\nsubseteq W$. 	

Then for any $0\neq w'\in V_w$ there exists $g\in A^{[k]}$ \ST\ $gw'\notin W$.
\IP\ for any $w'\in O$ there exists $g\in A^{[k]}$ \ST\ $gw'\notin W$.
}

\doPROOF{
Note that the claim is trivially true for $w'\in V_w \SM W$ so we need 
to prove it for $0\neq w'\in V_w\cap W$.
\IP, if $V_w\cap W=0$ we are done.

\WLOG\ $W=\CUP[i] W_i$ is the decomposition of $W$ as a union of spaces.
Set for $1\leq i\leq m$, $O_i:=O\cap W_i$ and $V_i:=V_w \cap W_i$ and 
\doeqnn{W_{(0)}:=V_w \cap W=\bigcup_{i=1}^m V_i.}
\BTA\ for any $i\leq m$, $O \nsubseteq W_i$. 
\TF\ $$V_w =\Span(O)\nsubseteq W_i$$ so $V_i< V_w$ and $V_i\leq W_i$.
Now for any $g\in G$, $$gO_i=g(O\cap W_i)=O\cap gW_i$$
and $gV_i=g(V_w\cap W_i)=V_w\cap gW_i$.
Note that $O_i=\emptyset \iff V_i=0$.

If $V_i=0$ then for any $g\in G$, 
$$0=V_i\cap gV_i=V_w\cap W_i \cap gW_i<W_i.$$
Now suppose $O_i=O\cap W_i\neq \emptyset$ for some $i\leq m$ and let $$x_i \in O_i\subseteq V_i.$$
Since $Gx_i=O \nsubseteq W_i$ we get that
there exists $g_i\in G$ \ST\ $g_ix_i \notin W_i$ so $$g_iV_i\nsubseteq V_i.$$
\IOW\ $\Stab(V_i)\neq G$. Therefore $V_i \cap g_iV_i < V_i$ so
\doeqnn{V_w \cap W_i \cap g_iW_i  < V_w \cap W_i\leq W_i.}
Since $\gen{A}=G$ we can choose $g_i$ to be $g_i\in A$.

Therefore if $\dim(W)>0$ and $\dim(W_i)=\dim(W)$
then there exists $a_1\in A$ \ST\ 
$$V_i\cap a_1V_i = V_w \cap W_i \cap a_1W_i < W_i$$
and for all other $j\leq m$, $V_j\cap a_1V_j\leq W_j$.
Set for any $1\leq j\leq m$, $W_{1j}:=V_j\cap a_1V_j=V_w \cap W_j \cap a_1W_j$ 
and \doeqnn{W_{(1)}:=W_{(0)}\cap a_1W_{(0)}=\bigcup_{j=1}^m W_{1j}.}
Therefore 
\doeqnn{W_{(1)}\subsetneq W_{(0)} \subseteq  W}
so either 
\doeqnn{\dim(W_{(1)})<\dim(W_{(0)})\leq \dim(W)} 
or 
\doeqnn{\deg(W_{(1)})<\deg(W_{(0)})\leq \deg(W).}
Therefore by iterating the previous step either $W_{(1)}=0$ or we can find $a_2\in A$ 
\ST\ for $W_{(2)}:=V_w \cap W_{(1)}\cap a_2W_{(1)}$ we get 
$$\tq[either]\dim(W_{(2)})<\dim(W_{(1)})\qtq[or]\deg(W_{(2)})<\deg(W_{(1)}).$$
Therefore for some $k\leq mn$ we get that $W_{(k)}=0$ therefore 
\doeqnn{\bigcap_{g\in A^{[k]}}g(V_w\cap W)=0.}
Therefore for any $0\neq w'\in V_w\cap W$ there exists $g\in A^{[k]}$ \ST\ $w'\notin W$
\sWAD.
}


Now we will prove the following result. 

\doCORO{\l{lemma_escaping from a union of subspaces_}
There exists $k\in \N_+$ such that the following holds
for any finite field $\F$ of size $|\F|>3$, and for any subset of generators $A$ of $\SL_2(\F)$.
For any $u\in GL_2(\ol{\F})$, there exists $a\in A^{[k]}$, \ST\ $a^u$ has no zero entries.
}


\doPROOF{
Denote $G:=\SL_2(\F)$ and $V:=\M_2(\ol{\F})$ and for $1\le i,j \le 2$
\doeqnn{W_{ij}	:= \set{\dopmatrix{a_{11}& a_{12}\\a_{21}&a_{22}}\in V:a_{ij}=0}}
and $W=\CUP[i,j]W_{ij}$.
Equivalently, if $g=\dobmatrix{a_{21}& a_{22}\\a_{11}& a_{12}}\in V$ then 
\doeqnn{a_{ij}=0 \iff	g e_j	= \lam e_i \tx{ for some }\lam\in \ol{\F}} 	


Now we are going to use Lemma \r{lemma_escape from orbits_} with the group $G^u$ and the orbit $O=G^u$ of $w'=I$ and the linear variety $W$.
We can use Lemma \r{lemma_escape from orbits_} if we show that $G^u \nsubseteq W$. 
We will show that $|G^u \cap W|<|G|$ so $G^u \nsubseteq W$.

Let $u=(u_1,u_2)$ where $u_i$ are the columns of $u$.
Therefore for any $g\in G^u\cap W$ there exist $1\le i,j\leq 2$ \ST\ $g\ol{u_i}=\ol{u_j}$. 
I.e., $g u_i = \lam u_j$ for some $\lam \in \ol{\F}^{\times}$.
Denote $$G_{ij}:=\set{g\in G:g\ol{u_i}=\ol{u_j}}.$$
So $G^u\cap W= \CUP[i,j] G_{ij}$.
In order to prove $|G^u \cap W|<|G|$ we will bound $|\CUP[i,j] G_{ij}|$ from above.

Let us choose for any $i\in\set{1,2}$ some $u_i'\in \F^2\SM \set{0}$ \ST\ $u_i,u'_i$ are linear independent.  
Now if $g,g'\in G_{ij}$ then $g u_i=\lam u_j$ and $g' u_i=\lam' u_j$ for some $\lam,\lam' \in \ol{\F}$. 
Note that knowing $gu'_i$ and $gu_i$ determine $g$ therefore if 
$g,g'\in G_{ij}$ and $gu'_i=g'u'_i\in \F^2\SM \set{0}$ then we must have $\lam=\lam'$
since $\det(g)=\det(g')=1$. 
Therefore we conclude that for any $i,j$ we have $|G_{ij}|\leq |\F|^2-1$.
Therefore $|G^u\cap W|=|\bigcup G_{ij}|\leq 4(|\F|^2-1)-1$ since $I\in G_{11}\cap G_{22}$.
So if $|\F|=q\geq 4$ then 
\doeqnn{|G^u|=|SL_2(\F)|=q(q^2-1)>4(q^2-1)-1\geq |\bigcup G_{ij}|}
so in particular $G^u \nsubseteq W$. 

Therefore we can apply Lemma \r{lemma_escape from orbits_} to get the following.
For any $u\in GL_2(\ol{\F})$ there exist $a\in A^{[k]}$ \ST
\doeqnn{a^u\tx{ has no zero entries.}	\QED}
}

\ssec{\l{prf helfgott reduction} Reduction from matrices to traces}
%
%
%



\doDEF{
Let $\F$ be a field and let $g,h\in \SL_2(\F)$. We will say the $g$ and $h$ are \tb{entangled} (or \tb{simultaneously triangular})
if $$\tq[either]\Fix(h)\subseteq \Fix(g) \qtq[or] \Fix(g)\subseteq \Fix(h).$$
}

The following Lemma will be useful later \CITE[cf. ]{Lemmas 4.7, 4.9}{helf}.

\doLEM{[Helfgott]\l{lemma_sl2 elements} 
There exists $C>0$ \ST\ the following properties hold for any field $\F$. 
Let $g,h\in \SL_2(\F)$ and suppose they are not entangled. 
Then there exists $w\in \SL_2(\ol{\F})$ \ST
\doeq{g^w=\dobmatrix{a & x\\0 & a^{-1}}\qtq[and]h^w=\dobmatrix{b & 0\\y & b^{-1}}.\l{lemma_sl2 elements_bound_a}}
Moreover if $g\in G_u$ then $a=\pm 1$ and $x\neq 0$ (and similarly for $h\in G_u$).

Let $V\subseteq \SL_2(\F)$ be a finite subset of diagonal matrices and suppose $V\nsubseteq \pmI$.
Let $g\in \SL_2(\F)$.
If $g$ has no zero entries\footnote{i.e., $abcd \neq 0$ where $g=\dosmatrix{a & b\\c & d}$.}
then we have,
\doeq{|VgVg^{-1}V|\geq \fr{C}|V|^3. \l{lemma_sl2 elements_bound_b}}
If $U\subseteq \SL_2(\F)$ is a finite non empty subset which has no triangular matrices%
\footnote{i.e., $bc\neq 0$ where $u=\dosmatrix{a & b\\c & d}$.} 
then we have,
\doeq{|\Tr(UU^{-1})|\ge \fr{C}\frac{|U|}{|\Diag(U)|}. \l{lemma_sl2 elements_bound_c}}
}

%

\doPROOF{

By taking the two eigen vectors $w_1,w_2\in \ol{\F}^2$ of $g$ and $h$ \RESP\ \ST\ 
$\ol{w_1}\in \Fix(g)\SM \Fix(h)$ and $\ol{w_2}\in \Fix(h)\SM \Fix(g)$,
and normalize them if needed, we get  \re{lemma_sl2 elements_bound_a}. \good

Suppose $V=D_S$ i.e., $S:=\set{s\in \F:\dobmatrix{s & 0\\0 & s^{-1}}\in V}$.
For any $g'=\dobmatrix{x' & y'\\z' & w'}$ we get 
$V g' V= \set{\dobmatrix{stx' & st^{-1}y'\\s^{-1}tz' & s^{-1}t^{-1}w'}:s,t\in S}.$ 
\TF\ $\Prod(Vg'V)=\Prod(g')$.
Moreover, we see that unless $g'$ is diagonal or anti-diagonal\footnote{i.e., has the form 
$\dobmatrix{\neq 0 & 0\\0 & \neq 0}$ or 
$\dobmatrix{0 & \neq 0\\\neq 0 & 0}$.} 
we can recover from any element of
$Vg'V$ the values $s^2,t^2$ so $|Vg'V|\ge \fr{4}|V|^2$. 
Now let $g'\in V^g$ so
\doalign{
	g'	&=	\dobmatrix{x' & y'\\z' & w'}
			= \dobmatrix{s & 0\\0 & s^{-1}}^{\dosmatrix{a & b\\c & d}}\\
			&= \dobmatrix{ads-s^{-1}bc & (s-s^{-1})db\\(s^{-1}-s)ac & ads^{-1}-sbc}.
	} 
\TF\ if $s\neq \pm 1$ then $x'y'z'w'\neq 0$ so in particular $g'$ is neither diagonal nor anti diagonal.
%
%
Altogether we get that 
\doeqnn{|VV^gV|\ge \fr{4}|V|^2|V\SM \pmI|\ge \fr{12}|V|^3}
\sWAD\ with  \re{lemma_sl2 elements_bound_b} .\good 

For any $g\in U$ denote by $U_g$ the subset of all $g'\in U$ with the same diagonal as $g$. 
Consider the trace map $\Tr:g\lr{U_g}^{-1} \to \Tr(UU^{-1})$.
By calculating the trace $\Tr(gg'^{-1})$ one see that each fiber is of size at most 2. 
\TF\ for any $g\in U$ we have $|\Tr(UU^{-1})|\ge \fr{2}|U_g|$. 
Since there exists $g$ with 
\doeqnn{|U_g|\ge \frac{|U|}{|\Diag(U)|}}
we get,
\doeqnn{|\Tr(UU^{-1})|\ge \fr{2}|U_g|\geq \fr{2}\frac{|U|}{|\Diag(U)|}}
\sWADW  \re{lemma_sl2 elements_bound_c}.
}

The following Lemma is the corner stone which connects the Growth of matrices and the Growth of traces \CITE[cf. ]{Propositions 4.8, 4.10}{helf}.
\doLEM{[Helfgott]\label{sl2 generators}
There exist $k\in \N_+$ and $C\in \R_+$ such that the following holds
for any finite field $\F$.
Let $G=\SL_2(\F)$ and let $A\subseteq G$ be a subset of generators of $G$.
Then we have,
\doeq{\l{sl2 generators_traces_geq} |\Tr(A^{[k]})|> \fr{C}|A|^{1/3}}
There exist $V\subseteq A^{[k]}$ and $w\in \SL_2(\ol{\F})$ \ST\ $V^w$ are diagonal and 
\doeq{\l{sl2 generators_V_geq}|V|\geq \fr{C}\frac{|\Tr(A)||A|}{|A^{[k]}|.}}
We also have,
\doeq{\l{sl2 generators_traces_leq} |\Tr(A)|\leq C\frac{|A^{[k]}|^{4/3}}{|A|}.}
}

\doPROOF{
By Lemma \r{lem:many s.s_} there exists $k_0\in \N_+$ \ST\ for $A_0:=A^{[k_0]}$ we have
$$|A_0\cap G_s|\gg |A|.$$
Let $h\in A_0\cap G_s$ be a semi simple element in $A_0$ and let $\set{\ol{v},\ol{u}}=\Fix(h)$ be 
its two fix points in $\P(\ol{\F})$.
\WLOG\ 
\doeqm{(v,u)\in \SL_2(\ol{\F})}	
and let us write from now the $\SL_2(\F)$ elements \WRT\ the basis\footnote{We denote a basis of a space as a tuple of vectors and not as a set of vectors. \TF\ the notation $(v,u)$ has a double meaning either as matrix (a tuple of columns) or as tuple of vectors.} $(v,u)$ of $\ol{\F}^2$.

Denote by $H$ and $K$ the stabilizers of these points
\doeqnn{H:=\set{g\in G: g\ol{v}=\ol{v}}\qtq[and]K:=\set{g\in G: g\ol{u}=\ol{u}}.}
By Lemma \r{lemma_group_elements_} there exists $k_1\in \N_+$ \ST\ 
for $A_1:=A_0^{[k_1]}$ and $U:=A_1\SM (H\cup K)$ we have, 
\doeqnn{|U|\byeq{lemma_group_elements_bound_2}{\gg} |A|.}	
Since $U$ has no triangular matrices we get by Lemma \r{lemma_sl2 elements} 
some $k_2\in \N_+$ \ST\ for $A_2:=A_1^{[k_2]}$ and $D:=\Diag(U)$ we have,
\doeqnn{|\Tr(A_2)|\geq |\Tr(UU^{-1})|\byeq{lemma_sl2 elements_bound_c}{\gg} \frac{|U|}{|D|} \gg \frac{|A|}{|D|}.}

In order to complete the proof of  \re{sl2 generators_traces_geq} we will show that $|D|\leq |\Tr(A_2)|^2$.
Now for any $t\in \Tr(U)$ denote by $S_t$ the elements in $D$ with this sum
and by $U_t$ the elements in $U$ with this trace. 
Therefore for some $t\in \Tr(U)$ we have,
\doeqnn{|U_t|\ge |S_t|\ge \frac{|D|}{|\Tr(U)|}\geq \frac{|D|}{|\Tr(A_1)|}.}
Therefore in order to complete the proof of  \re{sl2 generators_traces_geq} we will show that for any $t\in \Tr(U)$,
\doeqnn{|U_t|\leq |\Tr(A_2)|}
Indeed since $h=\dobmatrix{r & 0\\0 & r^{-1}}\in A_0$ then for any 
$g=\dobmatrix{a & b\\c & d}\in U_t$ we have
\doeqnn{\Tr(hg)=ra+r^{-1}d=(r-r^{-1})a+r^{-1}t}
therefore the trace map $\Tr:hU_t \to \F$ is injective so 
\doeqnn{|U_t|=|\Tr(hU_t)| \le |\Tr(A_2)|}
\sWADW  \re{sl2 generators_traces_geq}.\good

Denote by $B:=A_0\cap G_s$ the semi simple elements in $A_0$.
As we seen before previously $|B| \gg |A|$. 
By corollary \r{corr_large centrelizer_} there exists $b\in B$ \ST\ for 
$A_3:=A_0^{[3]}$ and $V:=\opC_{A_3}(b)$ we get,	
\eqSplit{sl2 generators_C_geq}{
	|V|	&\eqB |\opC_{A_3}(b)|\\
			&\geB	|\opC_{BB^{-1}}(b)|\\
			&\geBY{remark large centrelizer}	\frac{|\Tr(B)||A_0|}{|A_0^{-1}A_0 A_0|}\\
			&\ggB \frac{(|\Tr(A_0)|-2)|A|}{|A_3|}\\
			&\ggB \frac{|\Tr(A)||A|}{|A_3|}.
	} 
Since $b$ is semi-simple there exists $w'\in \SL_2(\ol{\F})$ \ST
\doEQ{eq:V,w,A_3}{V^{w'}\tx{ are diagonal and }V\subseteq A_3}
\sWADW  \re{sl2 generators_V_geq}.\good


By  \re{sl2 generators_C_geq} and  \re{eq:V,w,A_3} there exists a basis $w'\in \SL_2(\ol{\F})$ and $V\subseteq A_3$ \ST\ $V^{w'}$ are all diagonal and 
\doeqnn{|\Tr(A)|\llBY{sl2 generators_C_geq}	|V|\frac{|A_3|}{|A|}}
where $A_3=A^{[k_3]}$ and $k_3:=3k_0$.

By corollary \r{lemma_escaping from a union of subspaces_} there exists $k_4\in \N_+$ 
and $g\in A_4:=A^{[k_4]}$ \ST\ $g^{w'}$ has no zero entries.
Now we set \doeqm{k_5:=5 \MAX{k_4,k_3}.}
Therefore $VV^gV \subseteq A_5:=A^{[k_5]}$. 
\TF\ by Lemma \r{lemma_sl2 elements} we get, 
\doEQsplit{sl2 generators_traces_leq_b}{
		|A_5| &\geB|VV^gV|\\
					&\ggBY{lemma_sl2 elements_bound_b} |V|^3.
}					
Therefore we conclude,
\doalign{
	|\Tr(A)|	&\llB	|V|\frac{|A_3|}{|A|}\\
						&\llBY{sl2 generators_traces_leq_b}	\frac{|A_5|^{4/3}}{|A|}
	}					
\sWADW  \re{sl2 generators_traces_leq}.\good
}



\ssec{Corollaries}


Let us collect the properties that we will exploit soon.
\doTHM{[Helfgott]\l{thm_helfgott reduction to traces_}	
There exist $k\in \N_+$ and $C\in\R_+$ such that the following holds for any finite field $\F$.
Let $A$ be a subset of generators of $\SL_2(\F)$.
Then we have,
\doAlign{
	\l{large traces} 			|\Tr(A^{[k]})|		&> \fr{C}|A|^{1/3}\\
	\l{back to growth} 		|\Tr(A)|					&<C\frac{|A^{[k]}|^{4/3}}{|A|}\\
	\l{many semi simples} |A^{[k]}\cap G_s|	&> \fr{C}|A|
	}
}

\doPROOF{
Parts  \re{large traces} and  \re{back to growth} were proved in Lemma \r{sl2 generators} parts  \re{sl2 generators_traces_geq} and  \re{sl2 generators_traces_leq}.   
Part  \re{many semi simples} was proved in Lemma \r{lem:many s.s_}.
}


Now let's see how Helfgott managed to reduce the Growth of 
$A^{[k]}$ to the Growth of $\Tr(A^{[k']})$ and
then to reduce the Growth of traces to the Growth of eigenvalues \cfr[\S 3 Proposition 3.3 and \S 4.4]{helf}.

\doTHM{[Helfgott] \l{helf reduction} 
There exist $k\in \N_+$ and $C \in \R_+$ such that \ti{for any} $\eps\in \R_+$ that following holds.
Let $\F$ be a finite field and let $A$ be a subset of generators of $\SL_2(\F)$.
Denote $A_1=A^{[k]}$ and $A_2=A^{[k^2]}$.
Suppose 
\doeq{\l{helf reduction_assumption_1} |A_2| < |A|^{1+\eps}.} 
Then we have,	
\doeq{\fr{C}|A|^{1/3} < |\Tr(A_1)| < C|A|^{1/3+C\eps} \l{helf reduction_bound_1a}}
and there exists a semi simple element $g\in A_1\cap G_s$ and $V:=\opC_{A_2}(g)$ with
\doeq{|V| > \fr{C}|\Tr(A_1)|^{1-C\eps}. \l{helf reduction_bound_1b}}
Moreover, if \doeq{\l{helf reduction_assumption_2}|A_2| < |A|^{1+\eps}\qtq[and]|\Tr(A_2)|<|\Tr(A_1)|^{1+\eps}}	
then there exists a semi simple element 
\doeqnn{g\in A_1\cap G_s\qtq[and]V:=\opC_{A_2}(g)}
\ST\  \re{helf reduction_bound_1b} holds and also
\doeq{\l{helf reduction_bound_2a}|\Tr(V^2)\ddot \Tr(V^2)|+|\Tr(V^2)+\Tr(V^2)| <C|\Tr(V^2)|^{1+C\eps}.}
}


\doPROOF{

In the following proof we will use always the notation $A_i=A^{[k^i]}$ but we will change few times the value of $k$ itself.
We will always increase its value so to fit to all the properties that we will need.
All the properties of subsets that we will use are ``monotone increasing'' in the sense that if $A^{[k]}$ has the property $\cP$ and 
$k\leq k'$ then $A^{[k']}\in \cP$ as well. Note that the hypothesis depend also in the value of $k$ (and they are ``monotone decreasing'' properties). 

By theorem \r{thm_helfgott reduction to traces_}
there exists $k_1\in \N_+$ \ST\ the following holds. 
From the assumption  \re{helf reduction_assumption_1} for $k\geq k_1$ we get,	
\doEQsplit{size of traces_geq}{
			|A|^{1/3} &\llBY{large traces}	|\Tr(A_1)|\\	
								&\llBY{back to growth} 	\frac{|A_1^{[k_1]}|^{4/3}}{|A_1|}\\
								&\leB \frac{|A_2|^{4/3}}{|A|}\\
								&\lnBY{helf reduction_assumption_1}	|A|^{(4/3)(1+\eps)-1}\\
								&\llB	|A|^{1/3+\Oeps}.	
	}
\sWADW  \re{helf reduction_bound_1a}. \good

We also get by theorem \r{thm_helfgott reduction to traces_} that,
\doeqnn{|A_1\cap G_s| \ggBY{many semi simples} |A|.} 			
Therefore by corollary \r{corr_large centrelizer_} 
there exists $k_2\in \N_+$ \ST\ the following holds.
If $k\geq \MAX{k_1,k_2}$ then there exists a semi simple $g\in B:=A_1\cap G_s$ with large centralizer:

\doEQsplit{eq:size of centrelizer}{
	|\opC_{A_2}(g)| 	&\ggBY{remark large centrelizer} |\Tr(B)|\frac{|A_1|}{|A_2|}\\		
								&\gnB (|\Tr(A_1)|-2)\frac{|A|}{|A_2|}\\
								&\ggBY{helf reduction_assumption_1} |\Tr(A_1)||A|^{-\eps}\\
								&\ggBY{size of traces_geq} |\Tr(A_1)|^{1-\Oeps}		
	} 
\sWAD\ with  \re{helf reduction_bound_1b}.\good

Let $g\in A_1\cap G_s$ be as in  \re{eq:size of centrelizer} with large centralizer
\doeq{\l{size of Tr_A_1_leq} |\Tr(A_1)| \llBY{eq:size of centrelizer} |\opC_{A_2}(g)|^{1+\Oeps}=|V|^{1+\Oeps}}
where $V:=\opC_{A_2}(g)$.
By conjugating $g$ with $u\in \SL_2(\F_{q^2})$ \ST\ $g^u$ is
diagonal we get that all $V^u$ are diagonal, since $g$ is regular semi simple. 
By corollary \r{lemma_escaping from a union of subspaces_} 
there exists $k_3\in \N_+$ \ST\ if $k\geq k_3$ then there 
exists $a\in A_1$ \ST\ $a^u$ has no zero entries.
Therefore if $k\geq \MAX{k_1,k_2,k_3}$ then 
\doeqnn{V^{[4]}V^{a[4]}:=V^{[4]}a^{-1}V^{[4]}a\subseteq A_2^{[10]} = A^{[10k^2]}\subseteq A^{[(10k)^2]}.} 
\TF\ if we take $k= 10\MAX{k_1,k_2,k_3}$ we get that
\doeqnn{V^{[4]}V^{a[4]} \subseteq A_2.}
Now suppose  \re{helf reduction_assumption_2} holds with $k= 10\MAX{k_1,k_2,k_3}$. 
Therefore we get,  
\doalign{
	|\Tr(V^{u[4]}(a^u)^{-1} V^{u[4]} a^u)|		&\eqB	|\Tr(V^{[4]}{V^{a[4]}})|\\
															&\leB |\Tr(A_2)|\\
															&\leBY{helf reduction_assumption_2}	|\Tr(A_1)|^{1+\eps}\\
															&\llBY{size of Tr_A_1_leq} |V|^{1+\Oeps}\\
															&\llB |\Tr(V)|^{1+\Oeps}
	}
By applying theorem \r{helfgott: reduction from trace function to sum product} with $V^u$ and $a^u$ we get, 
\doeqnn{|\Tr(V^2)\ddot \Tr(V^2)|+|\Tr(V^2)+\Tr(V^2)| \lnBY{prop.3.18.2.then} C|\Tr(V^2)|^{1+C\eps}}
\sWADW  \re{helf reduction_bound_2a}.\good
}

\sec{Main results}\l{sec:Main results}			

\ssec{From matrices to traces and back in finite fields}

\doPROP{\l{prop_escaping from impure almost-field}
There exists $C\in \R_+$ such that the following holds.
Let $\F$ be a finite field and $G=SL_2(\F)$ and let $\eps\in \R_+$ with $\eps< \fr{C}$.
Let $V\subseteq \SL_2(\F)$ be a subset of diagonal matrices of size $|V|>C$.

Suppose 
\doeq{\l{escaping from impure almost-field_if_1}
	\Tr(V)\tx{ is an impure }\eps\tx{-field}
}
and
\doeq{\l{escaping from impure almost-field_if_2}
	|\Tr(V^{[4]})|<|\Tr(V)|^{1+\eps}.
}

Then we have,
\doeq{\l{escaping from impure almost-field} 
	\Tr(V^{[4]})\tx{ is not }\eps\tx{-field}.
}
}

\doPROOF{
Set $N:=|\Tr(V)|$. 
\BTA\  \re{escaping from impure almost-field_if_1} there is some proper subfield $\E< \F$ and some $x\in\F^{\times}$ 
\ST\ we have 
\doeqnn{|\Tr(V)\SM x\E|<N^{\eps}\mqbq{and}|\E|<N^{1+\eps}.}
By the assumption \re{escaping from impure almost-field_if_1}, $\Tr(V)$ is an impure subfield so $|\Tr(V)\SM \E|>0$.
There are two cases to consider: either (1) $x\in \E$ or (2) $x\notin \E$.

Case (1): Suppose 
\doeqnn{x\in \E\qtq[and]0<|\Tr(V)\SM \E|<N^{\eps}}
and let $g\in V$ with $\Tr(g)\notin \E$.
Since $g(V\MID{\E^{\times}})\subseteq V^{[2]}$ we get by Lemma \r{lem:trace identity for matrices} that,
%
%
\doAlign{
	\l{bound_escape from E_}  |V^{[2]}\NMID{\E}| &\geB |(g(V\MID{\E^{\times}}))\NMID{\E}|\\
	\nn													&\geBY{eq:trace identity for matrices} \fr{2}|V\MID{\E^{\times}}|\\
	\nn  												&\geB \fr{2}(|V\MID{\E}|-2)\\
	\nn  												&\ggB |V\MID{\E}|\\
	\nn  												&\geB N-N^{\eps}\\
	\nn  												&\ggB N.  		 
	}	
\BTA\ $|\Tr(V^{[4]})|\leBY{escaping from impure almost-field_if_2} N^{1+\eps}$ 
so 
\doeq{\Tr(V^{[2]})\tx{ cannot be $\eps$-field}.\l{eq:not eps field}}
Indeed: the bound \re{bound_escape from E_} exclude the possibility of $\Tr(V^{[2]})$ to be $\E'$-field for $\E'=\E$ or any other coset $\E'=x\E$ of $\E$.  
Now for any other field $\E'\neq \E$ if $|\E'|\leq |\Tr(V^{[2]})|^{1+\eps}$ then
\doeqnn{|\E'|\leq N^{1+\Oeps}}
since $|\Tr(V^{[4]})|\leBY{escaping from impure almost-field_if_2} N^{1+\eps}$. 
Therefore the intersection of the field $\E$ with any coset $x'\E'$ is
\doeqnn{|\E\cap x'\E'| \leq |\E\cap \E'| \leq N^{\Oeps}.} 
So the intersection  is too small to contain $\Tr(V\MID{\E})$, since 
\doeqnn{|\Tr(V\MID{\E})|\geq \fr{2}|V\MID{\E}|\geq N-N^{\eps}\gg N.}
Therefore we are done with \re{eq:not eps field}.


Case (2): Suppose 
\doeqnn{\Tr(V)\subseteq x\E\qtq[with]|\E|\leq N^{1+\eps}\qtq[and]x\notin \E.} 
This case is proved in a similarly to Case (1):
By multiplying by some $g \in V\MID{x\E}$ we get by Lemma \r{lem:trace identity for matrices} 
that at least $\fr{2}|V\MID{x(\E^{\times})}|$
elements in $V^{[2]}$ have trace not in $x\E$.
Therefore, as was proved in \re{eq:not eps field} in Case (1), we find that 
$\Tr(V^{[2]})$ cannot be $\eps$-field. 

In both cases we get that $\Tr(V^{[2]})$ cannot be $\eps$-field
\sWADW  \re{escaping from impure almost-field}. 
}

\doPROP{\l{prop_escaping from a proper field}
There exist $C\in \R_+$ and $k\in \N_+$ with $k>C$ such that the following holds.
Let $\F$ be a finite field, $G=\SL_2(\F)$ and let $\eps\in \R_+$ with $\eps<\fr{C}$.

Let $\E<\F$ be a proper subfield, $A\subseteq \SL_2(\F)$ with $\gen{A}=G$. 
For $1\leq i\leq 2$ denote $A_i=A^{[k^i]}$ and suppose
\doeq{|A_3| < |A|^{1+\eps}.	\l{bound_small_A_3}}
Then there exists a semi simple element 
$$g\in A_1\cap G_s\qtq[and]V\subseteq \opC_{A_2}(g)$$ 
\ST
\doAlign{
		|\Tr(V)|	&>	\fr{C}|\Tr(A_2)|^{1-C\eps}	\l{escaping from a proper field_a} \\
		\Tr(V)		&\subseteq \F\SM \E		\l{escaping from a proper field_b} 
	}
}

\doPROOF{
In order to make the notations in the proof simpler, we will use the notation $A_i:=A^{[k^i]}$ and we will increase, during the proof, the value of $k$.

By Lemma \r{corr_escaping from all subfields_} there exists $k_1\in \N_+$ \ST\ for 
$k\geq k_1$ and $B:=A_1\NMID{\E}$ we have
\doeq{|B|=|A_1\NMID{\E}|\gg |A|. \l{bound_outside E_}}
Now let $g\in A_1\cap G_s$ be a semi simple element with 
\doeqnn{\Fix(g)=\set{x_1,x_2}\subseteq \P(\ol{\F}).} 
Suppose that for any $h\in A$ we have $\Fix(h)\subseteq \Fix(g)$.
Since $\gen{A}=G$ we have $\Fix(A)=\emptyset$ so we can find $h_1,h_2\in A$ \ST\ $\Fix(h_i)=\set{x_i}$ so
$\Fix(h_1h_2)\cap \Fix(g)=\emptyset$. 
\IAC\ there exist $h\in A^{[2]}$ \ST\ 
\doeqnn{\Fix(h)\SM \Fix(g)\neq \emptyset.}
Therefore by Lemma \r{large trace set} we get that if $k\geq \MAX{k_1,2}$ then 									 
\doEQsplit{bound_TrB_geq}{
	|B| &\llBY{remark large trace} |\Tr(B)||\Tr(gB)||\Tr(hB)|	\\	
			&\leB |\Tr(B)||\Tr(A_2)|^2.																	
}

Now by by theorem  \r{thm_helfgott reduction to traces_} there exists $k_2\in \N_+$ 
\ST\ if $k \ge \MAX{2,k_1,k_2}$ then we have,
\doEQsplit{bound_Tr_A_2_leq}{
	|\Tr(A_2)|	&\llBY{back to growth} \frac{|A_2^{[k_2]}|^{4/3}}{|A_2|}\\
							&\leB	\frac{|A_3|^{4/3}}{|A|}\\				
							&\llBY{bound_small_A_3} |A|^{1/3+\Oeps}.  	
	}
Therefore we conclude,
\doalign{	
	|\Tr(A_2)|^{3-\Oeps} &\llBY{bound_Tr_A_2_leq} |A|\\	
		&\llBY{bound_outside E_} |B|\\
		&\leBY{bound_TrB_geq} |\Tr(B)||\Tr(A_2)|^2\\		
		&\leB	|\Tr(A_2)|^3		
	}
%
Therefore we get 
\doAlign{
	|\Tr(B)|	&\ggB |\Tr(A_2)|^{1-\Oeps}	\l{bound_Tr_B}\\
	|\Tr(A_2)|	&\ggB |A|^{1/3}.	\l{bound_Tr_A_gg}		
	} 
Now suppose 
\doeqnn{k \geq \MAX{3,k_1,k_2}.}	
Therefore by corollary \r{corr_large centrelizer_} there exists $b\in B$ s.t.
\doEQsplit{bound_C}{
	|\opC_{B^{-1}B}(b)|	&\ggBY{remark large centrelizer}  \frac{|\Tr(B)||A_1|}{|A_1^{-1}A_1A_1|}\\
		&\geB \frac{|\Tr(B)||A|}{|A_3|}\\	
		&\geBY{bound_small_A_3} |\Tr(B)||A|^{-\eps}\\
		&\ggBY{bound_Tr_A_gg} |\Tr(B)||\Tr(A_2)|^{-\Oeps}\\
		&\ggBY{bound_Tr_B} |\Tr(A_2)|^{1-\Oeps}.
	} 
Let $b$ be as in  \re{bound_C} and set\footnote{I want to thank H.Helfgott for a very fruitful discussions related the following argument.}
\doalign{
	C		&:=	\opC_{B^{-1}B}(b)\\
	C'	&:=	C\NMID{0}\\
	C'' &:= C'\cup b C'\\
	V 	&:= C''\NMID{\E}.
	}
Note that $\Tr(b)\notin \E$ so $b$ is semi simple therefore we get that $\opC_G(b)$ are simultaneously diagonalizable 
therefore $|C'|\ge |C|-2$ and $|\Tr(V)|\ge \fr{2}|V|$.
Now by Lemma \r{lem:trace identity for matrices} we get that for any $c\in C'$ that,
\doeqnn{\mbq{either}\Tr(c)\notin \E \mqbq{or} \Tr(bc)\notin \E \mqbq{or} \Tr(bc^{-1})\notin \E.}
Altogether we get that
\doeqm{V \subseteq \opC_{A_2}(b),} since $V\subseteq A_1^{[3]} \subseteq A_2$, and 
\doalign{
	|\Tr(V)|	&\ggB |V|\\
		&\geBY{eq:trace identity for matrices} \fr{2}|C'|\\
		&\geB \fr{2}(|C|-2)\\
		&\ggB |C|\\
		&\ggBY{bound_C} |\Tr(A_2)|^{1-\Oeps}.	\QED
	}
}

\ssec{Conclusions}\l{ssec:main}

We extend the following key proposition of Helfgott \CITE[cf. ]{``Key Proposition'' in \S 1.2}{helf}.	
\doTHM{[Helfgott]	\label{thm_helfgott main_} 	
For any $\del\in \R_+$ there exist $\eps\in \R_+$ such that for any finite field $\F_p$ of prime order and any 
subset of generators $A$ of $G=\SL_2(\F_p)$ we have,
\doeqnn{|A|<|G|^{1-\del} \Rightarrow |A^{(3)}|>|A|^{1+\eps}.}	
Moreover, there exist absolute $k\in \N$ and $\del_0\in \R_+$ \ST
\doeqnn{|A|>|G|^{1-\del_0}\Rightarrow A^{[k]}=G.}	
}

The main result of this manuscript is the following extension of the theorem above.
\doTHM{[Theorem \r{thm:A^{(3)}:growth} from the Introduction]\l{main theorem}
\TES\ $\eps\in \R_+$ \ST\ \tFHS\ for any finite field $\F_q$. 
Let $G$ be the group $\SL_2(\F_q)$ and let $A$ be a generating set of $G$.
\TWH, 
\doeq{\l{eq:main theorem}|A^{(3)}|\ge \MIN{|A|^{1+\eps},|G|}.}
}


\doPROOF{
By theorem \r{BNP} there exists $C_0,\del_0\in \R_+$ \ST
\doeqnn{|A|\ge C_0|G|^{1-\del_0} >C_0q^{2\frac{2}{3}} \qu\RightarrowBY{BNP_bound_1a}\qu A^{(3)}=G.}
%
\TF\ if $|A|\ge C_0|G|^{1-\del_0}$ we are done with  \re{eq:main theorem}.
So we will assume from now 
\doeqnn{|A|\ll |G|^{1-\del_0}.} 

Let $3\leq k \in \N_+$, $\eps_0\in \R_+$ and $c_0\in \R_+$ with $c_0\leq 1$. 
By Lemma \r{lem:from growth to tripling} the following holds with $\eps'=\frac{\eps_0}{3k}$ and $c'=\frac{c_0}{2}$.
 For any group $G$ and any finite subset $A\subseteq G$ we have,
\doeqnn{|A^{[k]}|>c_0|A|^{1+\eps_0} \qu\RightarrowBY{eq:from growth to tripling:3}\qu |A^{(3)}|>c'|A|^{1+\eps'.}}
Now if $|A|^{\eps'/2}<\fr{c'}$ then $A$ is bounded but if $A$ is a subset of generators we get that
\doeqnn{|A^{(3)}|\geq |A|+2\geq |A|^{1+\eps''}}
for some $\eps''\in \R_+$. 
Therefore for any $\eps<\MIN{\eps'/2,\eps''}$ we get that,
\doeqnn{|A^{[k]}|>c_0|A|^{1+\eps_0} \RightarrowLB |A^{(3)}|>|A|^{1+\eps.}}


Therefore in order to prove  \re{eq:main theorem} it is enough to prove 
\doeqnn{|A^{[k]}|>c_0|A|^{1+\eps_0}}
for some absolute $3\leq k\in \N_+$ and $c_0,\eps_0\in \R_+$. 

We will use the notation 
\doeqnn{A_i:=A^{[k^i]}}
and we will prove that there exists $C\in \R_+$ and $i\in \N_+$ \ST\ the following holds.
There exists $k\in \N$ and $\eps\in \R_+$ with $k>C$ and $\eps<\fr{C}$ \ST\ we have
(provided $|A|\le C |G|^{1-\del_0}$) 
\doeqnn{|A_i|=|A^{[k^i]}|>\fr{C}|A|^{1+\eps^i}.} 

By Lemma \r{Tr(A^{[6]})} there exists $k_0\in \N_+$ 
\ST\ if $k>k_0$ then 
$\Tr(A_1)$ is not contained in any subfield i.e.,
\doeqnn{\gen{\Tr(A_1)}=\F_q.}
Set $\eps_1:=\fr{2}$. Note that if \doeqm{0<f<\eps_1}, 
then 
\doeqnn{1-f<\fr{1+f}<1-\fr{2}f<1-\Om(f)}
and similarly 
\doeqm{1+f<\fr{1-f}<1+2f<1+O(f).}

By theorem \r{helf reduction}  \re{helf reduction_bound_1a} 
there exists $k_1\in \N_+$ (and implicit $C_1>0$) \ST\ for any $\eps\in \R_+$ 
and \doeqm{k>\MAX{k_0,k_1}} we have either 
\doeqnn{|A_2|\geq |A|^{1+\eps}} (\sWAD) or 
\doeq{\l{eq:Tr(A_1):gg}|A|^{1/3} \ll |\Tr(A_1)| \ll |A|^{1/3+\Oeps}}
(explicitly: \doeqm{\fr{C_1}|A|^{1/3}<|\Tr(A_1)|< C_1|A|^{1/3+C_1\eps}}).
By applying again theorem \r{helf reduction}  \re{helf reduction_bound_1a} now for $A_1$, 
for any \doeqm{k>\MAX{k_0,k_1}} we have either 
\doeqnn{|A_3|\geq |A_1|^{1+\epss[2]}} (\sWAD) or 
\doeq{\l{eq:Tr(A_2):ll}|A_1|^{1/3} \ll |\Tr(A_2)| \ll |A_1|^{1/3+\Oepss}.}
Now if $|A_3|<|A|^{1+\eps^2}$ and in addition 
\doeqnn{|\Tr(A_1)|^{1+\eps}\leq |\Tr(A_2)|}
then both  \re{eq:Tr(A_1):gg} and  \re{eq:Tr(A_2):ll} hold and we get,
\doalign{		
	|A| &\llBY{eq:Tr(A_1):gg}	|\Tr(A_1)|^3\\
			&\lnB	|\Tr(A_2)|^{3(1-\fr{2}\eps)}\\
			&\llBY{eq:Tr(A_2):ll} |A_1|^{(1-\fr{2}\eps)(1+\Oepss)}\\
			&\lnB |A_1|^{1-\fr{2}\eps+\Oepss}\\
			&\leB |A_1|^{1-\Omeps}.
	}
\IOW\ $|A|^{1+\Omeps}\ll |A_1|$ \sWADW  \re{eq:main theorem}.

Restating the conclusion explicitly, there exist $C_2=k_2\in \N_+$ and $\eps_2=\fr{k_2}$ \ST\ 
if \doeqm{k>\MAX{k_i}} and \doeqm{\eps<\MIN{\eps_i}} then 
\doeqnn{|A|\ll |A_1|^{1-\fr{2}\eps+\Oepss}\leq |A_1|^{1-\fr{4}\eps}} 
i.e., \doeqm{\fr{C_2}|A|^{1+\fr{4}\eps}<|A_1|}.

Note that so far we have changed $k$ and $\eps$ \ti{independently} to get the required growth property  \re{eq:main theorem}.
We can summarize what have proved so far as follows.
There exist some $C\in \R_+$ \ST\ for any $\eps<\fr{C}$ and for any $k>C$ we have,
\doeqnn{|\Tr(A_2)|\geq |\Tr(A_1)|^{1+\eps} \Rightarrow |A_3|\geq \fr{C}|A|^{1+\eps^2}.}
We can restate this in the $\Om$-language as follows,
\doeq{\l{case:small growth of traces} |\Tr(A_2)|\gg |\Tr(A_1)|^{1+\Omeps} \Rightarrow |A_3|\gg |A|^{1+\Omepss}.}
Therefore in order to complete the proof, we can assume from now 
\doeq{|A_3|<|A|^{1+\eps^2} \mqbq{and} |\Tr(A_2)|<|\Tr(A_1)|^{1+\eps^2}. \l{case:growth:A_3 and Tr(A_2):leq}}
\IP\footnote{Note that there was nothing special in choosing $\eps^2$ above, and we can replace $\eps^2$ with any $f=o(\eps)$.} 
\doEQsplit{eq:Tr(A_2):leq}{
	|A_2|	&< |A_1|^{1+\eps^2}\\	
	|\Tr(A_2)|	&<|\Tr(A_1)|^{1+\eps^2} 
	}
so we can apply theorem \r{helf reduction}  \re{helf reduction_assumption_2}.
Therefore there exists a exists a semi simple element $g\in A_1\cap G_s$ and $V:=\opC_{A_2}(g)$ with
\doEQsplit{eq:V:geq}{
	|V| &> \fr{C_1}|\Tr(A_1)|^{1-C_1\eps^2}\\	
			&> \fr{C_1}|\Tr(A_1)|^{1-C_1\eps}
	}
and in addition 
\doEQsplit{eq:Tr(V) growth:leq}{
	|\Tr(V^2)\ddot \Tr(V^2)|+|\Tr(V^2)+\Tr(V^2)|&< 	C_1|\Tr(V^2)|^{1+C_1\eps^2}\\	
																							&<	C_1|\Tr(V^2)|^{1+C_1\eps}.
	}
Therefore we get,
\doEQsplit{eq:Tr(V):gg}{
	|\Tr(V)|	&\geB \fr{2}|V|\\
		&\ggBY{eq:V:geq} |\Tr(A_1)|^{1-\Oepss}\\
 		&\ggBY{eq:Tr(A_2):leq} |\Tr(A_2)|^{(1-\Oepss)(1-\eps^2)}\\
		&\ggB |\Tr(A_2)|^{1-\Oepss}	
	}
and also 
\doEQsplit{eq:Tr(V^2):geq}{
	|\Tr(V^2)|&\geB \fr{4}|V|\\
						&\ggBY{eq:Tr(V):gg} |\Tr(A_2)|^{1-\Oepss}\\
						&\ggBY{eq:Tr(A_2):ll} |A_1|^{1/3-\Oepss}\\
						&\geB k^{1/3-\Oepss}.	
}
Denote
\doalign{	
	V_1	&:=	V\\
	U_1	&:= \Tr(V_1^2)\\
	K_1	&:=	C_1|U_1|^{C_1\eps}.
}
\TF\ we get 
\doeq{|U_1\ddot U_1|+|U_1+U_1| \lnBY{eq:Tr(V) growth:leq} K_1|U_1| \l{eq:Tr(V) K growth:leq}}
and for some absolute $C_3\in \R_+$
\doeqnn{|U_1|\geBY{eq:Tr(V^2):geq} \fr{C_3}k^{1/3-C_3\eps}.}
Therefore by theorem \r{freiman thm} there exists (an absolute) $C\in\R_+$
 \ST\ either \doeqnn{|U_1|< CK_1^{C}} 
or for some subfield $\E_1\leq\F$ and $x_1\in\F^{\times}$ we have, 
\doeq{|U_1\SM x_1\E_1|\leq CK_1^{C}\qtq[and]|\E|\leq CK_1^C|U_1|	\l{eq:U_1_K field}.}
Now set $C_4:=2CC_1$ and $\eps_3=\fr{3CC_1C_3}$.
Since $CK_1^C=CC_1|U_1|^{CC_1\eps}$ we get that for any $\eps<\MIN{\eps_i}$ there exists $k>\MAX{k_i}$ \ST
\doeqnn{CK_1^C<|U_1|^{C_4\eps}<|U_1|.}
Therefore the alternative  \re{eq:U_1_K field} must hold and we get
\doeq{|U_1\SM x_1\E_1|\leq |U_1|^{C_4\eps}\qtq[and]|\E|\leq |U_1|^{1+C_4\eps}. \l{eq:U_1_eps field}}
\IP\ we get
\doalign{
	|A_2| &\ggBY{eq:Tr(A_2):ll} \Tr(A_2)^{3-\Oepss}\\
				&\ggB |U_1|^{3-\Oepss}\\
				&\geBY{eq:U_1_eps field} |\E|^{3-\Oeps}.
	}
Therefore for any $\del_0\in \R_+$ we can find $C_5\in \R_+$ large enough\footnote{\ST\ $C_5>\MAX{k_i}$ and $\fr{C_5}<\MIN{\eps_i}$.} 
and we can find $\eps<\fr{C_5}$ and $k>C_5$ \ST
\doeqnn{|A_2|>|\E|^{3-\del_0}.}
Therefore if $\E=\F$ then we are done by theorem \r{BNP} which guarantee bounded generation 
for large subsets of $\SL_2(\F)$.

\TF\ in order to complete the proof of  \re{eq:main theorem} we are left to treat the case that 
 for some \ti{proper} subfield $\E<\F$
\doeq{\Tr(V^2)\tx{ is }C_4\eps\tx{-field} \l{case:eps field}.}

Suppose first that,
\doeq{\l{case:impure field} \Tr(V^2) \mqbq{is an impure} \Oeps\mb{-field}.} 
By proposition \r{prop_escaping from impure almost-field}  \re{escaping from impure almost-field} we get that
\doeqnn{\Tr(V^{2[4]})\mqbq{is not}C_4\eps\mb{-field}.}
Denote $V_2 := V_1^{[4]}$ and $U_2 :=\Tr(V_2^2)$ and $K_2:=|U_2|^{C_4\eps}$ and 
\doeqnn{K_2' := (K_2/C)^{1/C} \gg |U_2|^{2C_1\eps}} 
Therefore by theorem \r{freiman thm} we get 
\doEQsplit{eq:TrV_2_2 K growth:geq}{
	|\Tr(V_2^2)\ddot \Tr(V_2^2)|+|\Tr(V_2^2)+\Tr(V_2^2)| &\gg K_2'|\Tr(V_2^2)|\\
							&\gg	|\Tr(V_2^2)|^{1+2C_1\eps} 
	}
Now by corollary \r{lemma_escaping from a union of subspaces_}
there exists $k_5 \in \N_+$ \ST\ the following hold for any $k>k_5$.
For any $w\in \GL_2(\ol{\F})$ there exists $g\in A_1$ \ST\ $g^w$ has no zero entries.
\IP\ we can apply this for the basis $v\in \GL_2(\ol{\F})$ for which $V^v$ are simultaneously diagonalizable.

Therefore by the bound  \re{eq:TrV_2_2 K growth:geq} we can apply theorem \r{helfgott: reduction from trace function to sum product}
and we get that for some absolute $C_6=k_6\in \N_+$ and for $k>\MAX{k_i}$ we have
\doalign{
	|\Tr(A_3)| &\ggB |\Tr(V_2^{[4]}V_2^{g[4]})|\\	
			&\ggBY{prop.3.18.2.then}	|\Tr(V_2)|^{1+\frac{2C_1}{C_6}\eps}\\
			&\ggB |\Tr(V_2)|^{1+\Omeps}\\
			&\geB	|\Tr(V_1^2)|^{1+\Omeps}\\	
			&\ggBY{eq:Tr(V^2):geq} |\Tr(A_1)|^{(1+\Omeps)(1-\Oepss)}\\
			&\ggB |\Tr(A_1)|^{1+\Omeps}.	
	}
Therefore we get 
\doeq{|\Tr(A_3)| \gg |\Tr(A_1)|^{1+\Omeps} \l{eq:Tr(A_3):geq:Tr(A_1)}}
and this imply that 
$$\tq[either]|\Tr(A_2)| \gg |\Tr(A_1)|^{1+\Omeps}\qtq[or]|\Tr(A_3)| \gg |\Tr(A_2)|^{1+\Omeps}.$$	
Therefore by what we have proved in  \re{case:small growth of traces} we get that 
$$\tq[either]|A_3| \gg |A|^{1+\Omepss}\qtq[or]|A_4| \gg |A_1|^{1+\Omepss}$$ 
\IOW\ by  \re{case:small growth of traces} we get
\doeq{|\Tr(A_3)| \gg |\Tr(A_1)|^{1+\Omeps} \Rightarrow |A_4| \gg |A|^{1+\Omepss}. \l{eq:Tr(A_3) growth}}
\TF\ if  \re{case:impure field} holds (the case of impure proper almost subfield) then by  \re{eq:Tr(A_3):geq:Tr(A_1)} we 
are done with the proof of  \re{eq:main theorem}. 

Therefore we are left to treat the second subcase of  \re{case:eps field} that
\doeq{Tr(V^2) \mqbq{is pure}\Oeps\mb{-field} \l{case:pure}}
for some proper subfield $\E<\F$.
Note that if $\Tr(V^{[4]})\nsubseteq \E$ then we are done with  \re{eq:main theorem} 
by a similar argument to  \re{case:impure field} which treated the case of impure $\Oeps$-field.

Therefore in order to complete  \re{eq:main theorem} we can assume in addition to  \re{case:pure} that
\eqSplit{traceV}{
	\Tr(V)	\subseteqB  &\Tr(V^{[4]}) 		\subseteqB \E \\		
	|\E|		\llB |&\Tr(V)|^{1+\Oeps}  \llBY{eq:Tr(V^2):geq} {|\Tr(A_1)|}^{1+\Oeps}.		
	}
Now suppose we can find $g\in A_1$ \ST\footnote{$v$ was a basis that $V^v$ were diagonal.} $\Prod(g^v)\notin \E$, then
by Lemma \r{lem:expansion of products} we get
\doeq{|\Tr(V)|^{2-\Oeps} \llBY{expansion of fields} |\Tr(VV^g)|\llB |\Tr(A_3)| \label{using expansion of fields}}
so $|\Tr(V)|\ll |\Tr(A_3)|^{\fr{2}+\Oeps}$. 
\OTOH $$|\Tr(V)|\ggBY{eq:Tr(V^2):geq} |\Tr(A_1)|^{1-\Oeps}$$
therefore 
\doalign{
	|\Tr(A_1)|	&\ll |\Tr(A_3)|^{(\fr{2}+\Oeps)(1+\Oeps)}\\
		&\ll |\Tr(A_3)|^{\fr{2}+\Oeps}\\
		&\ll |\Tr(A_3)|^{1-\Oeps}.  
	}
Therefore by  \re{eq:Tr(A_3) growth} we are done with the case  \re{case:pure}.

Therefore we are left to treat the case that  \re{case:pure} and  \re{traceV} hold and 
$\Prod(g^v)\in \E$ for any $g\in A_1$.
Therefore by fact \r{sfact:tr_g(x,y)} we get for any $g\in A_1$ that
\doeqnn{\Tr(VV^g)\subseteqBY{eq:tr_g(x,y)} \E.}
\IP\ by definition \r{def:[g,h]} we get,  
\doeq{\Tr([V,A_1]_{set})\subseteqBY{eq:[A,B]} \Tr(VV^{A_1})\subseteq \E. \l{traceVV^A}}

Therefore the only case that we are left to resolve, in order to complete  \re{eq:main theorem}, is:
\doEQsplit{case:pure:more}{
	\Tr(VV^{A_1}) &\subseteqBY{traceVV^A} \E\\	
	|\E| &\llB |\Tr(V)|^{1+\Oeps}\\		
	|\Tr(V)|	&\ggBY{eq:Tr(V):gg} |\Tr(A_2)|^{1-\Oeps}.	
} 

Now by proposition \r{prop_escaping from a proper field} there exists $C_7\in \R_+$ \ST\ the following holds
 with $k_7=C_7$ and $\eps_7=\fr{C_7}$.
Assume $k>\MAX{k_i}$ and $\eps<\MIN{\eps_i}$.
Since $|A_3|\lnBY{case:growth:A_3 and Tr(A_2):leq} |A|^{1+o(\eps)}<|A|^{1+\eps}$  
we get by proposition \r{prop_escaping from a proper field}
that there exists a semi simple element $h\in A_1\cap G_s$ and $U\subseteq \opC_{A_2}(h)$ with
\doalign{
	|\Tr(U)|	&\ggBY{escaping from a proper field_a} |\Tr(A_2)|^{1-\Oeps}\\	
	\Tr(U)		&\subseteqBY{escaping from a proper field_b} \F\SM \E.
}
Therefore there exists $u\in \SL_2(\F_{q^2})$ \ST\ $U^u$ are diagonal and
\doeq{\Tr(U) \cap \E = \emptyset.	\l{Tr(U):cap:E}} 

By repeating all the cases before  \re{case:pure:more}, 
but now with $\Tr(U)$ instead of $\Tr(V)$, we get that 
the only case that we need to treat, in order to complete the theorem, is: 
\doeqnn{\Tr(U)\mqbq{is}\Oeps\mb{-field}}
for some proper field $\E'<\F$ and also,
\doEQsplit{traceUU^A}{
	\Tr(UU^{A_1}) &\subseteq \E'\\	
	|\E'| &\ll |\Tr(U)|^{1+\Oeps}	\\	
	|\Tr(U)|	&\gg |\Tr(A_2)|^{1-\Oeps}. 
	} 

Let us check what we got so far. 
Denote $$N:=|\Tr(A_2)| \qtq[and] \E'':=\E \cap \E'.$$
By the construction of $U$ in  \re{Tr(U):cap:E} we get that $\E \neq \E'$ and
\doeqnn{N^{1-\Oeps} \ll \MIN{|\E|,|\E'|} < \MAX{|\E|,|\E'|} \ll N^{1+\Oeps}.} 
Therefore we get,
\doeq{|\E''|\ll N^{\Oeps}.\l{small intersection}}
\IP\ by combining  \re{case:pure:more} and  \re{traceUU^A} we get that 
\doeq{\Tr([U,V]_{set})\subseteq \E''. \l{traceUV}}

Now if $V$ and $U$ do not have a common fix point\footnote{i.e., $\Fix(g)\cap \Fix(h)=\Fix(U)\cap \Fix(V)=\emptyset$.}
, then by Lemma \r{lem:expansion of products} we get, 
$$|\Tr([U,V]_{set})|\ggBY{expansion of commutator} |\Tr(V)|\gg N^{1-\Oeps}.$$ 
Therefore by  \re{traceUV} and  \re{small intersection} we get a contradiction for $\eps$ small enough.

On the other hand, if $V$ and $U$ do have a common fix point, then denote
their eigen values $X$ and $Y$ \RESP.
So $\tr(X^{[4]})\subseteq \E$ and $\tr(Y^{[4]})\subseteq \E'$ and
 $X\subseteq \K$ and $Y\subseteq \K'$ where $\K$
and $\K'$ are the two quadratic extensions of $\E$ and $\E'$ \RESP.
Denote $\K''=\K \cap \K'$ so we get 
\doeqnn{|\K''|=|\E''|^2\ll N^{\Oeps}.} 
Therefore we get,
\doalign{
	|\Tr(A_3)|	&\ge |\Tr(A_2^{(2)})|\\
		&\ge |\Tr(UV)|\\
		&= |\tr(XY)|\\
		&\ge \fr{2}|XY|\\
		&\ge \fr{2}\frac{|X||Y|}{|\K''|}\\
		&\gg N^{2-\Oeps}\\
		&\gg |\Tr(A_2)|^{2-\Oeps}.
	} 
\TF\ by  \re{eq:Tr(A_3) growth} we are done  with the case  \re{case:pure}.
So the proof is complete.
}

%

%

\doCORO{[Corollary \r{corr:main:diameter} from the Introduction]\l{corr:diameter}
\TE\ $C,d\in \R_+$ \ST\ \tFHS\ for any finite field $\F_q$. 
Let $A$ be a subset of generators of $G=\SL_2(\F_q)$.
\TWH,	
\doAlign{
	&\diam^+(G,A) < C \log^d(|G|) \l{corr_helfgott diameter_bound_1a}\\
\addTX{and for any $\del\in \R_+$ we have,}
	&|A|>|G|^{\del} \Rightarrow \diam^+(G,A) < C\lr{\tfr{\del}}^d. \l{corr_helfgott diameter_bound_1b} 
	}	
}

\doPROOF{
First suppose $\gen{A}=G$ and $|A|>|G|^{\del}$. 
Then by theorem \r{main theorem} we get for some absolute $\eps_0\in \R_+$ that
\doEQ{eq:iteration}{ 
	|A^{(3)}| \geBY{eq:main theorem} \MIN{|A|^{1+\eps_0},|G|}\ge |G|^{\MIN{\del(1+\eps_0),1}}
	}
By iterating  \re{eq:iteration} we get for any $i\geq 0$ that
\doeqnn{A^{(3^i)}\ge |G|^{\MIN{\del(1+\eps_0)^i,1}}.}
Therefore by taking $i$ \ST\ $1-\del_0<\del(1+\eps_0)^i$ 
we get that $$A^{(3^{i+1})}=G.$$ 
Now if we take $i$ such that,
\doeqnn{(1-\del_0)\fr{\del}< (1+\eps_0)^i < (1+\eps_0)^{i+1} \leq 4(1-\del_0)\fr{\del}}
then for $d:=\log_{1+\eps_0}(3)$ and $C_1:= (4(1-\del_0))^d$ we get
\doeqnn{3^{i+1} = (1+\eps_0)^{(i+1)d} \leq C_1(\fr{\del})^d}
and $A^{(3^{i+1})}=G$ \sWAD\ with  \re{corr_helfgott diameter_bound_1b}.\good

Now for arbitrary subset of generators we have $|A|\geq 2$ therefore
\doeqnn{A^{[3^i]}\ge \MIN{2^{(1+\eps_0)^i},|G|}.}
Therefore if we take $\del=\fr{2}$ and $i$ such that 
\doeqnn{\del \log(|G|) \leq (1+\eps_0)^i \leq 2\del \log(|G|)} 
then for $C_2:= (2\del)^d$ we get,
\doeq{3^{i} = (1+\eps_0)^{i\ddot d} \le C_2 \log^d(|G|)	\l{corr_helfgott diameter_proof_g}}
and $|A^{(3^{i})}|\ge |G|^{\del}$.
Denote $m_1:=3^{i}\leq C_2 \log^d(|G|)$ \ST\  \re{corr_helfgott diameter_proof_g} holds.
Therefore if we apply the first argument to $A_1=A^{(m_1)}$ we get that there exists 
$m_2 \le C_1(\fr{\del})^d$ with $A_1^{[m_2]}=A^{[m_1m_2]}=G$. Therefore 
\doeqnn{m_1m_2\leq C_3 \log^d(|G|)}
where $C_3=C_1C_2(\fr{\del})^d$, \sWAD\ with  \re{corr_helfgott diameter_bound_1a}.\good
}

\sec{Further conjectures and questions} \l{sec:Further conjectures and questions} 

We include here some interesting questions that we encountered during this work.

\ssec{Trace generation}
In the proof of theorem \r{main theorem} we have seen in lines  \re{using expansion of fields}, 
that if one can prove that the set of products\footnote{see definition \r{def:Prod and Diag}.} $\Prod(A^{[k]})$ is not contained in any subfield, 
then by Lemma \r{lem:expansion of products}  \re{expansion of fields},
we could complete the proof with a much simpler argument. 

If it would be true, then we will get (in particular) that this property is preserved under
conjugation, since the generation property of the matrices is preserved. 
Is it true for any subset of generators?
\doquest{\l{quest:{Prod_A^{[k]}}}
Does there exists an absolute $k\in \N_+$ \ST\ for any finite field $\F_q$ and $A\subseteq \SL_2(\F_q)$ we have 
\doeqnn{\gen{A}=\SL_2(\F_q) \RightarrowB \gen{\Prod(A^{[k]})}=\F_q.}
}

By using the invariant argument of Lemma \r{trace generation} we have seen, 
\doeqnn{\gen{A}=\SL_2(\F_q) \RightarrowB \gen{\Tr(A^{[6]})}=\F_q.}
Can one use this property of $\Tr(A^{[k]})$ in order to prove question \r{quest:{Prod_A^{[k]}}}?


\ssec{Avoiding proper subfields}

We have seen in corollary \r{corr_escaping from all subfields_} how to escape from one subfield.
I.e., there are at least $c|A|$ elements in $A^{[k]}$ with trace outside this field, where $k$ and $c$ are absolute constants. 
Clearly we can always assume the subfield is a maximal subfield.
More precisely: 
for any proper $\E<\F_q$ we have\footnote{See definition \r{def_A_mid_} of $A\NMID{\E}$.},
\doeqnn{\gen{A}=\SL_2(\F_q) \Rightarrow |A^{[k]}\NMID{\E}|\gg |A|.}


Therefore if $\F_{p^n}$ has only one maximal subfield (i.e., $n$ is a prime power),
then we could complete the proof of theorem \r{main theorem} with a much simpler argument:

First, in the steps after equation  \re{case:growth:A_3 and Tr(A_2):leq} we would invoke proposition \r{prop_escaping from a proper field},
instead of theorem \r{helf reduction}  \re{helf reduction_assumption_2}.
By this proposition we would get that $\Tr(V^2)\subseteq \F\SM \E$
so in particular $\Tr(V^2)$ cannot be a pure $\Oeps$ field.
Under these terms, the proof will be much shorter since its second half  \re{case:pure},
which deals with the case of a pure $\Oeps$-fields, is no longer needed.

Now suppose $\F_{p^n}$ and $n=p_1^{n_1}\cdots p_m^{n_m}$. 
Therefore if we could escape from all the $m$ maximal subfields \ti{simultaneously} 
i.e., finding $c|A|$ elements in $A^{[k]}$ with traces which are primitive generators of $\F_q$,
we could simplify the proof of theorem \r{main theorem} as above.
Is it possible? Is it possible to avoid with the traces simultaneity a \ti{bounded} number of subfields?

\doquest{\l{quest:{Tr_A^{[k]}}}
Let $\F_q$ be a finite field and let $\E_1,\ldots,\E_m<\F_q$ be proper subfields of $\F$ and denote 
$$\mbb{W}=\CUP[1\leq i\leq m]{\E_i}.$$
Is is possible to find $k\in \N_+$ and $c\in \R_+$ (which may depend in $m$) \ST
\doeqnn{\gen{A}=\SL_2(\F_q) \RightarrowB  |A^{[k]}\NMID{\mbb{W}}|\geq c|A|?}
Is is possible to find an \ti{absolute} $k\in \N_+$ and $c\in \R_+$ as above?
}


\ssec{Growth of trace functions}

In the proof of theorem \r{main theorem} in equations  \re{traceV} and  \re{traceUU^A}, 
we got two subsets $V_i\subseteq \SL_2(\F)$ and two $v_i\in \SL_2(\ol{\F})$ 
\ST\ $V_i^{v_i}$ are diagonal matrices. 
Denote the eigen values of $V_i$ by $X_i$. 
I.e., $X_i$ are the diagonal entries of $V_i^{v_i}=D_{X_i}$. 
Denote $T_i:=\Tr(V_i)$ and note that $|T_i|\sim |X_i|$.
In the course of the proof of theorem \r{main theorem} we found two proper subfields $\E_i<\F$ \ST,
\doalign{
	N^{1-\Oeps}	\ll |T_1| &\le |T_2|=N\\
		T_i &\subseteq \E_i\\
		X_i &\subseteq \K_i\\
	|\E_i| &\ll N^{1+\Oeps}\\
	\tq[and]\E_1 &\neq \E_2
}
where $\K_i$ are the quadratic extension of $\E_i$ \RESP.
\TF,
\doalign{
	|T_1\cap T_2|	&\le |\E_1\cap \E_2|	\ll N^{\Oeps}\\
	|X_1\cap X_2|	&\le |\K_1\cap \K_2|	\ll N^{\Oeps}
}
and so we have, 
\doeqnn{|X_1 X_2|,|X_1+X_2|,|X_1 T_2|,|X_1+T_2|\ldots \gg N^{2-\Oeps}.}
Let $g=\dopmatrix{a & b\\ c& d}$ and denote%
\footnote{%
See definitions \r{def:Tr_g}, \r{def:tr_g} and facts \r{fact:trace identity for scalars}, \r{sfact:tr_g(x,y)} in \Sr{ssec:Expansion functions in fields}.},
\doalign{
		D_x				&:= \dobdiag{x}\\
		\tr(x)		&:= \Tr(D_x) = x+x^{-1}\\
		\tr_g(x,y) &:= \Tr(D_x (D_y)^g)\\
							&= ad\ddot \tr(xy)-bc \ddot \tr(x/y).
}
Now set $u_i:=v_i^{-1}$ and $h:=u_2gv_1$ and define,
\doalign{
	\Tr_g(g_1,g_2)	&:= \Tr((g_1)(g_2)^g)\\
								&= \Tr\lr{\lr{D_{x_1}}^{u_1}\lr{D_{x_2}}^{u_2g}}\\
								&= \Tr\lr{\dosbdiag{x_1}\dosbdiag{x_2}^{u_2gv_1}}\\
								&= \tr_{h}(x_1,x_2).
}
Now define $T_g:X_1\times X_2\to \F$ by 
\doeqnn{T_g(x_1,x_2):=\tr_{u_2gv_1}(x_1,x_2).} 
Therefore if we denote $t=\Prod(h)$ then, 
\doalign{
	\op{Im}(T_g)	&= T_g(X_1,X_2)\\
			&= \Tr(V_1 V_2^g)\\
			&= \set{t\ddot \tr(x_1x_2)+(1-t)\ddot \tr(x_1/x_2):x_i\in X_i}
	}
Now in the proof of theorem \r{main theorem} we have seen that Trace Growth imply Growth of Matrices.
\TF\ if one can prove that 
\doeqnn{|\op{Im}(T_g)|=|\Tr(V_1 V_2^g)|\gg  N^{1+\Omeps}} 
then we could simplify the arguments in the proof of theorem \r{main theorem}.

\doquest{\l{quest:Im(T_g)}
Can one prove for some absolute $\delta>0$ that
\doeqnn{|\doIm{T_g}|\geq N^{1+\delta-\Oeps}}
with $T_g$ and $N$ as defined above.
}

%

Now denote 
\doeqnn{N_t(c):=|\set{(x,y):T_g(x,y)=c}|}
the number of solution for $T_g(x,y)=c$ where
$$T_g(x,y)=t\ddot \tr(xy)+(1-t)\ddot \tr(x/y)$$ as was defined above. 

What is the best upper bound $N_t(c)$ for a general $t$\ ?

\doquest{
Can one prove for some absolute $\delta>0$ we have,
\doeqnn{t\in \K_1\K_2\SM \set{0,1} \RightarrowLB |N_t(c)|\ll N^{1-\delta+\Oeps}.}
}

If the answer is affirmative then 
\doeqnn{|\opIm(T_g)|\gg \frac{|X_1||X_2|}{\mult(T_g)}\gg \frac{N^{2-\Oeps}}{N^{1-\delta+\Oeps}}\gg N^{1+\delta-\Oeps}} 
so question \r{quest:Im(T_g)} is resolved also.
Note that we have seen in Lemma \r{lem:expansion of products},
\doeqnn{t\notin \K_1\K_2 \RightarrowBY{expansion of fields} \mult(T_g)=1 \Rightarrow |\opIm(T_g)|=|X_1||X_2|\gg N^{2-\Oeps}.} 
And since $|X_1X_2|,|X_1X_2^{-1}|\gg N^{2-\Oeps}$ we also have
\doeqnn{t\in\set{0,1} \Rightarrow |\opIm(T_g)|\gg N^{2-\Oeps}.}

\end{document}